\documentclass[12pt,leqno,openany]{book}

\usepackage{makeidx}
\makeindex

\def\co{\mathop{\rm co}}

\def\dist{\mathop{\rm dist}}
\def\supp{\mathop{\rm supp}}

\newtheorem{theorem}{Theorem}
\newtheorem{lemma}[theorem]{Lemma}
\newtheorem{proposition}[theorem]{Proposition}
\newtheorem{sublemma}[theorem]{Sublemma}
\newtheorem{definition}[theorem]{Definition}
\newtheorem{corollary}[theorem]{Corollary}
\newtheorem{problem}[theorem]{Problem}
\newtheorem{remark}[theorem]{Remark}
\newtheorem{claim}[theorem]{Claim}
\newtheorem{assumptions}[theorem]{Assumptions}
\newtheorem{examples}[theorem]{Examples}
\newtheorem{question}[theorem]{Question}
\newtheorem{sassumptions}[theorem]{Standing Assumptions}
\newtheorem{sassumption}[theorem]{Standing Assumption}
\newtheorem{conjecture}[theorem]{Conjecture}

\newcommand{\begintheorem}{\addtocounter{equation}{1}\begin{theorem}}
\newcommand{\beginlemma}{\addtocounter{equation}{1}\begin{lemma}}
\newcommand{\beginproposition}{\addtocounter{equation}{1}\begin{proposition}}
\newcommand{\beginsublemma}{\addtocounter{equation}{1}\begin{sublemma}}
\newcommand{\begindefinition}{\addtocounter{equation}{1}\begin{definition}}
\newcommand{\begincorollary}{\addtocounter{equation}{1}\begin{corollary}}
\newcommand{\beginproblem}{\addtocounter{equation}{1}\begin{problem}}
\newcommand{\beginremark}{\addtocounter{equation}{1}\begin{remark}}
\newcommand{\beginclaim}{\addtocounter{equation}{1}\begin{claim}}
\newcommand{\beginassumptions}{\addtocounter{equation}{1}\begin{assumptions}}
\newcommand{\beginexamples}{\addtocounter{equation}{1}\begin{examples}}
\newcommand{\beginquestion}{\addtocounter{equation}{1}\begin{question}}
\newcommand{\beginsassumptions}{\addtocounter{equation}{1}\begin{sassumptions}}
\newcommand{\beginsassumption}{\addtocounter{equation}{1}\begin{sassumption}}
\newcommand{\beginconjecture}{\addtocounter{equation}{1}\begin{conjecture}}

\begin{document}

\frontmatter

\title{Notes on Topological Vector Spaces}

\author{Stephen Semmes \\
	Department of Mathematics \\
	Rice University}

\date{}

\maketitle

\chapter{Preface}

	In the notion of a topological vector space, there is a
very nice interplay between the algebraic structure of a vector
space and a topology on the space, basically so that the vector
space operations are continuous mappings.  There are also
plenty of examples, involving spaces of functions on various domains,
perhaps with additional properties, and so on.  Here we shall try
to give an introduction to these topics, without getting too much
into the theory, which is treated more thoroughly in the books
in the bibliography.  Knowledge of linear algebra is assumed, and
some familiarity with elementary topology would be useful, with
some of the relevant material being reviewed in the first chapter.

\tableofcontents

\mainmatter

\chapter{Some background information}
\label{chapter on background information}

\section{Real and complex numbers}
\label{section on real and complex numbers}

	The real numbers\index{real numbers} will be denoted ${\bf
R}$, and the complex numbers will be denoted ${\bf C}$.  Thus every
complex number can be written as $a + b \, i$, where $a$, $b$ are real
numebers and $i^2 = -1$.  If $z$ is a complex number with $z = a + b
\, i$, where $a, b \in {\bf R}$, then $a$, $b$ are called the real and
imaginary parts of $z$.  The integers\index{integers} are denoted
${\bf Z}$, and the positive integers are denoted ${\bf Z}_+$.  The
positive integers together with $0$ are called the nonnegative
integers.

	Of course the real and complex numbers are equipped with
the arithmetic operations of addition, subtraction, multiplication,
and division, satisfying the usual properties.  Let us mention
also the well-known facts that every polynomial on ${\bf R}$
with real coefficients and having odd degree takes the value $0$
somewhere on ${\bf R}$, and that every polynomial on ${\bf C}$
which is not a constant takes the value $0$ somewhere on ${\bf C}$.
By applying the second statement repeatedly, it follows that every
polynomial on ${\bf C}$ can be expressed as a constant times a 
product of linear factors.

	On the real numbers there is also the usual ordering with its
standard properties, such as the sum and product of positive real
numbers being positive real numbers.  If $a$, $b$ are real numbers
such that $a < b$, then we can define the open
interval\index{intervals of real numbers} $(a, b)$, the half-open,
half-closed intervals $[a, b)$, and $(a, b]$, and the closed interval
$[a, b]$, where the latter is also defined when $a = b$.  Namely, the
open interval $(a, b)$ consists of the real numbers $x$ such that $a <
x < b$, the interval $[a, b)$ consists of the real numbers $x$ such
that $a \le x < b$, the interval $(a, b]$ consists of the real numbers
$x$ such that $a < x \le b$, and $[a, b]$ consists of the real numbers
$x$ such that $a \le x \le b$.  We also allow $a$ or $b$ to be $\pm
\infty$ when appropriate, so that for example $(-\infty, \infty)$ is
the real line, and $[0, \infty)$ is the set of nonnegative real
numbers.

	If $E$ is a set of real numbers, then a real number $a$ is
said to be a \emph{lower bound} for $E$\index{lower bound for a set of
real numbers} if $a \le x$ for all $x \in E$.  Similarly, a real
number $b$ is said to be an \emph{upper bound} for $E$\index{upper
bound for a set of real numbers} if $x \le b$ for all $x \in E$.  A
real number $c$ is said to be the \emph{greatest lower
bound}\index{greatest lower bound for a set of real numbers} or
\emph{infimum}\index{infimum of a set of real numbers} of $E$ if $c$
is a lower bound of $E$, and if $a \le c$ whenever $a \in {\bf R}$ is
a lower bound for $E$.  It is easy to see from the definition that $c$
is unique if it exists.  A real number $d$ is said to be the
\emph{least upper bound}\index{least upper bound of a set of real
numbers} or \emph{supremum}\index{supremum of a set of real numbers}
of $E$ if $d$ is an upper bound for $E$ and $b \le d$ whenever $b \in
{\bf R}$ is an upper bound for $E$.  Again, it is easy to see that $d$
is unique if it exists.  A completeness property of the real numbers
is that every nonempty set $E$ of real numbers which has an upper
bound also has a least upper bound.  This implies that every nonempty
set $F$ of real numbers which has a lower bound also has a greatest
lower bound.  Indeed, the existence of the infimum of $F$ can be
derived from the existence of the supremum of $E = - F$, which
consists of the real numbers $-x$ for $x \in F$.  Alternatively, the
infimum of $F$ can be obtained as the supremum of the set of lower
bounds of $F$.

	If $x$ is a real number, then the \emph{absolute
value}\index{absolute value of a real number} is denoted $|x|$
and defined by $|x| = x$ if $x \ge 0$, $|x| = -x$ if $x \le 0$.
Thus the absolute value of a real number is always a nonnegative real
number, and one can check that
\begin{equation}
	|x + y| \le |x| + |y|, \qquad |x y| = |x| \, |y|
\end{equation}
for all $x, y \in {\bf R}$.

	If $z$ is a complex number, $z = x + y \, i$ with $x, y \in
{\bf R}$, then the \emph{complex conjugate}\index{complex conjugate of
a complex number} of $z$ is denoted $\overline{z}$ and defined by
\begin{equation}
	\overline{z} = x - y \, i.
\end{equation}
It is easy to see that
\begin{equation}
	\overline{z + w} = \overline{z} + \overline{w}, \qquad
		\overline{z \, w} = \overline{z} \, \overline{w}
\end{equation}
for all $z, w \in {\bf C}$, and that the complex conjugate of the complex
conjugate of $z$ is equal to $z$.  Also,
\begin{equation}
	z \, \overline{z} = x^2 + y^2,
\end{equation}
and in particular $z \, \overline{z}$ is always a nonnegative real
number.

	The \emph{modulus}\index{modulus of a complex number} or
\emph{absolute value}\index{absolute value of a complex number} of a
complex number $z$ is defined to be $|z| = \sqrt{z \, \overline{z}}$.
It is well-known that
\begin{equation}
	|z + w| \le |z| + |w|
\end{equation}
for all $z, w \in {\bf C}$, and we also have that
\begin{equation}
	|z \, w| = |z| \, |w|
\end{equation}
for all $z, w \in {\bf C}$.

	As a special case, if $\theta$ is a complex number such that
$|\theta| = 1$, then
\begin{equation}
	|\theta \, z| = |z|
\end{equation}
for all $z \in {\bf C}$.  In effect, multiplication by $\theta$
represents a rotation on ${\bf C}$.  

	Let $\{z_j\}_{j=1}^\infty$ be a sequence of complex numbers,
and let $z$ be another complex number.  We say that $\{z_j\}_{j=1}^\infty$
converges to $z$ if for every $\epsilon > 0$ there is a positive
integer $N$ such that
\begin{equation}
	|z_j - z| < \epsilon \ \hbox{ for all } j \ge N.
\end{equation}
In this event we write
\begin{equation}
	\lim_{j \to \infty} z_j = z,
\end{equation}
and call $z$ the limit of the sequence $\{z_j\}_{j=1}^\infty$.  It is
easy to see that the limit of a convergent sequence is unique.  

	It is easy to verify that a sequence $\{z_j\}_{j=1}^\infty$ of
complex numbers converges to a complex number $z$ if and only if the
sequence of real parts of the $z_j$'s converges to the real part of
$z$ and the sequence of imaginary parts of the $z_j$'s converges to
the imaginary part of $z$.  Similarly, $\{z_j\}_{j=1}^\infty$
converges to $z$ if and only if $\{\overline{z}_j\}_{j=1}^\infty$
converges to $\overline{Z}$ of $z$.  If $\{z_j\}_{j=1}^\infty$ is a
sequence of complex numbers which converges to the complex number $z$,
then $\{|z_j|\}_{j=1}^\infty$ converges to $|z|$.

	Suppose that $\{z_j\}_{j=1}^\infty$, $\{w_j\}_{j=1}^\infty$
are sequences of complex numbers which converge to $z, w \in {\bf C}$,
respectively.  It is well known that the sequences $\{z_j +
w_j\}_{j=1}^\infty$ and $\{z_j \, w_j\}_{j=1}^\infty$ converge
to $z + w$ and $z \, w$, respectively.  If $\{z_j\}_{j=1}^\infty$
is a sequence of nonzero complex numbers which converges to the
nonzero complex number $z$, then $\{1/z_j\}_{j=1}^\infty$ converges
to $1/z$.

	A sequence of complex numbers $\{z_j\}_{j=1}^\infty$ is said
to be a \emph{Cauchy sequence} if for every $\epsilon > 0$
there is a positive integer $N$ such that
\begin{equation}
	|z_j - z_l| < \epsilon \ \hbox{ for all } j, l \ge N.
\end{equation}
It is easy to check that if $\{z_j\}_{j=1}^\infty$ is a convergent
sequence of complex numbers, then it is also a Cauchy sequence.
A completeness property of the complex numbers states that
every Cauchy sequence in ${\bf C}$ converges in ${\bf C}$.

	Suppose that $p$ is a positive real number.  If $p \ge 1$, then
the function $t^p$ on $[0, \infty)$ is \emph{convex}, which means that
\begin{equation}
	(\lambda \, a + (1 - \lambda) \, b)^p 
		\le \lambda \, a^p + (1 - \lambda) \, b^p
\end{equation}
for all nonnegative real numbers $a$, $b$, and $\lambda$ with $\lambda
\le 1$.  This is a well-known fact from calculus.

	If $p \le 1$, then another useful inequality states that
\begin{equation}
	(a + b)^p \le a^p + b^p.
\end{equation}
This can also be analyzed in terms of calculus, and one can give
a more direct derivation as well.

\section{Vector spaces}
\label{section on vector spaces}

	Let $V$ be a vector space.  In this monograph we make the
standing assumption that all vector spaces use either the real or the
complex numbers as scalars, and we say ``real vector spaces'' and
``complex vector spaces'' to specify whether real or complex numbers
are being used.  To say that $V$ is a vector space means that $V$ is a
nonempty set with a distinguished element called $0$ and operations of
addition and scalar multiplication which satisfy the usual properties.
Recall that a subset $L$ of $V$ is said to be a \emph{linear
subspace}\index{linear subspace of a vector space} if $0 \in L$,
$v + w \in L$ whenever $v, w \in L$, and $\zeta \, v \in L$ whenever
$\zeta$ is a scalar, which is to say a real or complex number,
as appropriate, and $v \in L$.  Thus $L$ is then a vector space
too, with the same choice of scalars, and using the restriction
of the vector space operations from $V$ to $L$.

	If $V_1$, $V_2$ are two vector spaces, both real or both
complex, then a mapping $f : V_1 \to V_2$ is said to be
\emph{linear}\index{linear mappings between vector spaces} if
\begin{equation}
	f(v + w) = f(v) + f(w)
\end{equation}
for all $v, w \in V_1$, and
\begin{equation}
	f(\alpha \, v) = \alpha \, f(v)
\end{equation}
for all scalars $\alpha$ and $v \in V_1$.  If $f$ is a one-to-one
mapping from $V_1$ onto $V_2$, so that the inverse mapping $f^{-1} :
V_2 \to V_1$ is defined, then $f^{-1}$ is a linear mapping from $V_2$
to $V_1$ if $f$ is a linear mapping from $V_1$ to $V_2$.  In this
case we say that $f$ defines an \emph{isomorphism}\index{isomorphism
between vector spaces} from $V_1$ onto $V_2$ as vector spaces, and
that $V_1$, $V_2$ are isomorphic vector spaces.

	The space of linear mappings from $V_1$ to $V_2$ is denoted
$\mathcal{L}(V_1, V_2)$.  It is easy to see that the sum of two
elements of $\mathcal{L}(V_1, V_2)$ defines an element of
$\mathcal{L}(V_1, V_2)$, and that the product of a scalar and an
element of $\mathcal{L}(V_1, V_2)$ defines an element of
$\mathcal{L}(V_1, V_2)$.  Thus $\mathcal{L}(V_1, V_2)$ is a vector
space in a natural way, with the same scalars as for $V_1$, $V_2$.
For a single vector space $V$ we may write $\mathcal{L}(V)$ instead of
$\mathcal{L}(V, V)$.

	If $V$ is a real vector space and $f$ is a linear mapping from
$V$ to ${\bf R}$, or if $V$ is a complex vector space and $f$ is a
linear mapping from $V$ to ${\bf C}$, then $f$ is called a
\emph{linear functional}\index{linear functional on a vector space} on
$V$.  The space of linear functionals on $V$ is called the
\emph{dual}\index{dual of a vector space} of $V$ and is denoted $V'$.
Thus $V' = \mathcal{L}(V, {\bf R})$ and $V'$ is a real vector space
when $V$ is a real vector space, and $V' = \mathcal{L}(V, {\bf C})$
and $V'$ is a complex vector space when $V$ is a complex vector space.

	If $V_1$, $V_2$, and $V_3$ are vector spaces with the same
scalars, and if $f_1 : V_1 \to V_2$ and $f_2 : V_2 \to V_3$ are linear
mappings, then the composition $f_2 \circ f_1 : V_1 \to V_3$, defined
by $(f_2 \circ f_1)(v) = f_2(f_1(v))$ for $v \in V_1$, defines a
linear mapping from $V_1$ to $V_3$.  For a single vector space
$\mathcal{L}(V)$, we have that the composition of two elements of
$\mathcal{L}(V)$ is also an element of $\mathcal{L}(V)$.

	If $n$ is a positive integer, then ${\bf R}^n$ and ${\bf
C}^n$, which consist of $n$-tuples of real and complex numbers,
respectively, are real and complex vector spaces with respect to
coordinatewise addition and scalar multiplication.  If $V$ is a real
vector space with positive finite dimension $n$, then $V$ is
isomorphic to ${\bf R}^n$, and if $V$ is a complex vector space of
positive finite dimension $n$, then $V$ is isomorphic to ${\bf C}^n$.

	More generally, if $n$ is a positive integer greater than or
equal to $2$, and $V_1, V_2, \ldots, V_n$ are vector spaces, all real
or all complex, then the Cartesian product 
\begin{equation}
	V_1 \times V_2 \times \cdots \times V_n
\end{equation}
defines a vector space in a natural way, with the same choice of
scalars.  Specifically, this Cartesian product consists of $n$-tuples
$(v_1, v_2, \ldots, v_n)$, with $v_j \in V_j$ for $j = 1, 2, \ldots, n$.
If $(v_1, v_2, \ldots, v_n)$ and $(w_1, w_2, \ldots, w_n)$ are
two elements of the Cartesian product, then their sum is defined
coordinatewise, which is to say that the sum is equal to
\begin{equation}
	(v_1 + w_1, v_2 + w_2, \ldots, v_n + w_n).
\end{equation}
Similarly, if $(v_1, v_2, \ldots, v_n)$ is an element of the
Cartesian product and $\zeta$ is a scalar, then the scalar product
of $\zeta$ with $(v_1, v_2, \ldots, v_n)$ is defined coordinatewise
and is equal to
\begin{equation}
	(\zeta \, v_1, \zeta \, v_2, \ldots, \zeta \, v_n).
\end{equation}
Actually, the vector space that results in this manner is called the
\emph{direct sum}\index{direct sum of vector spaces} of $V_1, V_2,
\ldots, V_n$.  If each $V_j$ has finite dimension, then the direct
sum also has finite dimension, and the dimension of the direct
sum is equal to the sum of the dimensions of the $V_j$'s.

	If $V_1$ and $V_2$ are vector spaces, both real or both
complex, and if $f$ is a mapping from $V_1$ to $V_2$, then we can
associate to $f$ its graph in the direct sum of $V_1$ and $V_2$, which
is the set of ordered pairs $(v, f(v))$ with $v \in V_1$.  One can
verify that $f$ is a linear mapping from $V_1$ to $V_2$ if and only if
the graph of $f$ is a linear subspace of the direct sum of $V_1$ and
$V_2$.

	Let $n$ be a positive integer, $n \ge 2$, and suppose that
$V_1, V_2, \ldots, V_n$ and $W_1, W_2, \ldots, W_n$ are vector
spaces, all with the same choice of scalars.  Also assume that
for each $j = 1, 2, \ldots, n$ we have a linear mapping
$f_j : V_j \to W_j$.  Then we get a linear mapping $F$ from
the direct sum of the $V_j$'s to the direct sum of the $W_j$'s
in a natural way, following the rule
\begin{equation}
	(v_1, v_2, \ldots, v_n) \mapsto (f(v_1), f(v_2), \ldots, f(v_n)).
\end{equation}

	A complex vector space can be viewed as a real vector space
by ``restriction of scalars'', which amounts to just using real
numbers as scalars even if complex multiplication by $i$ is also
defined.  If $V$ is a complex vector space of finite dimension $n$,
then when we view it as a real vector space in this manner, it has
dimension $2n$.  If $V$ is a real vector space, then we can
``complexify'' it by taking the Cartesian product $V \times V$,
initially as a real vector space, and then extending the scalar
multiplication to complex numbers by taking
\begin{equation}
	i \, (v_1, v_2) = (-v_2, v_1).
\end{equation}
In other words, if $\widehat{V}$ denotes the complexification of $V$,
then we can think of $\widehat{V}$ as consisting of vectors of the
form $v_1 + v_2 \, i$, where $v_1, v_2 \in V$, and where scalar
multiplication by $i$ is defined in the obvious manner.  If $V$ is a
real vector space with finite dimension $n$, then the complexification
$\widehat{V}$ of $V$ is a complex vector space with dimension $n$.

	Here is another point of view with similar ideas.  Let $W$
be a real vector space, and let $J$ be a linear mapping from $W$
to itself whose square as a linear mapping is equal to minus the
identity mapping on $W$, i.e.,
\begin{equation}
	J(J(w)) = -w
\end{equation}
for all $w \in W$.  Then we can view $W$ as a complex vector space,
where $i \, w$ is defined to be $J(w)$ for all $w \in W$.  If $W$ has
finite dimension $n$ as a real vector space, then $n$ is necessarily
even in this case, and $W$ has dimension $n/2$ as a complex vector
space.

	Let $V$ be a vector space, real or complex, and let $E$ be a
subset of $V$.  We say that $E$ is \emph{convex}\index{convex subset
of a vector space} if $\lambda \, v + (1 - \lambda) w \in E$ whenever
$v, w \in E$ and $\lambda \in [0, 1]$.  If $A$ and $B$ are nonempty
subsets of $V$, then we write $A + B$ for the subset of $V$ consisting
of vectors of the form $v + w$, where $v \in A$ and $w \in B$.  If $A$
is a nonempty subset of $V$ and $\zeta$ is a scalar, i.e., a real or
complex number, as appropriate, then we write $\zeta \, A$ for the
subset of $V$ consisting of vectors of the form $\zeta \, v$, where $v
\in A$.  The statement that $E$ is convex can be rewritten as saying
that
\begin{equation}
	\lambda \, E + (1 - \lambda) \, E \subseteq E
\end{equation}
for all $\lambda \in [0, 1]$.

	If $A$ and $B$ are nonempty subsets of $V$ which are convex,
then $A + B$ is convex too.  If $A$ is a convex nonempty subset of $V$
and $\zeta$ is a scalar, then $\zeta \, A$ is a convex subset of $V$.

	If $E$ is a nonempty convex subset of $V$ and $f(v)$ is a
real-valued function on $E$, then we say that $f(v)$ is
\emph{convex}\index{convex function on a convex subset of a vector
space} if
\begin{equation}
	f(\lambda \, v + (1 - \lambda) \, w) 
		\le \lambda \, f(v) + (1 - \lambda) \, f(w)
\end{equation}
for all $v, w \in E$ and $\lambda \in [0, 1]$.  The sum of two convex
functions is convex, and the product of a convex function by a
nonnegative real number is convex.

	The property of convexity of $f(v)$ on $E$ can be
characterized geometrically as follows.  Inside the Cartesian
product $V \times {\bf R}$, consider the set
\begin{equation}
	\{(v, t) \in V \times {\bf R} : v \in E, \ t \ge f(v)\}.
\end{equation}
The function $f(v)$ on $E$ is convex if and only if this is
a convex subset of $V \times {\bf R}$.

	If $A$ is a nonempty subset of $V$, then we say that $A$ is
\emph{symmetric}\index{symmetric subset of a vector space} if $-A =
A$.  A subset $A$ of $V$ which contains $0$ is said to be
\emph{starlike around $0$}\index{starlike subset of a vector space
around $0$} if $t \, A \subseteq A$ for each real number $t$ such that
$0 \le t \le 1$.  If $V$ is a complex vector space and $A$ is a
nonempty subset of $V$, then $A$ is said to be
\emph{circular}\index{circular subset of a complex vector space} if
$\theta \, A = A$ for all complex numbers $\theta$ such that $|\theta|
= 1$.  Thus a nonempty subset $A$ of $V$ is starlike about $0$ and
circular if and only if $\zeta \, A \subseteq A$ for all complex
numbers $\zeta$ such that $|\zeta| \le 1$.

	If $V$ is a vector space and $A$ is a nonempty subset of $V$,
then the \emph{convex hull}\index{convex hull of a subset of a vector
space} of $A$ is denoted $\co(A)$ and is the subset of $V$ consisting
of all convex combinations of elements of $A$.  More precisely, a
vector $w \in V$ lies in $\co(A)$ if there is a positive integer $n$,
vectors $v_1, \ldots, v_n \in A$, and real numbers $\lambda_1, \ldots,
\lambda_n \in [0, 1]$ such that
\begin{equation}
	w = \sum_{i=1}^n \lambda_i \, v_i \quad\hbox{and}\quad
		\sum_{i=1}^n \lambda_i = 1.
\end{equation}
It is easy to see from the definitions that $\co(A) = A$ when $A$ is a
convex subset of $V$, and that $\co(A)$ is always a convex subset of
$V$.  Also, if $V$ is a real vector space of finite dimension $m$,
then the convex hull of a nonempty subset $A$ of $V$ is equal to
the set of points in $V$ which can be expressed as convex combinations
of at most $m+1$ elements of $A$.

\section{Topological spaces}
\label{section on topological spaces}

	Let $X$ be a nonempty set.  If $A$ is another nonempty set,
then $\{E_\alpha\}_{\alpha \in A}$ defines a \emph{family of subsets
of $X$}\index{family of subsets of a given set} if $E_\alpha$ is
a subset of $X$ for each $\alpha \in A$.  The set $A$ is called the
\emph{index} set of the family.  The union and intersection of the
sets in the family are denoted
\begin{equation}
	\bigcup_{\alpha \in A} E_\alpha, \qquad
		\bigcap_{\alpha \in A} E_\alpha,
\end{equation}
and the union consists of the points in $X$ which lie in $E_\alpha$
for at least one $\alpha \in A$, while the intersection consists
of the points in $X$ which lie in $E_\alpha$ for every $\alpha \in A$.

	If $F$ is a subset of $X$, then the complement of $F$ in $X$
is denoted $X \backslash F$ and it is defined to be the set of points
in $X$ which do not lie in $F$.  Of course the complement of the
complement of $F$ in $X$ is $F$ itself.  If $\{E_\alpha\}_{\alpha \in
A}$ is a family of subsets of $X$ as above, then
\begin{equation}
	X \backslash \biggl(\bigcup_{\alpha \in A} E_\alpha\biggr)
		= \bigcap_{\alpha \in A} (X \backslash E_\alpha),
\end{equation}
and
\begin{equation}
	X \backslash \biggr(\bigcap_{\alpha \in A} E_\alpha\biggr)
		= \bigcup_{\alpha \in A} (X \backslash E_\alpha).
\end{equation}

	Let $X$ be a nonempty set, and let $\tau$ be a collection of
subsets of $X$.  We say that $\tau$ defines a \emph{topology} on
$X$\index{topology on a set} if the empty set $\emptyset$ and $X$
itself lie in $\tau$, if the union of any family of elements of $\tau$
is again an element of $\tau$, and if the intersection of any finite
collection of elements of $\tau$ is also an element of $\tau$.
The combination of $X$ and $\tau$ is called a topological space.

	For any set $X$, we might choose $\tau$ so that it contains
only $\emptyset$ and $X$, i.e., this satisfies the conditions just
mentioned.  This is the smallest possible topology on $X$.  At the
other extreme we might take $\tau$ to be the collection of all subsets
of $X$, which is the largest possible topology on $X$, also known as
the discrete topology on $X$.  We shall normally be interested in
topologies between these, although the discrete topology is sometimes
useful.

	When a nonempty set $X$ is equipped with a topology $\tau$ in
this way, the elements of $\tau$ are called the \emph{open}\index{open
subsets of a topological space} subsets of $X$.  A subset $F$ of $X$
such that $X \backslash F$ is open is said to be a
\emph{closed}\index{closed subset of a topological space} subset of
$X$.  By the conditions mentioned before, the empty set $\emptyset$
and $X$ itself are automatically closed subsets of $X$, the
intersection of any family of closed subsets of $X$ is again a closed
subset of $X$, and the union of finitely many closed subsets of $X$ is
a closed subset of $X$.

	If $(X, \tau)$ is a topological space, $E$ is a subset of $X$,
and $p$ is a point in $X$, then $p$ is said to be a \emph{limit
point}\index{limit point of a subset of a topological space} of $E$ if
for each open subset $U$ of $X$ such that $p \in U$ there is a point
$q \in E$ such that $q \in U$ and $q \ne E$.  One can check that a
subset $F$ of $X$ is closed if and only every limit point of $F$ in
$X$ is also an element of $F$.  

	If $E$ is a subset of the topological space $X$, then the
\emph{closure}\index{closure of a subset of a topological space}
of $E$ is denoted $\overline{E}$ and defined to be the union of $E$
and the set of limit points of $E$.  Thus $p \in \overline{E}$ if and
only if for every open subset $U$ of $X$ such that $p \in U$ there is
a point $q \in E$ which satisfies $q \in U$.  

	A subset $D$ of $X$ is said to be \emph{dense}\index{dense
subset of a topological space} if the closure of $D$ is equal to $X$.
This is the same as saying that $D$ is dense in $X$ if every nonempty
open subset of $X$ contains an element of $D$.

	The \emph{interior}\index{interior of a subset of a topological
space} of a subset $A$ of $X$ is denoted $A^\circ$ and is defined to
be the set of points $x \in A$ such that there is an open subset 
$U$ of $X$ which satisfies $x \in U$ and $U \subseteq A$.  Because
the union of any family of open subsets of $X$ is an open subset
of $X$, the interior of $A$ is always an open subset of $X$, which
may be the empty set.  If $E$ is a subset of $X$, then one can check that
\begin{equation}
	X \backslash \overline{E} = (X \backslash E)^\circ,
\end{equation}
which is to say that the complement of the closure of $E$ in $X$
is equal to the interior of the complement of $E$ in $X$.  As a result,
the closure of $E$ is always a closed subset of $X$.

	It is often natural to assume \emph{separation conditions} on
a topological space.  One of the simplest of these asks that each
subset of the topological space with exactly one element is a closed
subset of the topological space.  This implies automatically that
every finite subset of the topological space is a closed subset.  If
$X$ is a topological space which has this property, if $E$ is a subset
of $X$, and if $p$ is a limit point of $E$ in $X$, then for every open
subset $U$ of $X$ such that $p \in U$, the set $E \cap U$ has
infinitely many elements.

	A topological space $(X, \tau)$ is said to be
\emph{Hausdorff}\index{Hausdorff topological spaces} if for every pair
of distinct points $p, q \in X$ there are open subsets $U$, $W$ of
$X$ such that $p \in U$, $q \in W$, and $U \cap W = \emptyset$.
This condition implies the one described in the previous paragraph.

	If $X$ is a nonempty set and $d(x,y)$ is a real-valued
function on the Cartesian product $X \times X$ of $X$ with itself,
then $d(x,y)$ is said to be a \emph{semimetric}\index{semimetric on a
set} if $d(x,y) \ge 0$ for all $x, y \in X$, $d(x, x) = 0$ for all $x
\in X$,
\begin{equation}
	d(x, y) = d(y, x)
\end{equation}
for all $x, y \in X$, and
\begin{equation}
	d(x, z) \le d(x, y) + d(y, z)
\end{equation}
for all $x, y, z \in X$.  If also $d(x,y) > 0$ when $x \ne y$,
then $d(\cdot, \cdot)$ is said to be a \emph{metric}.\index{metric on
a set}

	If $d(x,y)$ is a semimetric on $X$, then we can define a
topology on $X$ by saying that a subset $U$ of $X$ is open if and only
if for every $x \in U$ there is a positive real number such that
$y \in U$ for all $y \in X$ which satisfy $d(x,y) < r$.  This topological
space has the property that subsets of $X$ with exactly one element
are closed subsets of $X$ if and only if $d(\cdot, \cdot)$ is a metric.
Conversely, if $d(\cdot, \cdot)$ is a metric, then the associated
topology on $X$ is Hausdorff.

	A basic observation about semimetrics is that if $d(x, y)$
is a semimetric on $X$, then for each $p \in X$ and each positive
real number $t$ the open ball centered at $p$ with radius $t$ with
respect to the semimetric, defined by
\begin{equation}
	\{z \in X : d(z, p) < t\},
\end{equation}
is an open subset of $X$ with respect to the topology that comes from
$d(x,y)$.  This is not hard to show, using the triangle inequality.

	As basic examples, the real line ${\bf R}$ has the standard
metric $|x - y|$, and on the complex numbers ${\bf C}$ we have the
standard metric $|z - w|$.  The topologies on ${\bf R}$ and ${\bf C}$
associated to these metrics are called their standard topologies.

	If $\mathcal{F}$ is a nonempty family of semimetrics on $X$,
then we can define a topology on $X$ by saying that a subset $U$ of
$X$ is open if and only if for each point $x \in U$ there are
semimetrics $d_1(\cdot, \cdot), \ldots, d_n(\cdot, \cdot)$ in the
family $\mathcal{F}$ and positive real numbers $r_1, \ldots, r_n$ such
that $y \in U$ when $y \in X$ satisfies $d_j(y,x) < r_j$ for $j = 1,
\ldots, n$.  With this topology, one-element subsets of $X$ are closed
if and only if for every pair of distinct points $p, q \in X$ there is
a semimetric $d(\cdot, \cdot)$ in the family $\mathcal{F}$ such that
$d(p, q) > 0$, and conversely this condition implies that $X$ is
Hausdorff.

	Notice that if $d(x, y)$ is a semimetric on $X$, and $a$, $b$
are positive real numbers such that $b \le 1$, then
\begin{equation}
	\min(d(x,y), 1), \enspace d(x,y)^b, \enspace \frac{d(x,y)}{1 + d(x,y)}
\end{equation}
are semimetrics on $X$ as well.  In terms of defining topologies, as
above, these semimetrics are all equivalent to the original semimetric
$d(x, y)$.

	If $(X, \tau)$ is a topological space and $Y$ is a nonempty
subset of $X$, then we get a natural topology $\tau_Y$ on $Y$ induced
from the one on $X$ by taking $\tau_Y$ to be the collection of subsets
of $Y$ of the form $Y \cap U$, where $U \in \tau$.  In other words,
a subset of $Y$ is considered to be open relative to $Y$ if it is
the intersection of $Y$ with an open subset of $X$.  With respect
to this topology on $Y$, a subset of $Y$ is closed relative to $Y$
if it is the intersection of $Y$ with a closed subset of $X$.
If $X$ has the property that one-element subsets are closed, then
so does $Y$, and if $X$ is Hausdorff, then so is $Y$ too.

	Now suppose that $n$ is a positive integer greater than or
equal to $2$, and that $(X_1, \tau_1), (X_2, \tau_2), \ldots, (X_n,
\tau_n)$ are topological spaces.  Consider the Cartesian product
\begin{equation}
	X_1 \times X_2 \times \cdots \times X_n,
\end{equation}
which consists of $n$-tuples of the form $(x_1, x_2, \ldots, x_n)$
with $x_j \in X_j$ for $j = 1, 2, \ldots, n$.  There is a natural
topology on the Cartesian product of the $X_j$'s, in which a subset
$W$ of the Cartesian product is open if for each point $(x_1, x_2,
\ldots, x_n)$ in $W$ there are open subsets $U_j$ of $X_j$, $1 \le j
\le n$, such that $x_j \in U_j$ for each $j$ and 
\begin{equation}
	U_1 \times U_2 \times \cdots \times U_n \subseteq W.
\end{equation}
In other words, a subset $W$ of $X_1 \times X_2 \times \cdots \times
X_n$ is considered an open subset if it can be expressed as the union
of Cartesian products of open subsets of the $X_j$'s.  This topology
on the Cartesian product of the $X_j$'s is called the \emph{product
topology} associated to the topologies on the $X_j$'s individually.

	In particular, Cartesian products of open subsets of the
$X_j$'s define open subsets of the Cartesian product of the $X_j$'s.
One can also check that Cartesian products of closed subsets of the
$X_j$'s define closed subsets of the Cartesian product.  As a result,
if each $X_j$ has the property that one-element subsets are closed
subsets, then the same holds for the Cartesian product of the $X_j$'s.
Similarly, if each $X_j$ is Hausdorff, then one can verify that the
Cartesian product of the $X_j$'s is Hausdorff.

	Let $X$ and $Y$ be sets, and suppose that $f$ is a mapping
from $X$ to $Y$.  If $A$ is a subset of $X$, then the image of
$A$ under $f$ in $Y$ is denoted $f(A)$ and defined by
\begin{equation}
	f(A) = \{y \in Y : y = f(x) \hbox{ for some } x \in A\}.
\end{equation}
If $\{A_\iota\}_{\iota \in I}$ is a family of subsets of $X$, then
\begin{equation}
	f\biggl(\bigcup_{\iota \in I} A_\iota \biggr)
		= \bigcup_{\iota \in I} f(A_\iota).
\end{equation}
For intersections we have the inclusion
\begin{equation}
	f\biggl(\bigcap_{\iota \in I} A_\iota \biggr)
		\subseteq \bigcap_{\iota \in I} f(A_\iota),
\end{equation}
and equality does not have to hold in general.

	If $E$ is a subset of $Y$, then the inverse image of $E$ under
$f$ in $X$ is denoted $f^{-1}(E)$ and defined by
\begin{equation}
	f^{-1}(E) = \{x \in X : f(x) \in E\}.
\end{equation}
The inverse image behaves nicely in the sense that
for each family $\{E_\beta\}_{\beta \in B}$ of subsets of $Y$
we have that
\begin{equation}
	f^{-1}\biggl(\bigcup_{\beta \in B} E_\beta\biggr)
			= \bigcup_{\beta \in B} f^{-1}(E_\beta)
\end{equation}
and
\begin{equation}
	f^{-1}\biggl(\bigcap_{\beta \in B} E_\beta\biggr)
			= \bigcap_{\beta \in B} f^{-1}(E_\beta),
\end{equation}
and also
\begin{equation}
	f^{-1}(Y \backslash E) = X \backslash f^{-1}(E)
\end{equation}
for each subset $E$ of $Y$.

	Let $(X, \sigma)$ and $(Y, \tau)$ be topological spaces, and
let $f$ be a mapping from $X$ to $Y$.  If $p$ is a point in $X$, then
we say that $f$ is \emph{continuous at $p$}\index{continuity at a
point for a mapping between topological spaces} if for every open
subset $W$ of $Y$ such that $f(p) \in W$ there is an open subset $U$
of $X$ such that $U \subseteq f^{-1}(W)$.  We say that $f : X \to Y$
is \emph{continuous}\index{continuous mappings between topological
spaces} if $f$ is continuous at every point $p$ in $X$.  This is
equivalent to saying that $f^{-1}(W)$ is an open subset of $X$
for every open subset $W$ of $Y$, or that $f^{-1}(E)$ is a closed 
subset of $X$ for every closed subset $E$ of $Y$.

	If $f : X \to Y$ is continuous, then notice that
\begin{equation}
	f(\overline{A}) \subseteq \overline{f(A)}
\end{equation}
for every subset $A$ of $X$.  In other words, the image of the closure
of $A$ in $X$ under $f$, which is a subset of $Y$, is contained in the
closure of the image of $A$ under $f$ in $Y$.  Another nice property
of continuous mappings is that if $f : X \to Y$ is continuous, then
the graph of $f$ in $X \times Y$, defined by
\begin{equation}
	\{(x, y) \in X \times Y : y = f(x)\},
\end{equation}
is a closed subset of $X \times Y$ with respect to the topology on $X
\times Y$ obtained from the topologies on $X$, $Y$ as above.  The
space of continuous mappings from $X$ to $Y$ will be denoted
$\mathcal{C}(X, Y)$.

	Suppose that $n$ is a positive integer, $n \ge 2$, and that
$(X_j, \sigma_j)$, $(Y_j, \tau_j)$, $1 \le j \le n$, are topological
spaces.  Assume also that for each $j = 1, 2, \ldots, n$ we have a
mapping $f_j : X_j \to Y_j$.  We can combine these mappings in a
natural way to get a mapping $F$ from the Cartesian product of the
$X_j$'s into the Cartesian product of the $Y_j$'s, defined by
\begin{equation}
	(x_1, x_2, \ldots, x_n) \mapsto (f_1(x_1), f_2(x_2), \ldots, f_n(x_n)).
\end{equation}
If for each $j$ we have a point $p_j$ in $X_j$ and $f_j$ is continuous
at $p_j$, then $F$ is continuous as a mapping between the Cartesian
products at the point $(p_1, p_2, \ldots, p_n)$.  If $f_j : X_j \to Y_j$
is continuous for each $j$, then
\begin{equation}
	F : X_1 \times X_2 \times \cdots \times X_n \to
		Y_1 \times Y_2 \times \cdots \times Y_n
\end{equation}
is also continuous.

	Let $(X_1, \tau_1)$, $(X_2, \tau_2)$, and $(X_3, \tau_3)$
be topological spaces.  Assume that $f_1$ is a mapping from $X_1$
to $X_2$, and that $f_2$ is a mapping from $X_2$ to $X_3$.  If
$p$ is an element of $X_1$, $f_1$ is continuous at $p$, and $f_2$
is continuous at $f_1(p)$, then it is easy to check that the
composition $f_2 \circ f_1$, which is a mapping from $X_1$ to $X_3$,
is continuous at $p$.  If $f_1$ is continuous as a mapping from
$X_1$ to $X_2$ and $f_2$ is continuous as a mapping from $X_2$
to $X_3$, then $f_2 \circ f_1$ is continuous as a mapping from
$X_1$ to $X_3$.

	If $(X, \sigma)$ and $(Y, \tau)$ are topological spaces and
$f$ is a one-to-one mapping from $X$ onto $Y$, so that the inverse
mapping $f^{-1}$ from $Y$ to $X$ exists, then $f$ is said to be a
\emph{homeomorphism}\index{homeomorphism from one topological space
onto another} from $X$ onto $Y$ if $f$, $f^{-1}$ are both continuous
as mappings from $X$ to $Y$ and from $Y$ to $X$, respectively.
In this case $f^{-1}$ is also a heomeomorphism from $Y$ onto $X$.

	If $(X, \tau)$ is a topological space, and if $f_1$, $f_2$ are
two real-valued continuous functions on $X$, or two complex-valued
continuous functions on $X$, then the sum $f_1 + f_2$ and the product
$f_1 \, f_2$ are also continuous functions on $X$.  One way to look at
this is as a consequence of the fact that the mappings from ${\bf R}
\times {\bf R}$ to ${\bf R}$ and from ${\bf C} \times {\bf C}$ to
${\bf C}$ given by addition and multiplication are continuous.

\section{Countability conditions}
\label{section on countability conditions}

	Let $E$ be a set.  To be a bit formal, we say that $E$ is
finite if it is empty, in which case it has $0$ elements, or if there
is a positive integer $n$ such that there is a one-to-one
correspondence between $E$ and the set $\{1, 2, \ldots, n\}$ of the
first $n$ positive integers, in which case we say that $E$ has $n$
elements.  If there is a one-to-one correspondence between $E$ and the
set ${\bf Z}_+$ of positive integers, then we say that $E$ is
countably-infinite, and we say that $E$ is at most countable if $E$ is
either finite or countable.  Of course the union of finitely many
finite sets is a finite set, and if $\{E_\alpha\}{\alpha \in A}$ is a
family of subsets of a set $X$ such that the index set $A$ is at most
countable and each $E_\alpha$, $\alpha \in A$, is at most countable,
then the union $\bigcup_{\alpha \in A} E_\alpha$ is also at most countable.

	Now let $(X, \sigma)$ be a topological space, and let
$\mathcal{B}$ be a collection of open subsets of $X$.  We say that
$\mathcal{B}$ is a \emph{basis for the topology of $X$}\index{basis
for the topology of a topological space} if every open subset of $X$
is the union of some collection of open sets in $\mathcal{B}$.  If $p$
is a point in $X$, then we say that $\mathcal{B}$ is a \emph{local
basis for the topology of $X$ at $p$}\index{local basis for the
topology of a topological space at a point} if for every open subset
$U_1$ of $X$ such that $p \in U_1$ there is an open subset $U_2$ of
$X$ such that $p \in U_2$, $U_2 \subseteq U_1$, and $U_2 \in
\mathcal{B}$.  It is easy to see that $\mathcal{B}$ is a basis for
the topology of $X$ if and only if $\mathcal{B}$ is a local basis for the
topology of $X$ at each point in $X$.

	If $p$ is a point in $X$, then we say that $X$ satisfies the
\emph{first axiom of countability at $p$}\index{first axiom of
countability at a point in a topological space} if there is a local
basis for the topology of $X$ at $p$ which has at most countably many
element.  We say that $X$ satisfies the \emph{first axiom of
countability}\index{first axiom of countability for a topological
space} if $X$ satisfies the first axiom of countability at each point
in itself.  We say that $X$ satisfies the \emph{second axiom of
countability}\index{second axiom of countability for a topological
space} if there is a basis for the topology of $X$ which has at most
countably many elements.  Notice that the second axiom of countability
implies the first axiom of countability.

	Suppose that $\{p_j\}_{j=1}^\infty$ is a sequence of points in
$X$, and that $p$ is a point in $X$.  We say that
$\{p_j\}_{j=1}^\infty$ \emph{converges}\index{convergence of a
sequence in a topological space} to $p$ if for every open subset
$U$ of $X$ such that $p \in U$, there is a positive integer
$N$ so that $p_j \in U$ when $j \ge N$.  In this case we write
\begin{equation}
	\lim_{j \to \infty} p_j = p.
\end{equation}
If $X$ has the property that one-element subsets of $X$ are closed
subsets of $X$, then a sequence of points in $X$ can have at most one
limit, i.e., the limit is unique when it exists.

	Let $p$ be an element of $X$, let $E$ be a subset of $X$, and
assume that $X$ satisfies the first axiom of countability at $p$.  The
$p$ is an element of the closure of $E$ if and only if there is a
sequence $\{p_j\}_{j=1}^\infty$ of points in $E$ which converges to
$p$.  If $p \in E$, we can simply take $p_j = p$ for all $j$ here.

	Suppose that $(X, \sigma)$ and $(Y, \tau)$ are topological
spaces, $p$ is a point in $X$, and $f$ is a mapping from $X$ to $Y$.
If $X$ satisfies the first axiom of countability at $p$, then $f$ is
continuous at $p$ if and only if for every sequence
$\{p_j\}_{j=1}^\infty$ of points in $X$ which converges to $p$,
we have that $\{f(p_j)\}_{j=1}^\infty$ converges to $f(p)$ in $Y$.

	Let $(X, \sigma)$ be a topological space, and let
$\mathcal{B}$ be a collection of open subsets of $X$.  Define
$\widehat{\mathcal{B}}$ to be the collection of subsets of $X$ which
are intersections of finitely many subsets of $X$ in $\mathcal{B}$.
Thus $\widehat{\mathcal{B}}$ is at most countable if $\mathcal{B}$ is
at most countable.

	The collection $\mathcal{B}$ of open subsets of $X$ is said
to be a \emph{sub-basis for the topology of $X$}\index{sub-basis for
the topology of a space} if $\widehat{\mathcal{B}}$ is a basis for
the topology of $X$.  Similarly, if $p$ is a point in $X$, then
$\mathcal{B}$ is said to be a \emph{local sub-basis for the topology
of $X$ at $p$}\index{local sub-basis for the topology of a space
at a point} if $\widehat{\mathcal{B}}$ is a local basis for the
topology of $X$ at $p$.  It is easy to see that $\mathcal{B}$ is
a sub-basis for the topology of $X$ if and only if $\mathcal{B}$
is a local sub-basis for the topology of $X$ at each point in $X$.

	One can go in the other direction and start with a nonempty
set $X$ and a collection $\mathcal{B}$ of subsets of $X$, define
$\widehat{\mathcal{B}}$ as before, and then define a topology
on $X$ so that $\mathcal{B}$ is a sub-basis for that topology,
namely, where a subset of $X$ is open if it is a union subsets of
$X$ in $\widehat{\mathcal{B}}$.  For this to work one should assume
that the union of the subsets of $X$ in $\mathcal{B}$ is all of $X$,
and one can interpret the empty subset of $X$ as the empty union
of subsets of $X$ in $\widehat{\mathcal{B}}$, although it is also
easy for the empty set to occur as an intersection of finitely
many subsets of $X$ in $\mathcal{B}$.

	A basic situation of this type arises with families of
semimetrics, as discussed in the previous section.  Specifically, if
$\mathcal{F}$ is a nonempty family of semimetrics on $X$, then the
topology on $X$ associated to this topology is the same as the
topology for which the collection $\mathcal{B}$ of open balls in $X$
with respect to the semimetrics in $\mathcal{F}$ is a sub-basis.

	If $\mathcal{F}$ is at most countable, then it is easy to see
that the corresponding topology on $X$ satisfies the first axiom of
countability.  In fact, if $\mathcal{F}$ is at most countable, then
there are simple tricks for defining a single semimetric on $X$ with
the same topological information as the family $\mathcal{F}$ of
semimetrics.  However, it may be that the semimetrics in $\mathcal{F}$
have some nice features and that one would like to keep them.

	Notice that if $(X_1, \tau_1), \ldots, (X_n, \tau_n)$ are
topological spaces, $p_1, \ldots, p_n$ are points in $X_1, \ldots,
X_n$, respectively, and $X_j$ satisfies the first axiom of
countability at $p_j$ for $j = 1, \ldots, n$, then the product $X_1
\times \cdots \times X_n$, with the product topology discussed in the
previous section, satisfies the first axiom of countability at the
point $(p_1, \ldots, p_n)$.  Hence if $X_1, \ldots, X_n$ satisfy the
first axiom of countability at all of their points, then so does the
product $X_1 \times \cdots \times X_n$.  Similarly, if $X_1, \ldots,
X_n$ satisfy the second axiom of countability, then the product $X_1
\times \cdots \times X_n$ does too.

	A topological space $(X, \tau)$ is said to be
\emph{separable}\index{separable topological space} if it has a subset
which is at most countable and also dense.  If $X$ satisfies the
second axiom of countability, then $X$ is separable.  Conversely, if
the topology on $X$ comes from an at most countable family of
semimetrics, and if $X$ is separable, then $X$ satisfies the second
axiom of countability.

	Notice that the notions of basis and sub-basis and the
countability conditions for a topological space behave well in terms
of restricting from that space to a subspace.  To be more precise, let
$(X, \tau)$ be a topological space, let $Y$ be a nonempty subset of
$X$, and let $\tau_Y$ be the topology on $Y$ induced from the one on
$X$ by saying that a subset $E$ of $Y$ is open relative to $Y$ if it
is of the form $U \cap Y$, where $U$ is an open subset of $X$.  Suppose 
that $\mathcal{B}$ is a collection of open subsets of $X$, and put
\begin{equation}
	\mathcal{B}_Y = \{ U \cap Y : U \in \mathcal{B}\}.
\end{equation}
If $p$ is a point in $Y$ and $\mathcal{B}$ is a local basis for the
topology of $X$ at $p$, then $\mathcal{B}_Y$ is a local basis for the
topology of $Y$ at $p$.  Similarly, if $\mathcal{B}$ is a local
sub-basis for the topology of $X$ at $p$, then $\mathcal{B}_Y$ is a
local sub-basis for the topology of $Y$ at $p$.  If $\mathcal{B}$ is a
basis for the topology of $X$, then $\mathcal{B}_Y$ is a basis for the
topology of $Y$, and if $\mathcal{B}$ is a sub-basis for the topology
of $X$, then $\mathcal{B}_Y$ is a sub-basis for the topology of $Y$.
As a result, if $p \in Y$ and $X$ satisfies the first axiom of
countability at $p$, then $Y$ also satisfies the first axiom of
countability at $p$.  If $X$ satisfies the first axiom of
countability, then $Y$ satisfies the first axiom of countability, and
if $X$ satisfies the second axiom of countability, then $Y$ satisfies
the second axiom of countability.

\section{Compactness}
\label{section on compactness}

	Let $(X, \tau)$ be a topological space, and let $E$ be a
subset of $X$.  By an \emph{open covering}\index{open covering of a
subset of a topological space} of $E$ in $X$ we mean a family
$\{U_\alpha\}_{\alpha \in A}$ of open subsets of $X$ such that $E
\subseteq \bigcup_{\alpha \in A} U_\alpha$.  We say that $E$ is
\emph{compact}\index{compact subsets of a topological space} if for
every open covering of $E$ in $X$ there is a finite subcovering of $E$
from that open covering.  In other words, if $\{U_\alpha\}_{\alpha \in
A}$ is an open covering of $E$ in $X$, then we ask that there be a
finite subset $A_1$ of $A$ such that $E \subseteq \bigcup_{\alpha \in
A_1} U_\alpha$.

	If $(X, \sigma)$, $(Y, \tau)$ are topological space,
$f$ is a continuous mapping from $X$ to $Y$, and $E$ is a compact
subset of $X$, then $f(E)$ is a compact subset of $Y$.

	There are a number of related notions and variants of
compactness that one sometimes likes to consider.  Let us say that a
subset $E$ of $X$ is \emph{countably compact}\index{countably compact
subset of a topological space} if for every open covering
$\{U_\alpha\}_{\alpha \in A}$ of $E$ in $X$ there is a subcovering of
$E$ from this covering which is at most countable, which is to say that
there is a subset $A_2$ of $A$ which is at most countable such that
$E \subseteq \bigcup_{\alpha \in A_2} U_\alpha$.

	A basic fact is that if $(X, \tau)$ is a topological space,
$Y$ is a nonempty subset of $X$, and $E$ is a subset of $Y$, then $E$
is compact as a subset of $X$ if and only if $E$ is compact as a
subset of $Y$, using the topology on $Y$ induced from the one on $X$.
In the same way, if $E$ is countably compact as a subset of $X$,
then $E$ is countably compact as a subset of $Y$.

	It is well-known and not hard to show that if $X$ satisfies
the second axiom of countability, then every subset of $X$ is
countably compact.  Thus if $(X, \tau)$ is a topological space which
satisfies the second axiom of countability and if $E$ is a subset of
$X$, then $E$ is compact if and only if for every open covering
$\{U_\alpha\}_{\alpha \in A}$ of $E$ with index set $A$ which is at
most countable, there is a finite subcovering of $E$ from this open
covering.

	If $(X, \tau)$ is a topological space and $E$ is a subset of
$X$, then $E$ has the \emph{limit point property}\index{limit point
property for a subset of a topological space} if every infinite subset
$L$ if $E$ has a limit point in $X$ which is also an element of $E$.
If $E$ is a compact subset of $X$, then $E$ satisfies the limit point
property.  To see this, assume for the sake of a contradiction that
$E$ does not satisfy the limit point property.  This means that there
is an infinite subset $L$ of $E$ such that no element of $E$ is a
limit point of $L$.  This implies in turn that for each $p \in E$
there is an open subset $U_p$ of $X$ such that $p \in U_p$ and $U_p
\cap L$ is either empty or contains $p$ only.  Thus $\{U_p\}{p \in E}$
is an open covering of $E$ in $X$, and the compactness of $E$ implies
that there is a finite subset $E_1$ of $E$ such that $E \subseteq
\bigcup_{p \in E_1} U_p$.  Because $L \subseteq E$ and $L \cap U_p$
has at most one element for each $p \in E$, it follows that $L$
is finite, a contradiction.  Thus compactness implies the limit
point property.

	As for compactness, if $(X, \tau)$ is a topological
space, $Y$ is a nonempty subset of $X$, and $E$ is a subset of $Y$,
then $E$ has the limit point property as a subset of $X$ if and only
if $E$ has the limit point property as a subset of $Y$, equipped with
the topology induced from the one on $X$.  This is because a point
in $E$ is a limit point of a subset $L$ of $E$ in the topological
space $X$ if and only if this holds in the topological space $Y$.

	Let $(X, \tau)$ be a topological space and let $E_1$, $E_2$ be
subsets of $X$.  If $E_1$, $E_2$ are compact, then so is $E_1 \cup
E_2$.  If $E_1$, $E_2$ have the limit point property, then $E_1 \cup
E_2$ has the limit point property too.  If $E$, $F$ are subsets of $X$
such that $F \subseteq E$, $F$ is closed, and $E$ is compact, then $F$
is compact.  If $E$, $F$ are subsets of $X$ such that $F \subseteq E$,
$F$ is closed, and $E$ has the limit point property, then $F$ has the
limit property.  It is immediate from the definitions that finite sets
are compact and have the limit point property.

	Here is a reformulation of compactness which is sometimes
useful.  A topological space $(X, \tau)$ is itself compact if and only
if for every family $\{F_\alpha\}_{\alpha \in A}$ of closed subsets of
$X$ such that $\bigcap_{\alpha \in A_0} F_\alpha \ne \emptyset$
whenever $A_0$ is a nonempty finite subset of $A$, we have that
$\bigcap_{\alpha \in A} F_\alpha \ne \emptyset$.  If $(X, \tau)$
satisfies the second axiom of countability, then we can restrict our
attention to families of closed subsets of $X$ which are at most
countable, and which can then be arranged in a sequence.  In fact we
get that if $(X, \tau)$ satisfies the second axiom of countability,
then $X$ is compact if and only if for every sequence
$\{F_j\}_{j=1}^\infty$ of nonempty closed subsets of $X$ such that
$F_{j+1} \subseteq F_j$ for all $j$ we have that $\bigcap_{j=1}^\infty
F_j \ne \emptyset$.

	Suppose that $(X, \tau)$ is a Hausdorff topological space.
Let us check that if $E$ is a compact subset of $X$, then $E$ is
closed.  Assume for the sake of a contradiction that $p$ is a
limit point of $E$ which is not an element of $E$.  For each point
$q \in E$, let $U(q)$, $V(q)$ be disjoint open subsets of $X$
such that $q \in U(q)$ and $p \in V(q)$.  Thus $\{U(q)\}_{q \in E}$
is an open covering of $E$, and therefore there is a finite
subset $E_1$ of $E$ such that 
\begin{equation}
	E \subseteq \bigcup_{q \in E_1} U(q).
\end{equation}
This implies that $E$ is disjoint from $\bigcap_{q \in E_1} V(q)$,
which is an open subset of $X$ that contains $p$.  This contradicts
the assumption that $p$ is a limit point of $E$, and it follows
that $E$ is a closed subset of $X$.

	A subset $E$ of a topological space $X$ is said to be
\emph{sequentially compact}\index{sequentially compact subset of a
topological space} if for every sequence $\{x_j\}_{j=1}^\infty$ of
points in $E$ there is a subsequence $\{x_{j_l}\}_{l=1}^\infty$ and a
point $x \in E$ such that $\{x_{j_l}\}_{l=1}^\infty$ converges to $x$.
It is easy to see that finite sets are sequentially compact, and that
if $E$ is sequentially compact, then $E$ has the limit point property.
The union of two sequentially compact subsets of $X$ is also
sequentially compact, and a closed subset of a sequentially compact
subset of $X$ is sequentially compact.  If $(X, \tau)$ is a
topological space, $Y$ is a nonempty subset of $X$, and $E$ is a
subset of $Y$, then $E$ is sequentially compact as a subset of $X$
if and only if $E$ is sequentially compact as a subset of $Y$,
equipped with the topology induced from the one on $X$.

	Suppose that $(X, \tau)$ is a topological space which
satisfies the first axiom of countability and in which one-element
subsets are closed.  Let $E$ be a subset of $X$ which satisfies the
limit point property, which holds in particular when $E$ is compact,
and let us show that $E$ is sequentially compact.  Let
$\{x_j\}_{j=1}^\infty$ be a sequence of points in $E$, and let $L$
denote the subset of $X$ consisting of the points in this sequence.
If $L$ is a finite set, then there is a subsequence of
$\{x_j\}_{j=1}^\infty$ in which all of the terms are the same, and
this subsequence converges trivially.  Thus we assume that $L$ is
infinite.  Because $E$ satisfies the limit point property, there is a
point $x \in E$ which is a limit point of $L$.  It is not difficult to
show that there is then a subsequence of $\{x_j\}_{j=1}^\infty$ which
converges to $x$.

	If $(X, \tau)$ is a topological space which satisfies the
first axiom of countability and in which one-element subsets are closed,
and if $E$ is a sequentially-compact subset of $X$, then $E$ is a
closed subset of $X$.

	If $(X, \tau)$ is a topological space which satisfies the
second axiom of countability and if $X$ is sequentially compact, then
$X$ is compact.  Indeed, as mentioned above, it suffices to show that
if $\{F_j\}_{j=1}^\infty$ is a sequence of nonempty subsets of $X$
such that $F_{j+1} \subseteq F_j$ for all $j$, then
$\bigcap_{j=1}^\infty F_j \ne \emptyset$.  For each positive integer
$j$, let $x_j$ be an element of $F_j$.  This gives a sequence
$\{x_j\}_{j=1}^\infty$ of points in $X$, and sequential compactness
implies that there is a subsequence of this sequence which converges
to a point $x \in X$.  It is easy to see that $x \in F_j$ for all $j$,
as desired.

	Suppose that $(X_1, \tau_1)$, $(X_2, \tau_2)$ are topological
spaces, and consider the Cartesian product $X_1 \times X_2$ equipped
with the product topology coming from the topologies on $X_1$, $X_2$.
If $E_1$, $E_2$ are compact subsets of $X_1$, $X_2$, respectively,
then $E_1 \times E_2$ is a compact subset of $X_1 \times X_2$.  If
$E_1$, $E_2$ are sequentially compact subsets of $X_1$, $X_2$, then
$E_1 \times E_2$ is a sequentially compact subset of $X_1 \times X_2$.
These statements are not too difficult to show.

	A fundamental result states that closed and bounded subsets
of ${\bf R}^n$ are compact.

\chapter{Topological vector spaces}
\label{chapter on topological vector spaces}

\section{Basic notions}
\label{section on basic notions, tvs's}

	To say that $V$ is a \emph{topological vector
space}\index{topological vector space} means that $V$ is a vector
space, that $V$ is also equipped with a topology, and that with
respect to this topology the vector space operations
\begin{equation}
	\hbox{addition} : V \times V \to V,
\end{equation}
and
\begin{equation}
	\hbox{scalar multiplication} : {\bf R} \times V \to V
\end{equation}
if $V$ is a real vector space or
\begin{equation}
	\hbox{scalar multiplication} : {\bf C} \times V \to V
\end{equation}
if $V$ is a complex vector space, are continuous.

	It is easy to see from the definition that if $V$ is
a topological vector space and $v_0$ is an element of $V$,
then the translation mapping 
\begin{equation}
	v \mapsto v + v_0
\end{equation}
defines a homeomorphism from $V$ onto itself.  Similarly,
if $\alpha$ is a nonzero scalar, then the dilation mapping
\begin{equation}
	v \mapsto \alpha \, v
\end{equation}
defines a homeomorphism from $V$ onto itself.

	The continuity of the binary operation of vector addition
at $(0,0)$ in $V \times V$ is equivalent to the statement that
for each open subset $U_1$ of $V$ such that $0 \in U_1$, there
is an open subset $U_2$ of $V$ such that $0 \in U_2$ and
\begin{equation}
	U_2 + U_2 \subseteq U_1.
\end{equation}
This condition is often useful as a way of saying that $U_2$ is about
``half'' the size of $U_1$, or smaller.  If $U_1$ is convex, then one
can simply take 
\begin{equation}
	U_2 = \frac{1}{2} \, U_1,
\end{equation}
but this does not work in general.  

	If $E$ is a nonempty subset of $V$ and $U$ is an open subset
of $V$ such that $0 \in U$, then consider the set
\begin{equation}
	E + U.
\end{equation}
This is an open subset of $V$ which contains $E$.  In fact, one can
check that the closure of $E$ is contained in $E + U$.  If
$\mathcal{B}_0$ is a local basis for the topology of $V$ at $0$,
then one can also check that
\begin{equation}
	\overline{E} = \bigcap_{U \in \mathcal{B}_0} (E + U).
\end{equation}

	If $U_1$ is an open subset of $V$ such that $0 \in U_1$,
then there is an open subset $U_2$ of $V$ such that $0 \in U_2$ and
the closure of $U_2$ is contained in $U_1$.  Indeed, by the remarks
of the preceding paragraph, it is enough to choose $U_2$ so that
$U_2 + U_2 \subseteq U_1$, as before.

	The continuity of scalar multiplication at $(0,0)$ is equivalent
to the statement that for each open subset $W$ of $V$ such that $0 \in W$
there is a positive real number $r$ and an open subset $U$ of $V$
such that $0 \in U$ and $\alpha \, v \in W$ whenever $\alpha$ is a
scalar which satisfies $|\alpha| < r$ and $v \in U$.  If we put
\begin{equation}
	\widehat{U} = \bigcup_{|\alpha| < r} \alpha U,
\end{equation}
where the union is taken over scalars $\alpha$, then $\widehat{U}$
is an open subset of $V$ and 
\begin{equation}
	\widehat{U} \subseteq W.
\end{equation}
Moreover, $U$ is starlike around $0$, symmetric when $V$ is a real
vector space, and circular when $V$ is a complex vector space.

	Let $U$ be an open subset of $V$ such that $0 \in U$,
and let $v$ be a nonzero vector in $V$.  Using the continuity of
scalar multiplication, it follows that there is a positive real
number $r$ so that $\alpha \, v \in U$ when $\alpha$ is a scalar
which satisfies $|\alpha| < r$.  As a result,
\begin{equation}
	\bigcup_{n=1}^\infty n \, U = V.
\end{equation}

\beginlemma
If $V$ is a topological vector space and $\{0\}$ is a closed
subset of $V$, then $V$ is Hausdorff as a topological space.
\end{lemma}

	From now on in this monograph we make the standing assumption
that $\{0\}$ is a closed subset in a topological vector space, so that
topological vector spaces are Hausdorff, unless otherwise stated in
some situation

	If $V$ is a topological vector space and $\{0\}$ is a closed
subset of $V$, then every one-element subset of $V$ is closed, because
translations on $V$ are homeomorphisms.  To show that $V$ is
Hausdorff, let $p$, $q$ be a pair of distinct elements of $V$.  We
would like to show that there are disjoint open subsets of $V$ that
contain $p$, $q$.  We may as well assume that $p = 0$, so that $q$ is
any element of $V$ which is different from $0$, by using translations
again.

	Let $W$ be the open subset of $V$ defined by $W = V \backslash \{q\}$.
Because $V$ is a topological vector space, there is an open subset $U$ of
$V$ such that $0 \in U$ and 
\begin{equation}
	U - U \subseteq W.
\end{equation}
It follows that $U$ and $U + q$ are disjoint open subsets of $V$
containing $0$, $q$, respectively, as desired.  This proves the lemma.

	Suppose that $V$ is a topological vector space and that $L$
is a linear subspace of $V$.  Then $L$ also inherits a topology from
the one on $V$.  It is easy to see that $L$ is a topological vector
space in its own right, with respect to this topology.

	Now suppose that $V_1, \ldots, V_n$ are topological vector
spaces, all real or all complex.  Consider the direct sum of the
$V_j$'s, which is defined as a vector space by taking the Cartesian
product
\begin{equation}
	V_1 \times V_2 \times \cdots \times V_n
\end{equation}
and using coordinatewise addition and scalar multiplication.  The
topologies on the $V_j$'s lead to the product topology on the direct
sum, and one can check that the direct sum becomes a topological
vector space in this way.

	If $V$ is a real topological vector space, then one can
complexify it first as a vector space by taking the direct sum of $V$
with itself, and defining multiplication by $i$ on the direct sum by
$i (v_1, v_2) = (-v_2, v_1)$.  The topology on $V$ again leads to the
product topology on the direct sum, and one can verify that the
complexification of $V$ becomes a complex topological vector space.

	Let $X$ be a topological space and $V$ a topological vector
space, and recall that $\mathcal{C}(X, V)$ denotes the space of
continuous mappings from $X$ to $V$.  One can check that sums and
scalar multiples of continuous mappings from $X$ to $V$ are again
continuous mappings.  In other words, $\mathcal{C}(X, V)$ is a vector
space in a natural way, which is real if $V$ is real and complex if
$V$ is complex.

\section{Norms, seminorms, and local convexity}
\label{section on norms, seminorms, and local convexity}

	Let $V$ be a vector space.  A real-valued function $N(v)$
on $V$ is said to be a \emph{seminorm}\index{seminorms on vector spaces}
if $N(v) \ge 0$ for all $v \in V$, if the homogeneity condition
\begin{equation}
	N(\alpha \, v) = |\alpha| \, N(v)
\end{equation}
holds for all scalars $\alpha$, which is to say for all real or
complex numbers $\alpha$, as appropriate, and for all $v \in V$, and
if the triangle inequality
\begin{equation}
	N(v + w) \le N(v) + N(w)
\end{equation}
holds for all $v, w \in V$.  If also $N(v) > 0$ when $v \ne 0$, then
we say that $N(v)$ defines a \emph{norm}\index{norms on vector spaces}
on $V$.

\beginremark
{\rm If $N(v)$ is a nonnegative real-valued function on $V$ which
satisfies the homogeneity condition above, then $N(v)$ satisfies the
triangle inequality, and hence is a seminorm, if and only if $N(v)$
is a convex function on $V$.  This is also equivalent to the associated
open unit ball
\begin{equation}
	\{v \in V : N(v) < 1 \}
\end{equation}
being a convex subset of $V$.  These statements are not difficult to
verify.
}
\end{remark}

	If $N(v)$ is a seminorm on $V$, then
\begin{equation}
	|N(v) - N(w)| \le N(v - w)
\end{equation}
for all $v, w \in V$, as one can check using the triangle inequality.
If $V$ is a topological vector space, then it follows that $N(v)$ is
continuous as a real-valued function on $V$ if and only if
the open unit ball with respect to $N(v)$ is an open subset of $V$.

	If $N(v)$ is a seminorm on $V$, then
\begin{equation}
	d(v, w) = N(v - w)
\end{equation}
defines a semimetric on $V$.  If $N(v)$ is a norm, then
$d(v, w)$ as just defined is a metric on $V$.

	A vector space $V$ equipped with a choice of norm $N(v)$ is
called a \emph{normed vector space}.\index{normed vector space} The
norm leads to a metric and hence a topology on $V$, and it is easy to
see that $V$ becomes a topological vector space in this manner.  It is
not too difficult to check that two norms $N_1(v)$, $N_2(v)$ on a
vector space $V$ determine the same topology on $V$ if and only if
there is a positive real number such that $N_1(v) \le C \, N_2(v)$ and
$N_2(v) \le C \, N_1(v)$ for all $v \in V$.

	Suppose that $\mathcal{F}$ is a nonempty family of 
seminorms on $V$.  This leads to a family of semimetrics on $V$
as above.  Let us assume that for each $v \in V$ such that $v \ne 0$
there is a seminorm $N \in \mathcal{F}$ such that $N(v) \ne 0$.
Then the family of semimetrics on $V$ associated to the seminorms
in $\mathcal{F}$ defines a Hausdorff topology on $V$, and $V$
becomes a topological vector space in this way.

	A topological vector space $V$ is said to be \emph{locally
convex}\index{locally convex topological vector space} if for each
open subset $W$ of $V$ such that $0 \in W$ there is an open subset $U$
of $V$ such that $0 \in U$, $U \subseteq V$, and $U$ is convex.  This
is equivalent to saying that there is a local basis for the topology
of $V$ at $0$ consisting of convex open subsets of $V$, and it is also
equivalent to saying that there is a basis for the topology of $V$
consisting of convex open subsets of $V$.  If $V$ is a topological
vector space where the topology can be defined by a family of
seminorms, then $V$ is locally convex.

	Let $V$ be a vector space, and let $E$ be a subset of $V$.
Assume that $0 \in E$, that $E$ is starlike around $0$, and that $E$
is symmetric when $V$ is a real vector space and $E$ is circular when
$V$ is a complex vector space.  Assume also that
\begin{equation}
	\bigcup_{n=1}^\infty n \, E = V.
\end{equation}
Under these conditions, the \emph{Minkowski function associated
to $E$} is the function $\mu_E(v)$ defined on $V$ by
\begin{equation}
	\mu_E(v) = \inf \{ t > 0 : t^{-1} \, v \in E\}.
\end{equation}
Thus $\mu_E(v)$ is a nonnegative real-valued function on $V$ such that
$\mu_E(0) = 0$, $\mu(\alpha \, v) = |\alpha| \, \mu_E(v)$, $\mu(v) \le
1$ for all $v \in E$, and $v \in E$ when $\mu(v) < 1$.  If $E$ is
convex, then $\mu_E(v)$ defines a seminorm on $V$.

	Now suppose that $V$ is a topological vector space
and that $U$ is an open subset of $V$ such that $0 \in U$ and
$U$ is convex.  In the case where the scalars are real numbers,
if we put 
\begin{equation} 
	U_1 = U \cap (-U),
\end{equation}
then $U_1$ is an open subset of $V$ such that $0 \in U_1$, $U_1
\subseteq U$, $U_1$ is convex, and $U_1$ is symmetric.  In the case
where the scalars are complex numbers, if we put
\begin{equation}
	U_1 = \{v \in V : \theta \, v \in U \hbox{ for all } \theta \in {\bf C}
					\hbox{ such that } |\theta| = 1 \},
\end{equation}
then one can check that $U_1$ is an open subset of $V$ such that $0
\in U_1$, $U_1 \subseteq U$, $U_1$ is convex, and $U_1$ is circular.
In both cases we have that the Minkowski function $\mu_{U_1}(v)$ is a
continuous function on $V$ whose open unit ball is contained in $U$.
From this it follows that every locally convex topological vector space
can have the topology described by a family of seminorms.

\section{Bounded subsets of a topological vector space}
\label{section on bounded subsets of topological vector spaces}

	Let $V$ be a topological vector space, and let $E$ be a subset
of $V$.  We say that $E$ is \emph{bounded}\index{bounded subsets of
topological vector spaces} if for every open subset $U$ of $V$ such
that $0 \in U$, there is a positive real number $t$ so that $E
\subseteq t \, U$.  Notice that a compact subset of $V$ is bounded.

	If $V$ is a normed vector space, so that the topology of
$V$ comes from a norm $N(v)$ on $V$, then it is not hard to see that
a subset $E$ of $V$ is bounded if and only if the set of real numbers
$N(v)$, $v \in E$, is bounded from above.  More generally, if the topology
on $V$ comes from a nonempty family $\mathcal{F}$ of seminorms on $V$,
then $E$ is a bounded subset of $V$ if and only if for each seminorm
$N$ in $\mathcal{F}$ the set of real numbers $N(v)$, $v \in E$,
is bounded.

	Suppose that $E$ is a bounded subset of $V$.  It is easy to
see that $\alpha \, E$ is then bounded for all scalars $\alpha$.  Let
us check that the closure $\overline{E}$ of $E$ is also bounded.  Let
$W$ be any open subset of $V$ such that $0 \in W$.  We would like to
show that there is a positive real number $t$ such that $E \subseteq t
\, W$.  To do this we use the fact that there is an open subset $U$ of
$V$ such that $0 \in U$ and $\overline{U} \subseteq W$.  Because $E$
is bounded, there is a positive real number $t$ such that $E \subseteq
t \, U$, and that implies that $\overline{E} \subseteq t \, W$.

	Now assume that $E_1$, $E_2$ are bounded subsets of $V$, and
let us check that $E_1 \cup E_2$ is bounded.  Let $W$ be any open subset
of $V$ such that $0 \in W$.  We may as well assume that $W$ is starlike
around the origin.  Since $E_1$, $E_2$ are bounded, there are positive
real numbers $t_1$, $t_2$ such that 
\begin{equation}
	E_1 \subseteq t_1 \, W, \quad E_2 \subseteq t_2 \, W.
\end{equation}
If we set $t = \max(t_1, t_2)$, then 
\begin{equation}
	E_1 \cup E_2 \subseteq t \, W,
\end{equation}
which is what we wanted.

	Next, if $E_1$, $E_2$ are nonempty bounded subsets of $V$,
then the sum $E_1 + E_2$ is a bounded subset of $V$ too.  To see this,
let $W$ be an open subset of $V$ such that $0 \in W$, and let $U_1$,
$U_2$ be open subsets of $V$ such that $0 \in U_1, U_2$, each of
$U_1$, $U_2$ is starlike around the origin, and $U_1, U_2 \subseteq
W$.  Because $E_1$, $E_2$ are bounded, there are positive real numbers
$t_1$, $t_2$ such that
\begin{equation}
	E_1 \subseteq t_1 \, U_1, \quad E_2 \subseteq t_2 \, U_2.
\end{equation}
The starlikeness conditions imply that if $t = \max(t_1, t_2)$, then
\begin{equation}
	E_1 \subseteq t \, U_1, \quad E_2 \subseteq t \, U_2.
\end{equation}
Hence
\begin{equation}
	E_1 + E_2 \subseteq t \, U_1 + t \, U_2 \subseteq t \, W,
\end{equation}
which is what we wanted.

	In a locally convex topological vector space, the convex hull
of a bounded set is also a bounded set.

	Suppose that $V$ is a topological vector space, and that $U$
is a nonempty open subset of $V$ which is bounded.  This is a strong
condition.  By translation we may assume that $0 \in U$.  One can
check that the family of sets of the form
\begin{equation}
	\frac{1}{n} \, U, \quad n \in {\bf Z}_+,
\end{equation}
forms a local basis for the topology of $V$ at $0$.  If there is a
bounded nonempty convex open subset of $V$, then it follows that there
is a norm on $V$ so that the topology defined by the norm is the same
as the topology that is given.

\section{Countability conditions and topological vector spaces}
\label{section on countability conditions and topological vector spaces}

	When a topological space satisfies the first axiom of
countability, it has the nice property that sequences are enough to
describe the topology, in the sense that closed subsets can be
characterized in terms of convergent sequences, for instance.  In the
context of topological vector spaces, we shall be concerned with two
basic scenarios, where the space itself satisfies the first axiom of
countability, and when suitably restricted subsets satisfy the first
axiom of countability as topological spaces themselves.  In this
section we shall concentrate on the first case, where a topological
vector space itself satisfies the first axiom of countability.

	Let $V$ be a topological vector space, and assume that $V$
satisfies the first axiom of countability at $0$, so that there is a
sequence of open sets in $V$ which contain $0$ and form a basis for
the local topology at $0$.  This implies that $V$ satisfies the
first axiom of countability at all of its points, since one can
translate the sequence of open sets around $0$ to any point in $V$.
In other words, if $\mathcal{B}_0$ is a local basis for the topology
of $V$ at $0$, then
\begin{equation}
	\mathcal{B}_v = \{v + U : U \in \mathcal{B}_0\}
\end{equation}
is a local basis for the topology of $V$ at $v \in V$.  If the
topology on $V$ can be defined by a norm, or if $V$ has a nonempty
bounded open subset, or if the topology on $V$ can be defined by
an at most countable family of seminorms, then $V$ satisfies the
first axiom of countability.

	A basic result states that when $V$ is a topological vector
space which satisfies the first axiom of countability, there is a
metric on $V$, which is in fact invariant under translations, such
that the topology determined by the metric is the same as the one on
$V$ already.  In general this metric may not come from a norm, or
something like that, however.  I like to not take this route, unless
there is a nice metric with some nice features, or something along
those lines, but instead pretty much do the same things directly, in
terms of the objects available, which may have quite nice properties.
Of course if there is a nice metric around, translation-invariant in
particular, say, then one would like to do things in a way which is
compatible with the metric.

	As a basic instance of this, let us say that a sequence
$\{z_j\}_{j=1}^\infty$ in $V$ is a \emph{Cauchy sequence}\index{Cauchy
sequence in a topological vector space} if for each open subset
$U$ of $V$ such that $0 \in U$ there is a positive integer $N$
such that
\begin{equation}
	z_j - z_l \in U \hbox{ for all } j, l \ge N.
\end{equation}
This definition makes sense in any topological vector space, but behaves
better when something like the first axiom of countability holds.
Observe that if $\{z_j\}_{j=1}^\infty$ converges to some point in $V$,
then $\{z_j\}_{j=1}^\infty$ is a Cauchy sequence.  Note that if $V$
is equipped with a translation-invariant metric which defines the
topology on $V$, then this notion of a Cauchy sequence and the usual
one defined using the metric are equivalent.

	A topological vector space which satisfies the first axiom of
countability and in which every Cauchy sequence converges is said to
be \emph{complete}.\index{complete topological vector space satisfying
the first axiom of countability} A normed vector space which is
complete is called a \emph{Banach space}.\index{Banach space} A vector
space equipped with a countable family of seminorms yielding a
complete topological vector space is called a \emph{Fr\'echet
space}.\index{Fr\'echet space}

	Let $V$ be a topological vector space which satisfies the
first axiom of countability at $0$, and let $\mathcal{B}_0$ be a
countable collection of open subsets of $V$ which contain $0$ and form
a local basis for the topology at $0$.  Suppose that $V$ is separable
as a topological space, so that there is a countable dense subset $E$
of $V$.  Consider the countable family of open subsets of $V$ given by
\begin{equation}
	\{x + U : x \in E, U \in \mathcal{B}_0\}.
\end{equation}
Let us show that this family is a basis for the topology of $V$,
and hence that $V$ satisfies the second axiom of countability.

	Let $y$ be any point in $V$ and let $W$ be an open subset of
$V$ such that $y \in W$.  We would like to find an $x \in E$ and a
$U \in \mathcal{B}_0$ such that
\begin{equation}
	y \in x + U, \quad x + U \subset W.
\end{equation}
Since $-y + W$ is an open subset of $V$ which contains $0$, there is a
$U \in \mathcal{B}_0$ so that $U + U \subseteq -y + W$, which is to
say that $y + U + U \subseteq W$.  Since $y + U$ is a nonempty open
subset of $V$, there is an $x \in E$ such that $x \in y + U$.
Hence 
\begin{equation}
	x + U \subseteq y + U + U \subseteq W,
\end{equation}
as desired.

\section{Continuous linear mappings}
\label{section on continuous linear mappings}

	Let $V_1$, $V_2$ be topological vector spaces, both real or
both complex.  Since $V_1$, $V_2$ are vector spaces, it makes sense to
talk about linear mappings from $V_1$ to $V_2$, and since $V_1$, $V_2$
are topological spaces it makes sense to talk about continuous
mappings from $V_1$ to $V_2$.  We can also talk about continuous
linear mappings,\index{continuous linear mappings between topological
vector spaces} which we shall consider now.  

	Suppose that $f$ is a linear mapping from $V_1$ to $V_2$.  To
say that $f$ is continuous at $0$ in $V_1$ means that for every open
subset $W$ of $V_2$ which contains $0$ there is an open subset $U$ of
$V_1$ which contains $0$ such that
\begin{equation}
	f(U) \subseteq W.
\end{equation}
It is easy to see that this implies that $f$ is continuous everywhere
on $V_1$, because of linearity.

	For any mapping $h$ from $V_1$ to $V_2$, we can define
uniform continuity to mean that for every open subset $W$ of $V_2$
such that $0 \in W$ there is an open subset $U$ of $V_1$ such that
$0 \in U$ and
\begin{equation}
	h(x) - h(y) \in W \ \hbox{ for all } \ x, y \in V_1 
						\hbox{ with } x - y \in U.
\end{equation}
This is analogous to the notion of uniform continuity for mappings between
metric spaces, and if $V_1$, $V_2$ are equipped with translation-invariant
metrics which define their topologies, then the two notions of uniform
continuity are equivalent.  A continuous linear mapping from $V_1$
to $V_2$ is automatically uniformly continuous.  Just as for continuous
mappings, notice that sums and scalar multiples of uniformly
continuous mappings are uniformly continuous.

	The vector space of continuous linear mappings from $V_1$ to
$V_2$ is denoted $\mathcal{CL}(V_1, V_2)$.  For a single topological
vector space $V$ we may write $\mathcal{CL}(V)$ instead of
$\mathcal{CL}(V, V)$.  Note that if $V_1$, $V_2$, $V_3$ are
topological vector spaces, all real or all complex, and if $f_1$ is a
continuous linear mapping from $V_1$ to $V_2$ and $f_2$ is a
continuous linear mapping from $V_2$ to $V_3$, then the composition
$f_2 \circ f_1$ is a continuous linear mapping from $V_1$ to $V_3$.
In particular, for a single topological vector space, composition
of elements of $\mathcal{CL}(V)$ again lie in $\mathcal{CL}(V)$.

	If $V$ is a topological vector space, then the
\emph{continuous dual}\index{continuous dual of a topological vector
space} of $V$ is denoted $V^*$ and defined to be $\mathcal{CL}(V, {\bf
R})$ when $V$ is a real vector space and to be $\mathcal{CL}(V, {\bf
C})$ when $V$ is a complex vector space.  Thus $V^*$ is a real vector
space when $V$ is a real vector space and $V^*$ is a complex vector
space when $V$ is a complex vector space.

	Let $V_1$, $V_2$ be topological vector spaces again, both real
or both complex.  A linear mapping $f$ from $V_1$ to $V_2$ is said
to be \emph{bounded}\index{bounded linear mappings between topological
vector spaces} if for every bounded subset $E$ of $V_1$ we have
that $f(E)$ is a bounded subset of $V_2$.  It is easy to see that
continuous linear mappings are bounded.

	The space of bounded linear mappings from $V_1$ to $V_2$ is
denoted $\mathcal{BL}(V_1, V_2)$.  As before, for a single topological
vector space $V$ we may write $\mathcal{BL}(V)$ rather than
$\mathcal{BL}(V, V)$.  Notice that sums and scalar multiples of
bounded linear mappings from $V_1$ to $V_2$ are again bounded, so that
$\mathcal{BL}(V_1, V_2)$ is a vector space in a natural way.  Since
continuous linear mappings are bounded, $\mathcal{CL}(V_1, V_2)$ is a
linear subspace of $\mathcal{BL}(V_1, V_2)$.

	If $V_1$, $V_2$, $V_3$ are topological vector spaces, all real
or all complex, and if $f_1$ is a bounded linear mapping from $V_1$ to
$V_2$ and $f_2$ is a bounded linear mapping from $V_2$ to $V_3$, then
the composition $f_2 \circ f_1$ is a bounded linear mapping from $V_1$
to $V_3$.  This is easy to see from the definitions.  In particular,
for a single topological vector space $V$, compositions of elements of
$\mathcal{BL}(V)$ are again elements of $\mathcal{BL}(V)$.

	If $V$ is a topological vector space, then the \emph{bounded
dual}\index{bounded dual of a topological vector space} is denoted
$V^b$ and defined to be $\mathcal{BL}(V, {\bf R})$ when $V$ is a real
vector space and to be $\mathcal{BL}(V, {\bf C})$ when $V$ is a
complex vector space.  Thus $V^b$ is a real vector space when $V$ is a
real vector spave and $V^b$ is a complex vector space when $V$ is a
complex vector space.  In both cases the continuous dual $V^*$ is a
linear subspace of the bounded dual $V^b$.

	If $V_1$, $V_2$ are topological vector spaces, both real or
both complex, and if $V_1$ contains a nonempty bounded open subset,
then it is easy to check that
\begin{equation}
	\mathcal{BL}(V_1, V_2) = \mathcal{CL}(V_1, V_2).
\end{equation}
In particular, this holds when $V_1$ is a normed vector space.  If
both $V_1$, $V_2$ are normed vector spaces, with norms $N_1$, $N_2$,
then the condition of boundedness for a linear mapping $f$ from $V_1$
to $V_2$ can be expressed succinctly by the statement that there is a
nonnegative real number $C$ such that
\begin{equation}
	N_2(f(v)) \le C \, N_1(v)
\end{equation}
for all $v \in V_1$.

	Suppose that $V_1$, $V_2$, and $V_3$ are vector spaces,
all real or all complex, equipped with norms $N_1$, $N_2$, $N_3$,
respectively.  Let $f_1$ be a bounded linear mapping from $V_1$
to $V_2$, and define the corresponding operator norm of $f_1$ by
\begin{equation}
	\|f_1\|_{op, 12} = \sup \{N_2(f_1(v)) : v \in V_1, N_1(v) \le 1\}.
\end{equation}
It is not difficult to show that this does define a norm on the vector
space $\mathcal{BL}(V_1, V_2)$.  Let $f_2$ be a bounded linear mapping
from $V_2$ to $V_3$, and define the corresponding operator norm
$\|f_2\|_{op, 23}$ analogously.  The operator norm $\|\cdot \|_{op,
13}$ for bounded linear mappings from $V_1$ to $V_3$ can be defined in
the same way, and one can check that
\begin{equation}
	\|f_2 \circ f_1\|_{op, 13} \le \|f_1\|_{op, 12} \, \|f_2\|_{op, 23}.
\end{equation}

	If $V_1$, $V_2$ are vector spaces, both real or both complex,
equipped with norms $N_1$, $N_2$, and if $V_2$ is complete with
respect to this norm, which is to say that it is a Banach space, then
$\mathcal{BL}(V_1, V_2)$, equipped with the operator norm described in
the previous paragraph is also complete, and hence defines a Banach
space.  For if $\{T_j\}_{j=1}^\infty$ is a Cauchy sequence of bounded
linear mappings from $V_1$ to $V_2$, then $\{T_j(v)\}_{j=1}^\infty$
is a Cauchy sequence of elements of $V_2$ for each $v \in V_1$,
and hence converges in $V_2$.  If we denote the limit as $T(v)$, then
one can check that $T(v)$ is a bounded linear mapping from $V_1$ to $V_2$,
and that $\{T_j\}_{j=1}^\infty$ converges to $T$ with respect to the
operator norm.

	What about analogous topologies on $\mathcal{BL}(V_1, V_2)$,
$\mathcal{CL}(V_1, V_2)$ for topological vector spaces $V_1$, $V_2$,
both real or both complex, more generally?  This can be a somewhat
tricky issue, and so we only mention a few basic points.

	If there is a nonempty bounded open subset of $V_1$, then
$\mathcal{BL}(V_1, V_2) = \mathcal{CL}(V_1, V_2)$, and one can follow
similar ideas as for the case of normed vector spaces.  Namely, the
operator norm topology in the case of normed vector spaces basically
looks at the size of a bounded linear operator on, say, the unit ball
in the domain $V_1$, which then limits the size on all of $V_1$ by
linearity.  If $V_1$ has a nonempty open subset $U$, which we may as
well assume contains $0$, then we can look at the size of a bounded
linear operator $T$ from $V_1$ to $V_2$ in terms of the size of $T$ on
$U$, noting that $T(U)$ will be a bounded subset of $V_2$.  This is
especially simple if $V_2$ also contains a nonempty bounded open
subset $W$, which we may assume contains $0$, because we can then look
at for which positive real numbers $r$ we have $T(U) \subseteq r \,
W$.  In general, we can look at the size of $T$ in terms of which
open subsets of $V_2$ contain $T(U)$, which dilations of them do, etc.

	In general one can try to measure the size of a bounded linear
operator $T$ from $V_1$ to $V_2$ in terms of the size of $T$ on
bounded subsets of $V_1$.  That is, if $E$ is a bounded subset of
$V_1$ and $W$ is an open subset of $V_2$ which contains $0$, one can
ask whether $T(E) \subseteq W$, look at the set of $r > 0$ such that
$T(E) \subseteq r \, W$, and so on.  However, one may need to consider
families of choices of $E$ and $W$.

\chapter{Examples, 1}
\label{chapter on examples, 1}

\section{$\ell^p$ spaces, $1 \le p \le \infty$, and $c_0$}
\label{section on ell^p spaces, 1 le p le infty, and c_0}

	If $p$ is a real number such that $1 \le p < \infty$,
define $\ell^p({\bf R})$, $\ell^p({\bf C})$, respectively, to
be the space of sequences $x = \{x_j\}_{j=1}^\infty$ with values
in ${\bf R}$, ${\bf C}$, respectively, such that
\begin{equation}
	\sum_{j=1}^\infty |x_j|^p < \infty.
\end{equation}
For $x = \{x_j\}_{j=1}^\infty$ in $\ell^p({\bf R})$, $\ell^p({\bf
C})$, we put
\begin{equation}
	\|x\|_p = \biggl(\sum_{j=1}^\infty |x_j|^p \biggr)^{1/p}.
\end{equation}
When $p = \infty$, we define $\ell^\infty({\bf R})$, $\ell^\infty({\bf
C})$, respectively, to be the space of sequences $x =
\{x_j\}_{j=1}^\infty$ with values in ${\bf R}$, ${\bf C}$,
respectively, which are bounded, i.e., whose values are contained
in a bounded set.  In this case we set
\begin{equation}
	\|x\|_\infty = \sup \{|x_j| : j \in {\bf Z}_+\}
\end{equation}
for $x = \{x_j\}_{j=1}^\infty$ in $\ell^\infty({\bf R})$ or
$\ell^\infty({\bf C})$.  When $p = \infty$ we can also consider the
spaces $c_0({\bf R})$, $c_0({\bf C})$, respectively, of sequences
$x = \{x_j\}_{j=1}^\infty$ with values in ${\bf R}$, ${\bf C}$,
respectively, such that
\begin{equation}
	\lim_{j \to \infty} x_j = 0.
\end{equation}
Such a sequence is bounded in particular, so that
\begin{equation}
	c_0({\bf R}) \subseteq \ell^\infty({\bf R}), \quad
		c_0({\bf C}) \subseteq \ell^\infty({\bf C}).
\end{equation}

	It is easy to see that $\ell^p({\bf R})$, $\ell^p({\bf C})$
are real and complex vector spaces for all $p$, $1 \le p \le \infty$,
and that $c_0({\bf R})$, $c_0({\bf C})$ are linear subspaces of
$\ell^\infty({\bf R})$, $\ell^\infty({\bf C})$, respectively.
Clearly $\|x\|_p \ge 0$ for all $x$ in $\ell^p({\bf R})$,
$\ell^p({\bf C})$, with $\|x\|_p = 0$ if and only if $x = 0$,
and 
\begin{equation}
	\|\alpha \, x\|_p = |\alpha| \, \|x\|_p
\end{equation}
when $\alpha$ is a real or complex number and $x$ is an element of
$\ell^p({\bf R})$ or $\ell^p({\bf C})$, respectively.  When $p = 1$
one can check directly from the definition that $\|x\|_p$ satisfies
the triangle inequality, and hence defines a norm on $\ell^p({\bf
R})$, $\ell^p({\bf C})$.  One way to see this when $1 < p < \infty$ is
to check that the set of $x$ in $\ell^p({\bf R})$ or in $\ell^p({\bf
C})$ such that $\|x\|_p \le 1$ is a convex set.  This is not difficult
to do, using the fact that $t^p$ is a convex function on the set
of nonnegative real numbers.

	Thus $\ell^p({\bf R})$, $\ell^p({\bf C})$ are normed vector
spaces for $1 \le p \le \infty$.  Also $c_0({\bf R})$, $c_0({\bf C})$
are normed vector spaces with respect to the norm $\|x\|_\infty$.
One can show that if $1 \le p \le q < \infty$, then
\begin{eqnarray}
\label{ell^p(R) subseteq ell^q(R) subseteq c_0(R) subseteq ell^infty(R)}
	&& \ell^p({\bf R}) \subseteq \ell^q({\bf R}) \subseteq c_0({\bf R})
				\subseteq \ell^\infty({\bf R}), 	\\
\label{ell^p(C) subseteq ell^q(C) subseteq c_0(C) subseteq ell^infty(C)}
	&& \ell^p({\bf C}) \subseteq \ell^q({\bf C}) \subseteq c_0({\bf C})
				\subseteq \ell^\infty({\bf C}),
\end{eqnarray}
and
\begin{equation}
	\|x\|_\infty \le \|x\|_q \le \|x\|_p
\end{equation}
when $x$ is an element of $\ell^p({\bf R})$, $\ell^p({\bf C})$.

	Notice that $\ell^p({\bf R})$, $\ell^p({\bf C})$ are complete,
in the sense that every Cauchy sequence converges.  Let us sketch the
argument.  Suppose that $\{x(l)\}_{l=1}^\infty$ is a Cauchy sequence
in one of these $\ell^p$ spaces, so that for each $l$, $x(l) =
\{x_j(l)\}_{j=1}^\infty$ is itself a sequence of real or complex
numbers.  Because of the way that the norm is defined, one can check
that $\{x_j(l)\}_{l=1}^\infty$ is a Cauchy sequence as a sequence of
real or complex numbers for each positive integer $j$, and hence
converges as a sequence of real or complex numbers to some $x_j$.
Thus we get a sequence $x = \{x_j\}_{j=1}^\infty$, and one can show
that this $x$ lies in the $\ell^p$ space and is the limit of the
sequence $\{x(l)\}_{l=1}^\infty$ in the $\ell^p$ space.

	One can also show that $c_0({\bf R})$, $c_0({\bf C})$ are
closed linear subspaces of $\ell^\infty({\bf R})$, $\ell^p({\bf C})$,
and hence are complete as topological vector spaces in their own
right.  This is not hard to do.

	For the $\ell^p$ spaces when $1 \le p < \infty$ and for the
$c_0$ spaces, the linear subspaces consisting of sequences $x =
\{x_j\}_{j=1}^\infty$ such that $x_j = 0$ for all but at most finitely
many $j$ are dense.  This is not the case for the $\ell^\infty$
spaces, for which the closure of this subspace is the corresponding
$c_0$ space.  As a result, the $\ell^p$ spaces for $1 \le p < \infty$
and the $c_0$ spaces are separable, which is to say that they contain
countable dense subsets, and therefore satisfy the second axiom of
countability.

	The inclusions (\ref{ell^p(R) subseteq ell^q(R) subseteq
c_0(R) subseteq ell^infty(R)}), (\ref{ell^p(C) subseteq ell^q(C)
subseteq c_0(C) subseteq ell^infty(C)}) can be viewed as bounded
linear mappings between the $\ell^p$ and $c_0$ spaces.  Another
natural class of bounded linear mappings from the $\ell^p$ or $c_0$
spaces into the same space again are given by multiplication by a
bounded sequence.  One can also map $\ell^\infty$ and $c_0$ spaces
into $\ell^p$ spaces by multiplication by a sequence in $\ell^p$, to
get a bounded linear mapping.  Next we describe a class of bounded
linear mappings from $\ell^p$ spaces to other $\ell^r$ spaces, using
H\"older's inequality, which we review first.

	Suppose that $p_0$, $q_0$ are real numbers such that
$1 < p_0, q_0 < \infty$ and
\begin{equation}
	\frac{1}{p_0} + \frac{1}{q_0} = 1.
\end{equation}
Let $\{x_j\}_{j=1}^\infty$, $\{y_j\}_{j=1}^\infty$ be sequences of
nonnegative real numbers such that
\begin{equation}
	\sum_{j=1}^\infty x_j^{p_0}, \ \sum_{j=1}^\infty y_j^{q_0} < \infty.
\end{equation}
H\"older's inequality states that
\begin{equation}
	\sum_{j=1}^\infty x_j \, y_j 
		\le \biggl(\sum_{k=1}^\infty x_k^{p_0}\biggr)^{1/p_0} \,
			\biggl(\sum_{l=1}^\infty y_l^{q_0}\biggr)^{1/q_0},
\end{equation}
and in particular that the sum on the left converges.  

	To see this,
one can start with the inequality
\begin{equation}
	a \, b \le \frac{a^{p_0}}{p_0} + \frac{b^{q_0}}{q_0}
\end{equation}
for nonnegative real numbers.  We leave this inequality as an exercise.
As a consequence,
\begin{equation}
	\sum_{j=1}^\infty x_j \, y_j
		\le \frac{1}{p_0} \sum_{k=1}^\infty x_k^{p_0}
			+ \frac{1}{q_0} \sum_{l=1}^\infty y_l^{q_0}.
\end{equation}
One can derive H\"older's inequality from this using homogeneity
considerations, which is to say by multiplying the $x_j$'s and $y_l$'s
by positive real numbers.

	As a consequence of H\"older's inequality, if $p$, $q$, and $r$
are positive real numbers such that
\begin{equation}
	\frac{1}{r} = \frac{1}{p} + \frac{1}{q},
\end{equation}
and if $\{u_j\}_{j=1}^\infty$, $\{w_j\}_{j=1}^\infty$ are sequences
of nonnegative real numbers such that
\begin{equation}
	\sum_{j=1}^\infty u_j^p, \ \sum_{j=1}^\infty w_j^q < \infty,
\end{equation}
then
\begin{equation}
	\biggl(\sum_{j=1}^\infty (u_j \, w_j)^r \biggr)^{1/r}
		\le \biggl(\sum_{k=1}^\infty u_j^p \biggr)^{1/p} \,
			\biggl(\sum_{l=1}^\infty w_j^q \biggr)^{1/q},
\end{equation}
and in particular the sum on the left converges.  This leads to
bounded linear mappings from $\ell^p$ to $\ell^r$ by multiplication by
sequences in $\ell^q$ when $p$, $q$, and $r$ are as above and $r \ge
1$.

	Let us mention two other types of bounded linear operators on
the $\ell^p$ and $c_0$ spaces, namely, backward and forward shift
operators.  If $\{x_j\}_{j=1}^\infty$ is a sequence of real or complex
numbers, then the backward shift operator sends this sequence to
$\{x_{j+1}\}_{j=1}^\infty$, which has the effect of dropping the first
term.  On each $\ell^p$ space, and on the $c_0$ spaces, this defines a
bounded linear operator from the space onto the same space, with
operator norm $1$, and with a $1$-dimensional kernel or nullspace
consisting of the sequences $\{x_j\}_{j=1}^\infty$ such that $x_j = 0$
when $j \ge 2$.  The forward shift operator takes a sequence
$\{x_j\}_{j=1}^\infty$ and sends it to the sequence with first term
equal to $0$ and with $j$th term equal to $x_{j-1}$ when $j \ge 2$.
This defines a bounded linear operator on each $\ell^p$ space and on
the $c_0$ spaces, which is an isometry, in the sense that the norm is
preserved, and which maps the space onto the subspace of sequences
whose first term is equal to $0$.

	One can also consider doubly-infinite sequences
$\{x_j\}_{j=-\infty}^\infty$ of real and complex numbers, and the
corresponding $\ell^p$ and $c_0$ Banach spaces.  For these spaces we
have the same kind of inclusions and multiplication operators as
before.  The backward and forward shift operators now define isometric
linear mappings of these spaces onto themselves.
	
	If $V$ is any topological vector space, one can define
$\ell^p(V)$ in a natural way as a generalization of $\ell^p({\bf R})$,
$\ell^p({\bf C})$.  This is simplest when $V$ is equipped with a norm
$N$.  In this case we can define $\ell^p(V)$ to be the vector space of
sequences $\{x_j\}_{j=1}^\infty$ with terms in $V$ such that the
sequence $\{N(x_j)\}_{j=1}^\infty$ of norms of the $x_j$'s lies in
$\ell^p({\bf R})$.  We can define the $\ell^p(V)$ norm of such a
sequence to be the $\ell^p({\bf R})$ norm of
$\{N(x_j)\}_{j=1}^\infty$, and indeed one can check that this does
define a norm on $\ell^p(V)$.  Similarly, we can define $c_0(V)$ as a
closed subspace of $\ell^\infty(V)$, consisting of the sequences
$\{x_j\}_{j=1}^\infty$ in $V$ such that $\{N(x_j)\}_{j=1}^\infty$ lies
in $c_0({\bf R})$.

	In general for a topological vector space $V$, we can take
$\ell^\infty(V)$ to be the vector space of sequences
$\{x_j\}_{j=1}^\infty$ in $V$ which are bounded, in the sense that
there is a bounded subset of $V$ which contains all the terms of the
sequence.  We can define $c_0(V)$ to be the space of sequences in $V$
which tend to $0$ in $V$, and $\ell^p(V)$ to be the space of sequences
$\{x_j\}_{j=1}^\infty$ such that for every open subset $U$ of $V$
which contains $0$ there is a sequence $\{r_j\}_{j=1}^\infty$ of
positive real numbers in $\ell^p({\bf R})$ such that $x_j \in r_j \,
U$ for all $j$.  Alternatively one might define $c_0(V)$, $\ell^p(V)$
to be the spaces of sequences $\{x_j\}_{j=1}^\infty$ in $V$ which can
be expressed as $\{r_j \, v_j\}_{j=1}^\infty$, where
$\{r_j\}_{j=1}^\infty$ is a sequence of positive real numbers in
$c_0({\bf R})$, $\ell^p({\bf R})$, respectively, and
$\{v_j\}_{j=1}^\infty$ is a bounded sequence in $V$.

\section{$\ell^p$ spaces, $0 < p < 1$}
\label{section on ell^p spaces, 0 < p < 1}

	Just as when $p \ge 1$, for a real number $p$ such that $0 < p
< 1$, let us define $\ell^p({\bf R})$, $\ell^p({\bf C})$,
respectively, to be the spaces of sequences $x = \{x_j\}_{j=1}^\infty$
with values in ${\bf R}$, ${\bf C}$, respectively, such that
\begin{equation}
	\sum_{j=1}^\infty |x_j|^p < \infty.
\end{equation}
When $x = \{x_j\}_{j=1}^\infty$ is an element of $\ell^p({\bf R})$ or
$\ell^p({\bf C})$, we again set
\begin{equation}
	\|x\|_p = \biggl(\sum_{j=1}^\infty |x_j|^p \biggr)^{1/p}.
\end{equation}

	It is easy to see that $\ell^p({\bf R})$, $\ell^p({\bf C})$
are real and complex vector spaces, respectively, when $0 < p < 1$.
Also $\|x\|_p \ge 0$ for all $x$ in $\ell^p({\bf R})$, $\ell^p({\bf C})$,
$\|x\|_p = 0$ if and only if $x = 0$, and 
\begin{equation}
	\|\alpha \, x\|_p = |\alpha| \, \|x\|_p
\end{equation}
for all real or complex numbers $\alpha$ and all $x$ in $\ell^p({\bf
R})$ or $\ell^p({\bf C})$, respectively.  It is not the case that
$\|x\|_p$ is a norm on $\ell^p({\bf R})$ or $\ell^p({\bf C})$ when 
$0 < p < 1$, however, because the triangle inequality does not hold.
There is a substitute for this, which is that
\begin{equation}
	\|x + y\|_p^p \le \|x\|_p^p + \|y\|_p^p
\end{equation}
for $x$, $y$ in $\ell^p({\bf R})$ or in $\ell^p({\bf C})$, $0 < p \le 1$.
One sometimes calls this the \emph{$p$-triangle inequality}, and
says that $\|x\|_p$ is a \emph{$p$-norm}.

	The $p$-triangle inequality implies that
\begin{equation}
	\|x - y\|_p^p
\end{equation}
defines a metric on $\ell^p({\bf R})$, $\ell^p({\bf C})$, and it is
not difficult to see that $\ell^p({\bf R})$, $\ell^p({\bf C})$
become topological vector spaces with respect to this topology.
Analogous to the situation for a normed vector space, a subset $E$
of $\ell^p({\bf R})$, $\ell^p({\bf C})$ is bounded in the sense for
topological vector spaces if and only if it is bounded in the sense
for metric spaces, which is to say that the set of nonnegative
real numbers 
\begin{equation}
	\{\|x\|_p : x \in E\}
\end{equation} is bounded from above.  The open unit ball
\begin{equation}
	\{x \in \ell^p({\bf R}) : \|x\|_p < 1\}
\end{equation}
in $\ell^p({\bf R})$ is a bounded open subset of $\ell^p({\bf R})$
which contains $0$, is starlike around $0$, and is symmetric, and
the open unit ball
\begin{equation}
	\{x \in \ell^p({\bf C}) : \|x\|_p < 1\}
\end{equation}
in $\ell^p({\bf C})$ is a bounded open subset of $\ell^p({\bf C})$
which contains $0$, is starlike around $0$, and is circular.

	When $0 < p < 1$, one can check that the convex hull of the
open unit ball in $\ell^p({\bf R})$, $\ell^p({\bf C})$ contains the
set of real or complex sequences $y = \{y_j\}_{j=1}^\infty$,
respectively, such that $y_j = 0$ for all but finitely many $j$
and 
\begin{equation}
	\|y\|_1 = \sum_{j=1}^\infty |y_j| < 1.
\end{equation}
As a result, the convex hull of the open unit ball in $\ell^p({\bf
R})$, $\ell^p({\bf C})$ is not bounded.  Thus $\ell^p({\bf R})$,
$\ell^p({\bf C})$ are not locally convex when $0 < p < 1$, and
in particular the topology on them does not come from a norm.

	Let us mention that for $\ell^p({\bf C})$ there are natural
\emph{pseudoconvexity}\index{pseudoconvexity} properties of the unit
ball, in terms of complex analysis.  The function $\|v\|_p$ enjoys
corresponding \emph{plurisubharmonicity}\index{plurisubharmonicity}
properties rather than convexity.  Compare with \cite{Boc-M, Krantz2,
Hormander3}.

	If $0 < p \le q \le 1$, then
\begin{equation}
\label{ell^p(R) subseteq ell^q(R) subseteq ell^1(R)}
	\ell^p({\bf R}) \subseteq \ell^q({\bf R}) \subseteq \ell^1({\bf R})
\end{equation}
and
\begin{equation}
\label{ell^p(C) subseteq ell^q(C) subseteq ell^1(C)}
	\ell^p({\bf C}) \subseteq \ell^q({\bf C}) \subseteq \ell^1({\bf C}),
\end{equation}
and for $x$ in $\ell^p({\bf R})$, $\ell^p({\bf C})$, we have that
\begin{equation}
	\|x\|_1 \le \|x\|_q \le \|x\|_p.
\end{equation}
Just as when $p \ge 1$, one can show that these $\ell^p$ spaces are
complete.  Also, the linear subspaces of the $\ell^p$ spaces
consisting of sequences $x = \{x_j\}_{j=1}^\infty$ such that $x_j = 0$
for all but at most finitely many $j$'s are dense, and as a result the
$\ell^p$ spaces are separable, and satisfy the second axiom of
countability, when $0 < p < 1$.

	The inclusions (\ref{ell^p(R) subseteq ell^q(R) subseteq
ell^1(R)}) and (\ref{ell^p(C) subseteq ell^q(C) subseteq ell^1(C)})
can be interpreted again as defining bounded linear operators between
the $\ell^p$ spaces.  Just as before, one has mappings between
$\ell^p$ spaces coming from multiplications.  One also has backward
and forward shift operators as before.

	One can also consider $\ell^p(V)$ for topological vector
spaces $V$, in essentially the same manner as when $p \ge 1$.

\section{Continuous functions on ${\bf R}^n$}
\label{section on continuous functions on R^n}

	Let $n$ be a positive integer, and let $\mathcal{C}({\bf R}^n, {\bf
R})$, $\mathcal{C}({\bf R}^n, {\bf C})$ denote the real and complex vector
spaces of real and complex-valued continuous functions on ${\bf R}^n$,
respectively.

	For each positive integer $j$, define $N_j(\cdot)$ on
$\mathcal{C}({\bf R}^n, {\bf R})$, $\mathcal{C}({\bf R}^n, {\bf C})$ by
\begin{equation}
	N_j(f) = \sup \{|f(x)| : x \in {\bf R}^n, |x| \le j\}
\end{equation}
when $f$ is a continuous real or complex-valued function on ${\bf
R}^n$.  It is easy to see that this defines a seminorm on
$\mathcal{C}({\bf R}^n, {\bf R})$, $\mathcal{C}({\bf R}^n, {\bf C})$,
and that $N_j(f) = 0$ for all $j \in {\bf Z}_+$ if and only if $f
\equiv 0$ on ${\bf R}^n$.  Thus $N_j$, $j \in {\bf Z}_+$, determine
topologies on $\mathcal{C}({\bf R}^n, {\bf R})$, $\mathcal{C}({\bf
R}^n, {\bf C})$, in such a way that they become real and complex
locally-convex topological spaces.

	If $E$ is a subset of $\mathcal{C}({\bf R}^n, {\bf R})$ or
$\mathcal{C}({\bf R}^n, {\bf C})$, then $E$ is bounded with respect to
the topology on the corresponding vector space just defined if
and only if for each $j$ the set of nonnegative real numbers
\begin{equation}
	\{N_j(f) : f \in E\}
\end{equation}
is a bounded set of real numbers.  Let us emphasize that the upper
bounds for these sets are permitted to depend on $j$.

	Using this description of the bounded subsets of
$\mathcal{C}({\bf R}^n, {\bf R})$, $\mathcal{C}({\bf R}^n, {\bf C})$,
it follows that these topological vector spaces do not contain
nonempty open subsets which are bounded.  In particular, there is no
single norm on $\mathcal{C}({\bf R}^n, {\bf R})$ or $\mathcal{C}({\bf
R}^n, {\bf C})$ which defines the same topology.

	Let $\{f_k\}_{k=1}^\infty$ be a sequence of functions in
$\mathcal{C}({\bf R}^n, {\bf R})$ or $\mathcal{C}({\bf R}^n, {\bf
C})$, and let $f$ be another function in the same space.  Then
$\{f_j\}_{j=1}^\infty$ converges to $f$ in the space if and only if
\begin{equation}
   \lim_{k \to \infty} \sup \{|f_k(x) - f(x)| : x \in {\bf R}^n, |x| \le j\}
		= 0
\end{equation}
for all positive integers $j$.  This is the same as saying that
$\{f_k\}_{k=1}^\infty$ converges to $f$ uniformly on every bounded
subset of ${\bf R}^n$.

	Now suppose that $\{f_k\}_{k=1}^\infty$ is a sequence of
functions in $\mathcal{C}({\bf R}^n, {\bf R})$ or $\mathcal{C}({\bf
R}^n, {\bf C})$ which is a Cauchy sequence.  Explicitly, this means
that for every $\epsilon > 0$ and every positive integer $j$ there is
a positive integer $L$ such that
\begin{equation}
	|f_k(x) - f_l(x)| \le \epsilon 
		\quad\hbox{when } x \in {\bf R}^n, |x| \le j, 
				\hbox{ and } k, l \ge L.
\end{equation}
In particular, $\{f_k(x)\}_{k=1}^\infty$ is a Cauchy sequence of real
or complex numbers for each $x \in {\bf R}^n$, which therefore
converges to a real or complex number, as appropriate, which we can
denote $f(x)$.  One can then show that $f$ is a continuous function on
${\bf R}^n$, and that $\{f_k\}_{k=1}^\infty$ converges to $f$
uniformly on bounded subsets of ${\bf R}^n$.  Thus $\mathcal{C}({\bf
R}^n, {\bf R})$, $\mathcal{C}({\bf R}^n, {\bf C})$ are complete,
so that they are in fact Fr\'echet spaces.

	Let us write $\mathcal{P}({\bf R}^n, {\bf R})$,
$\mathcal{P}({\bf R}^n, {\bf C})$ for the vector spaces of
real and complex-valued polynomials on ${\bf R}^n$, respectively.
Thus,
\begin{equation}
	\mathcal{P}({\bf R}^n, {\bf R}) 
		\subseteq \mathcal{C}({\bf R}^n, {\bf R}),
	\quad \mathcal{P}({\bf R}^n, {\bf C}) 
		\subseteq \mathcal{C}({\bf R}^n, {\bf C}).
\end{equation}
In fact, $\mathcal{P}({\bf R}^n, {\bf R})$, $\mathcal{P}({\bf R}^n,
{\bf C})$ are dense in $\mathcal{C}({\bf R}^n, {\bf R})$,
$\mathcal{C}({\bf R}^n, {\bf C})$, respectively.  Indeed, if $f$ is a
continuous function on ${\bf R}^n$, then for each positive integer $j$
there is a polynomial $P_j$ on ${\bf R}^n$ such that
\begin{equation}
	|f(x) - P_j(x)| \le \frac{1}{j}
		\quad \hbox{for all } x \in {\bf R}^n, |x| \le j,
\end{equation}
since continuous functions on compact subsets of ${\bf R}^n$ can be
approximated uniformly by polynomials.  Hence $\{P_j\}_{j=1}^\infty$
converges to $f$ uniformly on bounded subsets of ${\bf R}^n$, which
implies the above-mentioned denseness result.

	Using this, one can show that $\mathcal{C}({\bf R}^n, {\bf
R})$, $\mathcal{C}({\bf R}^n, {\bf C})$ contain countable dense
subsets, so that they are separable.  This implies that they satisfy
the second axiom of countability.  This might seem a bit surprising,
and one might say that these spaces and the topologies on them balance
each other in a nice way.

	If $f(x)$ is a real or complex-valued function on ${\bf R}^n$,
then the support of $f$, denoted $\supp f$, is defined to be the
closure of the set of points $x \in {\bf R}^n$ such that $f(x) \ne 0$.
Let us write $\mathcal{C}_{00}({\bf R}^n, {\bf R})$, $\mathcal{C}({\bf
R}^n, {\bf C})$ for the real and complex vector spaces of real and
complex-valued continuous functions on ${\bf R}^n$ with compact
support.  Of course 
\begin{equation}
	\mathcal{C}_{00}({\bf R}^n, {\bf R}) 
		\subseteq \mathcal{C}({\bf R}^n, {\bf R}), \quad
	\mathcal{C}_{00}({\bf R}^n, {\bf C}) 
		\subseteq \mathcal{C}({\bf R}^n, {\bf C}).
\end{equation}

	For each positive integer $l$, let $\phi_l(x)$ be a continuous
real-valued function on ${\bf R}^n$ such that $\phi_l(x) = 1$ when
$|x| \le l$, $\phi_l(x) = 0$ when $|x| \ge l + 1$, and $0 \le
\phi_l(x) \le 1$ for all $x \in {\bf R}^n$.  If $f$ is a continuous
real or complex-valued function on ${\bf R}^n$, then the sequence of
products $\{\phi_l \, f\}_{l=1}^\infty$ is a sequence of continuous
functions with compact support on ${\bf R}^n$ which converges to $f$
uniformly on every bounded subset of ${\bf R}^n$.  Thus
$\mathcal{C}_{00}({\bf R}^n, {\bf R})$, $\mathcal{C}_{00}({\bf R}^n,
{\bf C})$ are dense subspaces of $\mathcal{C}({\bf R}^n, {\bf R})$,
$\mathcal{C}({\bf R}^n, {\bf C})$, respectively, and this can be used
to give an alternate approach to the separability of the latter
spaces.

	For each nonempty compact subset $K$ of ${\bf R}^n$, define
$\mathcal{C}(K, {\bf R})$, $\mathcal{C}(K, {\bf C})$ to be the real
and complex vector spaces of real and complex-valued continuous
functions on $K$, respectively.  On these spaces we have the supremum
norm $\|f\|_{sup}$ defined by
\begin{equation}
	\|f\|_{sup} = \sup \{|f(x)| : x \in K \}.
\end{equation}
It is easy to see that the supremum norm does indeed define a norm
on $\mathcal{C}(K, {\bf R})$, $\mathcal{C}(K, {\bf C})$, and it
is well known that $\mathcal{C}(K, {\bf R})$, $\mathcal{C}(K, {\bf C})$
are complete with respect to the supremum norm, and therefore become
Banach spaces.  Of course convergence of a sequence of functions on $K$
in the supremum norm is the same as uniform convergence.

	For each nonempty compact subset $K$ of ${\bf R}^n$ we get
linear mappings
\begin{equation}
\label{{C}(R^n, R) to {C}(K, R), {C}(R^n, C) to {C}(R^n, C)}
	\mathcal{C}({\bf R}^n, {\bf R}) \to \mathcal{C}(K, {\bf R}),
		\quad
	   \mathcal{C}({\bf R}^n, {\bf C}) \to \mathcal{C}({\bf R}^n, {\bf C})
\end{equation}
defined by taking a continuous function $f$ on ${\bf R}^n$ and simply
restricting it to a continuous function on $K$.  There are well-known
extension results which say that every continuous function on $K$ can
be realized as the restriction to $K$ of a continuous function on
${\bf R}^n$, so that the mappings in (\ref{{C}(R^n, R) to {C}(K, R),
{C}(R^n, C) to {C}(R^n, C)}) are surjections.  In fact there are
continuous linear extension operators
\begin{equation}
	\mathcal{C}(K, {\bf R}) \to \mathcal{C}({\bf R}^n, {\bf R}),
		\quad
	   \mathcal{C}(K, {\bf C}) \to \mathcal{C}({\bf R}^n, {\bf C})
\end{equation}
such that the compositions of these mappings with the restriction mappings
are equal to the identity mapping on $\mathcal{C}(K, {\bf R})$,
$\mathcal{C}(K, {\bf C})$, respectively.

	Let us be a bit more precise.  Let $K_1$ be a compact subset
of ${\bf R}^n$ such that $K$ is contained in the interior of $K_1$.
The linear extension operators just mentioned can be chosen so that
continuous functions on $K$ are extended to continuous functions on
${\bf R}^n$ which are supported in $K_1$.  This can be obtained as a
by-product of standard constructions of the extension operator, or arranged
afterwards simply by multiplying by a continuous function on ${\bf R}^n$
which is equal to $1$ on $K$ and has support contained in $K_1$.
It is easy to have the extension operators also have norm equal to $1$
with respect to the supremum metric on the domain and range.

	Notice that for each $a \in {\bf R}^n$ the translation operator
\begin{equation}
	f(x) \mapsto f(x - a)
\end{equation}
defines a continuous linear mapping from each of $\mathcal{C}({\bf
R}^n, {\bf R})$, $\mathcal{C}({\bf R}^n, {\bf C})$ to itself.
More generally, if $\rho$ is a continuous mapping from ${\bf R}^n$
to itself, then 
\begin{equation}
	f(x) \mapsto f(\rho(x))
\end{equation}
defines a continuous linear mapping from each of $\mathcal{C}({\bf
R}^n, {\bf R})$, $\mathcal{C}({\bf R}^n, {\bf C})$ to itself.  As
another class of operators, if $h(x)$ is a continuous real or
complex-valued function on ${\bf R}^n$, respectively, then
\begin{equation}
	f(x) \mapsto h(x) \, f(x)
\end{equation}
defines a continuous linear mapping from $\mathcal{C}({\bf R}^n, 
{\bf R})$, $\mathcal{C}({\bf R}^n, {\bf C})$ to itself, respectively.

	One can also consider continuous functions with values
in a topological vector space $V$.

\section{Rapidly decreasing continuous functions on ${\bf R}^n$}
\label{section on rapidly decreasing continuous functions on R^n}

	Let $f(x)$ be a real or complex-valued function on ${\bf
R}^n$.  We say that $f(x)$ is \emph{rapidly decreasing}\index{rapidly
decreasing functions on ${\bf R}^n$} on ${\bf R}^n$ if for every
positive integer $j$ there is a positive real number $C(j)$ such that
\begin{equation}
	|f(x)| \le C(j) \, (|x| + 1)^j   \quad\hbox{for all } x \in {\bf R}^n.
\end{equation}
In other words, $f(x) = O((|x| + 1)^{-j})$ for all positive integers
$j$.  The space of real-valued continuous rapidly decreasing functions
on ${\bf R}^n$ is denoted $\mathcal{RC}({\bf R}^n, {\bf R})$, and the
space of complex-valued rapidly decreasing continuous functions on
${\bf R}^n$ is denoted $\mathcal{RC}({\bf R}^n, {\bf C})$.

	Clearly $\mathcal{RC}({\bf R}^n, {\bf R})$, $\mathcal{RC}({\bf
R}^n, {\bf C})$ are real and complex vector spaces, with respect to
ordinary addition and scalar multiplication of functions.  On each of
these spaces and for each positive integer $j$ we can define the norm
\begin{equation}
	M_j(f) = \sup \{|f(x)| \, (|x| + 1)^j : x \in {\bf R}^n\},
\end{equation}
and this family of norms leads to a topology on each of
$\mathcal{RC}({\bf R}^n, {\bf R})$, $\mathcal{RC}({\bf R}^n, {\bf C})$
which makes these spaces locally convex topological vector spaces.

	Let $\mathcal{C}_{00}({\bf R}^n, {\bf R})$,
$\mathcal{C}_{00}({\bf R}^n, {\bf C})$ be the real and complex vector
spaces of real and complex-valued continuous functions on ${\bf R}^n$,
respectively, as in the previous section.  Also, for each positive
integer $l$, let $\phi_l(x)$ be a continuous real-valued function on
${\bf R}^n$ such that $\phi_l(x) = 1$ when $|x| \le l$, $\phi_l(x) =
0$ when $|x| \ge l + 1$, and $0 \le \phi_l(x) \le 1$ for all $x \in
{\bf R}^n$, again as in the previous section.  If $f(x)$ is a rapidly
decreasing continuous function on ${\bf R}^n$, then one can check that
the sequence of products $\{\phi_l \, f\}_{l=1}^\infty$ converges to
$f$ in $\mathcal{RC}({\bf R}^n, {\bf R})$ or $\mathcal{RC}({\bf R}^n,
{\bf C})$, as appropriate

	Thus $\mathcal{C}_{00}({\bf R}^n, {\bf R})$,
$\mathcal{C}_{00}({\bf R}^n, {\bf C})$ are dense linear subspaces of
$\mathcal{RC}({\bf R}^n, {\bf R})$, $\mathcal{RC}({\bf R}^n, {\bf
C})$, respectively.  One can use this observation to show that the
topological vector spaces $\mathcal{RC}({\bf R}^n, {\bf R})$,
$\mathcal{RC}({\bf R}^n, {\bf C})$ are separable.  As a result, they
satisfy the second axiom of countability.

	Suppose that $\{f_k\}_{k=1}^\infty$ is a sequence of real or
complex-valued rapidly decreasing sequence of functions on ${\bf R}^n$
which is a Cauchy sequence with respect to the family of norms $M_j$,
$j \in {\bf Z}_+$.  This means that for each $\epsilon > 0$ and 
each positive integer $j$ there is a positive integer $L$ such that
\begin{equation}
	\sup \{|f_k(x) - f_l(x)| \, (|x| + 1)^j : x \in {\bf R}^n\} 
		\le \epsilon \quad\hbox{ for all } k, l \ge L.
\end{equation}
For each $x \in {\bf R}^n$, the sequence $\{f_k(x)\}_{k=1}^\infty$ of
real or complex numbers is then a Cauchy sequence as well, and hence
converges.  If we denote the limit $f(x)$, then one can show that
$f(x)$ is a rapidly decreasing continuous function on ${\bf R}^n$ and
that $\{f_k\}_{k=1}^\infty$ converges to $f$ with respect to the norms
$M_j$ for each positive integer $j$.  Thus $\mathcal{RC}({\bf R}^n,
{\bf R})$, $\mathcal{RC}({\bf R}^n, {\bf C})$ are complete, and hence
are Fr\'echet spaces.

	Let $E$ be a subset of $\mathcal{RC}({\bf R}^n, {\bf R})$
or $\mathcal{RC}({\bf R}^n, {\bf C})$.  Then $E$ is bounded in the
sense of bounded subsets of topological vector spaces if and only
if for every positive integer $j$ the set of nonnegative real numbers
\begin{equation}
	\{M_j(f) : f \in E\}
\end{equation}
is bounded from above.  As a result, one can verify that nonempty open
subsets of $\mathcal{RC}({\bf R}^n, {\bf R})$ and $\mathcal{RC}({\bf
R}^n, {\bf C})$ are not bounded.

	The family of norms $M_j$, $j \in {\bf Z}_+$, has a nice
feature, which basically says that the property of boundedness for a
subset of $\mathcal{RC}({\bf R}^n, {\bf R})$, $\mathcal{RC}({\bf R}^n,
{\bf C})$.  Namely, if $\{f_k\}_{k=1}^\infty$ is a bounded sequence in
$\mathcal{RC}({\bf R}^n, {\bf R})$ or $\mathcal{RC}({\bf R}^n, {\bf
R})$, so that $\{M_j(f_k)\}_{k=1}^\infty$ is a bounded sequence of
real numbers for each $j$, and if $\{f_k\}_{k=1}^\infty$ converges
uniformly on compact subsets of ${\bf R}^n$ to a real or
complex-valued function $f$, as appropriate, then $f$ is rapidly
decreasing and $\{f_k\}_{k=1}^\infty$ converges to $f$ in the topology
of rapidly decreasing functions.  In other words, if $E$ is a bounded
subset of $\mathcal{RC}({\bf R}^n, {\bf R})$ or $\mathcal{RC}({\bf
R}^n, {\bf C})$, then the topology on $E$ inherited from the space of
rapidly decreasing functions is the same as the topology on $E$
inherited from $\mathcal{C}({\bf R}^n, {\bf R})$, $\mathcal{C}({\bf
R}^n, {\bf C})$, as appropriate.

	Let $K$ be a nonempty compact subset of ${\bf R}^n$, and
consider the restriction mappings
\begin{equation}
	\mathcal{RC}({\bf R}^n, {\bf R}) \to \mathcal{C}(K, {\bf R}),
		\quad
	\mathcal{RC}({\bf R}^n, {\bf C}) \to \mathcal{C}(K, {\bf C}),
\end{equation}
i.e., the linear mappings which take a rapidly-decreasing continuous
function $f$ on ${\bf R}^n$ and restrict it to $K$.  As in the
previous section, $\mathcal{C}(K, {\bf R})$, $\mathcal{C}(K, {\bf C})$
become normed vector spaces in a natural way by using the supremum
metric.  It is easy to see that the restriction mappings above are
continuous mappings, and in fact they map $\mathcal{RC}({\bf R}^n,
{\bf R})$, $\mathcal{RC}({\bf R}^n, {\bf C})$ onto $\mathcal{C}(K, {\bf
R})$, $\mathcal{C}(K, {\bf C})$.  As in the previous section again,
there are even bounded linear extension operators from continuous
functions on $K$ to continuous functions on ${\bf R}^n$ with support
contained in a compact subset $K_1$ of ${\bf R}^n$ such that $K$ is
contained in the interior of $K$, so that the extension operator from
functions on $K$ to functions on ${\bf R}^n$ composed with the
restriction operator from functions on ${\bf R}^n$ to functions on $K$
is the identity operator on functions on $K$.

	If $h(x)$ is a continuous real or complex-valued function
of moderate growth, in the sense that there is a positive integer $l$
and a nonnegative real number $C$ such that
\begin{equation}
	|h(x)| \le C \, (1 + |x|)^l
\end{equation}
for all $x \in {\bf R}^n$, then the multiplication operator
\begin{equation}
	f(x) \mapsto h(x) \, f(x)
\end{equation}
defines a continuous linear mapping from $\mathcal{RC}({\bf R}^n, {\bf
R})$, $\mathcal{RC}({\bf R}^n, {\bf C})$ to itself, respectively.
This is not hard to check from the definitions.  If $\phi(x)$ is a
real or complex-valued continuous function on ${\bf R}^n$ with compact
support, then the multiplication operator
\begin{equation}
	f(x) \mapsto \phi(x) \, f(x)
\end{equation}
defines a continuous linear mapping from $\mathcal{C}({\bf R}^n, {\bf
R})$, $\mathcal{C}({\bf R}^n, {\bf C})$ to $\mathcal{RC}({\bf R}^n,
{\bf R})$, $\mathcal{RC}({\bf R}^n, {\bf C})$, respectively.  Of
course the range of this mapping is contained in the space of real or
complex-valued continuous functions on ${\bf R}^n$, respectively, with
support contained in the support of $\phi$.  Also, the obvious
inclusions of $\mathcal{RC}({\bf R}^n, {\bf R})$, $\mathcal{RC}({\bf
R}^n, {\bf C})$ into $\mathcal{C}({\bf R}^n, {\bf R})$,
$\mathcal{C}({\bf R}^n, {\bf C})$, respectively, are continuous linear
mappings.

	If $a \in {\bf R}^n$, then the translation operator
\begin{equation}
	f(x) \mapsto f(x - a)
\end{equation}
defines a continuous linear operator on each of $\mathcal{RC}({\bf
R}^n, {\bf R})$, $\mathcal{RC}({\bf R}^n, {\bf C})$.  If $A$
is an invertible linear mapping of ${\bf R}^n$ onto itself,
then the composition operator
\begin{equation}
	f(x) \mapsto f(A(x))
\end{equation}
defines a continuous linear operator on each of $\mathcal{RC}({\bf
R}^n, {\bf R})$, $\mathcal{RC}({\bf R}^n, {\bf C})$.  Compositions
with continuous mappings from ${\bf R}^n$ to ${\bf R}^n$ do not
quite work in general, but they do work under reasonable conditions.

\section{Normed and topological algebras}
\label{section on normed and topological algebras}

	Let $A$ be a vector space.  Suppose that in addition to the
vector space operations there is another binary operation on $A$ which
one can call a ``product'' which satisfies the usual associative law,
the usual dsitributive laws with respect to addition on $A$, and the
usual compatibility conditions with respect to scalar multiplication,
so that scalar multiplication commutes with multiplication on $A$.
Then we say that $A$ is an \emph{algebra}.\index{algebra} For
instance, if $V$ is a vector space, then the vector space
$\mathcal{L}(V)$ of linear mappings from $V$ into itself is an
algebra, using composition of linear operators as multiplication.

	Suppose in addition that $A$ is a topological vector space,
and that the operation of multiplication is a continuous mapping from
$A$ to $A$.  Then $A$ is said to be a \emph{topological
algebra}.\index{topological algebra}\index{algebra!topological} If the
topology on $A$ is determined by a norm $\|\cdot\|$, and if
\begin{equation}
	\|a \, b\| \le \|a\| \, \|b\|
\end{equation}
for all $a, b \in A$, i.e., the norm of a product is less than or
equal to the product of the corresponding norms, then we say that $A$
is a \emph{normed algebra}.\index{normed
algebra}\index{algebra!normed} This condition implies that
multiplication is continuous on $A$.

	A number of examples of these notions occur in this chapter.
The algebra of bounded linear operators on a normed vector space is a
normed algebra, with respect to the corresponding operator norm.  Even
if there is not a norm available, there are various ways in which
algebras of bounded or continuous linear operators on a topological
vector space can be topological algebras.

\chapter{Examples, 2}
\label{chapter on examples, 2}

	If $f(x)$ is a real or complex-valued continuous function on a
topological space, like ${\bf R}^n$, then the
\emph{support}\index{support of a real or complex-valued function} of
$f$, denoted $\supp f$, is defined to be the closure of the set of $x
\in {\bf R}^n$ such that $f(x) \ne 0$.

\section{Continuous functions with compact support on ${\bf R}^n$}
\label{section on continuous functions with compact support on R^n}

	As before, we write $\mathcal{C}_{00}({\bf R}^n, {\bf R})$,
$\mathcal{C}_{00}({\bf R}^n, {\bf C})$ for the spaces of real and
complex-valued continuous functions on ${\bf R}^n$ with compact
support, respectively.  Let us also write $\mathcal{C}_0({\bf R}^n,
{\bf R})$, $\mathcal{C}_0({\bf R}^n, {\bf C})$ for the spaces of real
and complex-valued continuous functions $f(x)$ on ${\bf R}^n$ which
``vanish at infinity'', in the sense that for every $\epsilon > 0$
there is a compact subset $K$ of ${\bf R}^n$ such that
\begin{equation}
	|f(x)| < \epsilon \ \hbox{ for all } x \in {\bf R}^n \backslash K.
\end{equation}
These are all clearly real and complex vector spaces with respect to
the usual operations of addition and scalar multiplication of
functions, and
\begin{equation}
	\mathcal{C}_{00}({\bf R}^n, {\bf R}) 
			\subseteq \mathcal{C}_0({\bf R}^n, {\bf R}), \quad 
		\mathcal{C}_{00}({\bf R}^n, {\bf R}) 
			\subseteq \mathcal{C}_0({\bf R}^n, {\bf C}).
\end{equation}
Functions in $\mathcal{C}_0({\bf R}^n, {\bf R})$, $\mathcal{C}_0({\bf
R}^n, {\bf C})$ are bounded in particular, and thus we may define
their supremum norms as usual by
\begin{equation}
	\|f\|_{sup} = \sup \{|f(x)| : x \in {\bf R}^n\}.
\end{equation}

	The spaces $\mathcal{C}_0({\bf R}^n, {\bf R})$,
$\mathcal{C}_0({\bf R}^n, {\bf C})$ are complete with respect to the
supremum norm.  In other words, suppose that $\{f_j\}_{j=1}^\infty$ is
a Cauchy sequence in one of these spaces, so that for every $\epsilon
> 0$ there is a positive integer $L$ such that
\begin{equation}
	\|f_j - f_l\|_{sup} < \epsilon \ \hbox{ for all } j, l \ge L.
\end{equation}
This implies that for each $x \in {\bf R}^n$ the sequence
$\{f_j(x)\}_{j=1}^\infty$ is a Cauchy sequence of real or complex
numbers, and hence converges.  If we denote the limit $f(x)$, then one
can show that $f(x)$ is continuous function on ${\bf R}^n$ which
vanishes at infinity, and that $\{f_j\}_{j=1}^\infty$ converges to $f$
uniformly on ${\bf R}^n$.  Thus $\{f_j\}_{j=1}^\infty$ converges to
$f$ in $\mathcal{C}_0({\bf R}^n, {\bf R})$ or $\mathcal{C}_0({\bf
R}^n, {\bf C})$, as appropriate.

	The linear subspaces $\mathcal{C}_{00}({\bf R}^n, {\bf R})$,
$\mathcal{C}_{00}({\bf R}^n, {\bf C})$ are dense in
$\mathcal{C}_0({\bf R}^n, {\bf R})$, $\mathcal{C}({\bf R}^n, {\bf
C})$, respectively.  To see this, for each positive integer $j$,
choose a real-valued continuous function $\phi_j(x)$ on ${\bf R}^n$
such that $\phi_j(x) = 1$ when $|x| \le j$, $\phi_j(x) = 0$ when $|x|
\ge j + 1$, and $0 \le \phi_j(x) \le 1$ for all $x \in {\bf R}^n$.  If
$f(x)$ is a function in $\mathcal{C}_0({\bf R}^n, {\bf R})$ or
$\mathcal{C}({\bf R}^n, {\bf C})$, then the sequence $\{\phi_j \,
f\}_{j=1}^\infty$ converges to $f$ in the supremum norm, and of course
$\phi_j(x) \, f(x)$ is a continuous function on ${\bf R}^n$ with
compact support for each $j$ too.

	On the other hand, we can be interested in some kind of
topological structure on $\mathcal{C}_{00}({\bf R}^n, {\bf R})$,
$\mathcal{C}_{00}({\bf R}^n, {\bf C})$ so that they are already
complete, i.e., a ``finer'' topological structure.  Here are three
basic features of such a structure.  First, a subset $E$ of
$\mathcal{C}_{00}({\bf R}^n, {\bf R})$, $\mathcal{C}_{00}({\bf R}^n,
{\bf C})$ would be bounded if and only if there is a compact subset
$K$ of ${\bf R}^n$ such that each function $f$ in $E$ has support
contained in $E$, and the collection of supremum norms of functions
$f$ in $E$ is a bounded set of real numbers.  Second, a sequence of
functions $\{f_j\}_{j=1}^\infty$ in $\mathcal{C}_{00}({\bf R}^n, {\bf
R})$, $\mathcal{C}_{00}({\bf R}^n, {\bf C})$ would be considered to
converge to a function $f$ in the same space if there is a compact
subset $K$ of ${\bf R}^n$ such that the support of each $f_j$ and o
$f$ is contained in $K$, and if $\{f_j\}_{j=1}^\infty$ converges to
$f$ uniformly.  Third, a sequence $\{f_j\}_{j=1}^\infty$ in
$\mathcal{C}_{00}({\bf R}^n, {\bf R})$, $\mathcal{C}_{00}({\bf R}^n,
{\bf C})$ would be considered a Cauchy sequence if there is a compact
subset $K$ of ${\bf R}^n$ such that the support of each $f_j$ is
contained in $K$ and if $\{f_j\}_{j=1}^\infty$ satisfies the usual
Cauchy sequence condition with respect to the supremum norm, in which
case it follows that $\{f_j\}_{j=1}^\infty$ converges in the sense
mentioned before.

	In fact, the idea of bounded or continuous linear functionals
on $\mathcal{C}_{00}({\bf R}^n, {\bf R})$, $\mathcal{C}_{00}({\bf
R}^n, {\bf C})$ behaves well and is quite interesting.  Namely, a
linear functional $\lambda$ on one of these spaces is considered to
be bounded if for each compact subset $K$ of ${\bf R}^n$ there is
a nonnegative real number $C_K$ such that
\begin{equation}
	|\lambda(f)| \le C_K \, \|f\|_{sup}
\end{equation}
whenever $f$ is a continuous function with support in $K$.  We consider
$\lambda$ to be continuous if 
\begin{equation}
	\lim_{j \to \infty} \lambda(f_j) = \lambda(f)
\end{equation}
whenever $\{f_j\}_{j=1}^\infty$ converges to $f$ in the sense
described in the previous paragraph.  This is equivalent to saying
that for each compact subset $K$ of ${\bf R}^n$, $\lambda$ is
continuous on the normed vector space of continuous functions with
support contained in $K$ with respect to the supremum norm, and
it is also equivalent to the boundedness condition above.

	Suppose that $\lambda$ is a bounded linear functional on
$\mathcal{C}_{00}({\bf R}^n, {\bf R})$ or $\mathcal{C}_{00}({\bf R}^n,
{\bf C})$.  For each nonempty compact subset $K$ of ${\bf R}^n$,
define $N_K(\lambda)$ to be the supremum of $|\lambda(f)|$
over all real or complex-valued continuous functions $f$ on ${\bf R}^n$,
as appropriate, with 
\begin{equation}
	\supp f \subseteq K \ \hbox{ and } \ \|f\|_{sup} \le 1.
\end{equation}
The boundedness of $\lambda$ implies that $N_K(\lambda)$ is finite for
all compact subsets $K$ of ${\bf R}^n$.  This gives a nice family of
seminorms on the duals of $\mathcal{C}_{00}({\bf R}^n, {\bf R})$,
$\mathcal{C}_{00}({\bf R}^n, {\bf C})$, and in fact the countable
family $N_{\overline{B}(0,j)}$, $j \in {\bf Z}_+$ is sufficient, for
defining topologies on the duals of $\mathcal{C}_{00}({\bf R}^n, {\bf
R})$, $\mathcal{C}_{00}({\bf R}^n, {\bf C})$.  It is not difficult to
see that these dual spaces are complete, and hence are Fr\'echet
spaces.

	There is a natural way in which to define topologies on
$\mathcal{C}_{00}({\bf R}^n, {\bf R})$, $\mathcal{C}_{00}({\bf R}^n,
{\bf C})$ so that they become locally convex topological vector spaces
with the kind of properties described above.  These are basic examples
of \emph{inductive limit spaces}.\index{inductive limit spaces} This
is a bit tricky, and we shall not pursue this here.

\section{Bounded continuous functions on ${\bf R}^n$}
\label{section on bounded continuous functions on R^n}

	Let $\mathcal{C}_b({\bf R}^n, {\bf R})$, $\mathcal{C}_b({\bf
R}^n, {\bf C})$ denote the real and complex vector spaces of real and
complex-valued bounded continuous functions on ${\bf R}^n$,
respectively.  Because the functions are assumed to be bounded, the
supremum norm is defined on $\mathcal{C}_b({\bf R}^n, {\bf R})$,
$\mathcal{C}_b({\bf R}^n, {\bf C})$ as before, so that these spaces
become normed vector spaces.  Also, $\mathcal{C}_b({\bf R}^n, {\bf
R})$, $\mathcal{C}_b({\bf R}^n, {\bf C})$ are complete with respect to
the supremum norm, so that they are Banach spaces.  

	Of course
\begin{equation}
	\mathcal{C}_{00}({\bf R}^n, {\bf R}) 
		\subseteq \mathcal{C}_0({\bf R}^n, {\bf R})
		\subseteq \mathcal{C}_b({\bf R}^n, {\bf R})
		\subseteq \mathcal{C}({\bf R}^n, {\bf R})
\end{equation}
and
\begin{equation}
	\mathcal{C}_{00}({\bf R}^n, {\bf C})
		\subseteq \mathcal{C}_0({\bf R}^n, {\bf C})
		\subseteq \mathcal{C}_b({\bf R}^n, {\bf C})
		\subseteq \mathcal{C}({\bf R}^n, {\bf C}),
\end{equation}
and these inclusions are all continuous.  For each nonempty compact
subset $K$ of ${\bf R}^n$, we also have the continuous linear mapping
of restriction of continuous functions on ${\bf R}^n$ to continuous
functions on $K$, and we have mentioned before that for each compact
subset $K_1$ of ${\bf R}^n$ such that $K$ is contained in the interior
of $K_1$ there are continuous linear extension operators from
continuous functions on $K$ to continuous functions on ${\bf R}^n$
with support contained in $K_1$.  Actually, we have also discussed
rapidly decreasing continuous functions on ${\bf R}^n$, which include
continuous functions with compact support and which are included among
continuous functions which tend to $0$ at infinity in ${\bf R}^n$, and
in the next section we shall discuss continuous functions of
polynomial growth, which include bounded continuous functions and are
included in the space of continuous functions in general.

	Let us consider another way to think about topology and
so forth on the spaces of real and complex-valued bounded continuous
functions on ${\bf R}^n$.  Let $E$ be a collection of such functions,
and let us assume that $E$ is bounded in the usual sense, so that
\begin{equation}
	\{\|f\|_{sup} : f \in E \}
\end{equation}
is a bounded set of real numbers.  On this set, let us use the topology
from the space of continuous functions on ${\bf R}^n$ in general,
rather than the topology induced by the supremum norm.

	Notice that $\mathcal{C}_b({\bf R}^n, {\bf R})$,
$\mathcal{C}_b({\bf R}^n, {\bf C})$ are dense subspaces of
$\mathcal{C}({\bf R}^n, {\bf R})$, $\mathcal{C}({\bf R}^n, {\bf C})$,
respectively, and are not closed subspaces in particular.
However, for each positive real number $r$, the collections
\begin{equation}
	\{f \in \mathcal{C}_b({\bf R}^n, {\bf R}) : \|f\|_{sup} \le r\}
\end{equation}
and
\begin{equation}
	\{f \in \mathcal{C}_b({\bf R}^n, {\bf C}) : \|f\|_{sup} \le r\}
\end{equation}
are closed subsets of $\mathcal{C}({\bf R}^n, {\bf R})$,
$\mathcal{C}({\bf R}^n, {\bf C})$, respectively.  This is not
difficult to verify.

	Suppose that $\lambda$ is a linear functional on
$\mathcal{C}_{00}({\bf R}^n, {\bf R})$ or $\mathcal{C}_{00}({\bf R}^n,
{\bf C})$ which is bounded with respect to the supremum norm.
This means that there is a nonnegative real number $L$ such that
\begin{equation}
	|\lambda(f)| \le L \, \|f\|_{sup}
\end{equation}
for all $f$ in $\mathcal{C}_{00}({\bf R}^n, {\bf R})$ or
$\mathcal{C}_{00}({\bf R}^n, {\bf C})$, as appropriate.  This is
equivalent to saying that $\lambda$ extends to a continuous linear
functional on $\mathcal{C}_0({\bf R}^n, {\bf R})$ or
$\mathcal{C}_0({\bf R}^n, {\bf C})$, as appropriate.

	Let us assume that $L$ is the smallest nonnegative real number
for which the inequality above holds.  This is the same as saying that
$L$ is the supremum of $|\lambda(f)|$ for all real or complex-valued
continuous functions on ${\bf R}^n$, as appropriate, which have
compact support and satisfy $\|f\|_{sup} \le 1$.  Thus, for each
$\epsilon > 0$, there is a continuous function $h_\epsilon$ on ${\bf
R}^n$ with compact support $K_\epsilon$ such that
$\|h_\epsilon\|_{sup} \le 1$ and $|\lambda(h_\epsilon)| \ge L -
\epsilon$.

	Using these properties of $h_\epsilon$ it follows that
if $f$ is a continuous function on ${\bf R}^n$ such that
$f(x) = 0$ for all $x \in K_\epsilon$ and $\|f\|_{sup} \le 1$,
then
\begin{equation}
	|\lambda(f)| \le \epsilon.
\end{equation}
Indeed, if $\alpha$, $\beta$ are scalars such that 
\begin{equation}
	|\alpha|, |\beta| \le 1,
\end{equation}
then $\alpha \, f + \beta \, h_\epsilon$ is a continuous function
on ${\bf R}^n$ with compact support such that
\begin{equation}
	\|\alpha \, f + \beta \, h_\epsilon\|_{sup} \le 1.
\end{equation}
Hence 
\begin{equation}
	|\alpha \, \lambda(f) + \beta \, \lambda(h)| \le L,
\end{equation}
and one can use this to bound $|\lambda(f)|$.

	From here it is not too hard to show that for each positive
real number $r$ the restriction of $\lambda$ to the collection of real
or complex-valued continuous functions $f$ on ${\bf R}^n$ with compact
support and $\|f\|_{sup} \le r$ is continuous with respect to the
topology of $\mathcal{C}({\bf R}^n, {\bf R})$ or $\mathcal{C}({\bf
R}^n, {\bf C})$, as appropriate.  This topology is less restrictive
than the one associated to the supremum norm.  In other words, this is
a stronger continuity property than the one that was explicitly
assumed at the beginning.

	As a result, one can show that $\lambda$ extends to a linear
functional on $\mathcal{C}_b({\bf R}^n, {\bf R})$ or
$\mathcal{C}_b({\bf R}^n, {\bf C})$ which is also continuous on each
bounded subset with respect to the topology of $\mathcal{C}({\bf R}^n,
{\bf R})$, $\mathcal{C}({\bf R}^n, {\bf C})$, as appropriate.  We still
have that
\begin{equation}
	|\lambda(f)| \le L \, \|f\|_{sup}
\end{equation}
for all bounded continuous functions $f$ now, i.e., with the same
choice of $L$ as before.  This extension of $\lambda$ to bounded
continuous functions on ${\bf R}^n$ is unique, basically because every
bounded continuous function $f$ on ${\bf R}^n$ is the limit of a
sequence $\{f_j\}_{j=1}^\infty$ of countinuous functions on ${\bf
R}^n$ with compact support such that the supremum norms
$\|f_j\|_{sup}$ are uniformly bounded and the convergence is uniform
on bounded subsets of ${\bf R}^n$.

	If $f_1$, $f_2$ are bounded continuous functions on ${\bf
R}^n$, then the product $f_1 \, f_2$ is also a bounded continuous
function on ${\bf R}^n$, and
\begin{equation}
	\|f_1 \, f_2\|_{sup} \le \|f_1\|_{sup} \, \|f_2\|_{sup}.
\end{equation}
Thus $\mathcal{C}_b({\bf R}^n, {\bf R})$, $\mathcal{C}_b({\bf R}^n,
{\bf C})$ are normed algebras.  The product operation also behaves
well in terms of continuity with respect to the less restrictive
topologies from $\mathcal{C}({\bf R}^n, {\bf R})$, $\mathcal{C}({\bf
R}^n, {\bf C})$.

\section{Continuous functions on ${\bf R}^n$ of polynomial growth}
\label{section on continuous functions on R^n of polynomial growth}

	Let us define $\mathcal{PC}({\bf R}^n, {\bf R})$,
$\mathcal{PC}({\bf R}^n, {\bf C})$ to be the real and complex vector
spaces of real and complex-valued continuous functions $f(x)$ on ${\bf
R}^n$, respectively, which have \emph{polynomial
growth},\index{polynomial growth of a real or complex-valued function
on ${\bf R}^n$} which is to say that there is a nonnegative real
number $C$ and a positive integer $j$ such that
\begin{equation}
	|f(x)| \le C \, (1 + |x|)^j
\end{equation}
for all $x \in {\bf R}^n$.  Instead of defining topologies on
$\mathcal{PC}({\bf R}^n, {\bf R})$, $\mathcal{PC}({\bf R}^n, {\bf
C})$, let us make the conventions that a subset $E$ of one of these
spaces is bounded if there are choices of $C, j$ so that each $f \in
E$ satisfies the inequality above, a sequence $\{f_j\}_{j=1}^\infty$
in one of these spaces converges to a function $f$ in the same space
if the sequence is contained in a bounded subset and the $f_j$'s
converge to $f$ uniformly on compact subsets of ${\bf R}^n$, and that
a sequence $\{f_j\}_{j=1}^\infty$ in one of these spaces is a Cauchy
sequence if it is contained in a bounded subset of the space and it
satisfies the usual Cauchy condition with respect to the supremum norm
on any compact subset of ${\bf R}^n$.  By standard reasoning each such
Cauchy sequence converges uniformly on compact subsets of ${\bf R}^n$
to a continuous function on ${\bf R}^n$ of polynomial growth,
so that each Cauchy sequence converges.

	In other words, on each bounded subset of $\mathcal{PC}({\bf
R}^n, {\bf R})$, $\mathcal{PC}({\bf R}^n, {\bf C})$, we are using the
same topology as for $\mathcal{C}({\bf R}^n, {\bf R})$,
$\mathcal{C}({\bf R}^n, {\bf C})$.  Notice that for each nonnegative
real number $C$ and each positive integer $j$, the set of real or
complex-valued continuous functions $f(x)$ on ${\bf R}^n$ such that
$|f(x)|$ is bounded by $C \, (1 + |x|)^j$ on ${\bf R}^n$ is a closed
subset of $\mathcal{C}({\bf R}^n, {\bf R})$, $\mathcal{C}({\bf R}^n,
{\bf C})$, respectively.  Also, convergence of a sequence in such a set
uniformly on compact subsets of ${\bf R}^n$ is equivalent to convergence
with respect to the norm 
\begin{equation}
	\sup \{|f(x)| \, (1 + |x|)^{j+1} : x \in {\bf R}^n\}.
\end{equation}

	Sums and products of continuous functions of polynomial growth
on ${\bf R}^n$ are continuous functions with polynomial growth, and
these operations behave well in terms of bounded collections of
functions and convergence of sequences of functions in the senses
described above.  One can also think of multiplication as defining
binary operations
\begin{eqnarray}
  \mathcal{PC}({\bf R}^n, {\bf R}) \times \mathcal{RC}({\bf R}^n, {\bf R})
	& \to & \mathcal{RC}({\bf R}^n, {\bf R}),		\\
  \mathcal{PC}({\bf R}^n, {\bf C}) \times \mathcal{RC}({\bf R}^n, {\bf C})
	& \to & \mathcal{RC}({\bf R}^n, {\bf C}),
\end{eqnarray}
and these operations behave well in terms of bounded sets and
convergent sequences in these spaces too.  Ordinary polynomials
on ${\bf R}^n$ are continuous functions with polynomial growth,
and standard considerations for them fit in a nice way with the
spaces of general continuous functions of polynomial growth, e.g.,
a collection of polynomials is bounded in the sense discussed
here if their degrees and coefficients are bounded.

\section{More on spaces of sequences}
\label{section with more on spaces of sequences}

	Let $1 \le p, q \le \infty$ be conjugate exponents, so that
\begin{equation}
	\frac{1}{p} + \frac{1}{q} = 1.
\end{equation}
If $w = \{w_j\}_{j=1}^\infty$ is an element of $\ell^q({\bf R})$,
$\ell^q({\bf C})$, then
\begin{equation}
	\lambda_w(x) = \sum_{j=1}^\infty x_j \, w_j, 
					\quad x = \{x_j\}_{j=1}^\infty,
\end{equation}
defines a bounded linear functional on $\ell^p({\bf R})$, $\ell^p({\bf
C})$, respectively.  This uses H\"older's inequality to say that the
sum converges, and also that
\begin{equation}
	|\lambda_w(x)| \le \|w\|_q \, \|x\|_p.
\end{equation}

	We can turn this around a bit and say that for each $w$
in $\ell^q({\bf R})$, $\ell^q({\bf C})$, the linear functional
$\lambda_w$ leads to a seminorm
\begin{equation}
	|\lambda_w(x)|
\end{equation}
on $\ell^p({\bf R})$, $\ell^p({\bf C})$, and that the family of these
seminorms defines a topology on $\ell^p({\bf R})$, $\ell^p({\bf C})$
which gives the structure of a locally convex topological vector
space.  This topology is less restrictive than the one associated to
the norm $\|x\|_p$.  In the case of $p = 1$, there is a variant of this
where one only uses the family of seminorms corresponding to $w$
in $c_0({\bf R})$, $c_0({\bf C})$, respectively, and this version is
nicer in several respects.

	Let $r$ be a positive real number, and consider the subset of
$\ell^p({\bf R})$ or $\ell^p({\bf C})$ consisting of $x =
\{x_j\}_{j=1}^\infty$ such that
\begin{equation}
	\|x\|_p \le r.
\end{equation}
Of course this is a convex set which is closed with respect to the
topology associated to the norm $\|x\|_p$.  One can also check that
this set is closed with respect to the topology coming from linear
functionals as in the previous paragraph, using either $c_0$ or
$\ell^\infty$ when $p = 1$.

	Let $E$ be a dense subset of $\ell^q({\bf R})$ or $\ell^q({\bf
C})$, as appropriate, or of $c_0({\bf R})$, $c_0({\bf C})$ when $p =
1$.  Actually, it is enough to choose $E$ so that its span is dense
for the present purposes.  When $1 \le q < \infty$, or if one uses $c_0$
instead of $\ell^\infty$ when $p = 1$ and $q = \infty$, one can take
$E$ to be a countable set, which is a nice thing to do.

	On the set of $x = \{x_j\}_{j=1}^\infty$ in $\ell^p({\bf R})$
or $\ell^p({\bf C})$, as appropriate, such that $\|x\|_p \le r$, 
the topology that one gets from the seminorms $|\lambda_w(y)|$, $w \in E$,
is the same as the topology that one gets from using the seminorms
corresponding to all $w$ in $\ell^q$ or $c_0$, as appropriate.  This is
not too difficult to verify, and it would not work on the whole
$\ell^p$ space instead of just this bounded part of it.  This uses the
fact that if $w_1$, $w_2$ are close in $\ell^q$ or $c_0$, then 
\begin{equation}
	\lambda_{w_1} - \lambda_{w_2} = \lambda_{w_1 - w_2}
\end{equation}
is uniformly small on this bounded subset of $\ell^p$.

	Let us consider some other spaces of sequences, namely,
$r\ell({\bf R})$, $r\ell({\bf C})$, consisting of sequences which
decay rapidly, and $p\ell({\bf R})$, $r\ell({\bf C})$, consisting of
sequences with polynomial growth.  To be more precise, $r\ell({\bf
R})$, $r\ell({\bf C})$ are the real and complex vector spaces of real
and complex sequences $x = \{x_j\}_{j=1}^\infty$, respectively, such
that for each positive integer $k$ there is a nonnegative real number
$C(k)$ so that
\begin{equation}
	|x_j| \le C(k) \, j^{-k}
\end{equation}
for all $j \in {\bf Z}_+$.  Also, $p\ell({\bf R})$, $p\ell({\bf C})$
are the real and complex vector spaces of real and complex sequences
$x = \{x_j\}_{j=1}^\infty$ for which there is a positive integer $l$
and a nonnegative real number $C$ such that
\begin{equation}
	|x_j| \le C \, j^k
\end{equation}
for all $j \in {\bf Z}_+$.

	On $r\ell({\bf R})$, $r\ell({\bf C})$, we have for each positive
integer $k$ the seminorm defined by
\begin{equation}
	\sup \{j^k \, |x_j| : j \in {\bf Z}_k\},
\end{equation}
and this family of seminorms defines topologies on $r\ell({\bf R})$,
$r\ell({\bf C})$, so that they become locally convex topological
vector spaces.  A subset $E$ of $r\ell({\bf R})$, $r\ell({\bf C})$ is
bounded if and only if for each positive integer $k$ there is a
nonnegative real number $C(k)$ such that $|x_j| \le C(k) \, j^{-k}$
for all $x \in E$ and $j \in {\bf Z}_+$.  As usual, one can check that
$r\ell({\bf R})$, $r\ell({\bf C})$ are complete, so that they are
actually Fr\'echet spaces.

	In $p\ell({\bf R})$, $p\ell({\bf C})$, we can make the convention
that a subset $E$ is bounded if and only if there is a positive integer
$l$ and a nonnegative real number $C$ such that
\begin{equation}
	|x_j| \le C \, j^l
\end{equation}
for all $x = \{x_j\}_{j=1}^\infty$ in $E$ and all $j \in {\bf Z}_+$.
We can make the convention that a sequence of elements of $p\ell({\bf
R})$ or $p\ell({\bf C})$ converges to another element in the same
space if the sequence is contained in a bounded subset as in the
previous sentence and one has pointwise convergence.  We can also make
the convention that a sequence of elements of $p\ell({\bf R})$ or
$p\ell({\bf C})$ is a Cauchy sequence when it is contained in a
bounded set and one has Cauchy sequences pointwise, so that Cauchy
sequences always converge in these spaces.

\section{Weak topologies and bilinear forms}
\label{section on weak topologies and bilinear forms}

	Let $V$ be a topological vector space, and let $\lambda$ be a
linear functional on $V$.\index{continuous linear functionals on a
topological vector space} If $\lambda$ is continuous, then there is an
open subset $U$ of $V$ such that $0 \in U$ and
\begin{equation}
	|\lambda(v)| < 1
\end{equation}
for all $v \in U$.  The converse is also true, which is to say that if
there is such an open subset $U$ of $V$, then $\lambda$ is continuous,
as one can check using the linearity of $\lambda$.

	Now suppose that $V$ is a vector space and that $\mathcal{F}$
is a nonempty family of seminorms\index{seminorms on vector spaces} on
$V$.  In the topology on $V$ generated by $\mathcal{F}$, the sets
\begin{equation}
	\{v \in V : N(v) < r\},
\end{equation}
where $N$ is a seminorm on $V$ in the family $\mathcal{F}$ and $r$
is a positive real number, are open subsets of $V$, and in fact form
a sub-basis for the topology of $V$ at $0$.  In other words, if $U$
is an open subset of $V$ in this topology such that $0 \in U$, then
there are seminorms $N_1, \ldots, N_l$ in the family $\mathcal{F}$
and positive real numbers $r_1, \ldots, r_l$ such that
\begin{equation}
	\cap_{i=1}^l \{v \in V : N_i(v) < r_i\}  \subseteq U.
\end{equation}

	Let $\lambda$ be a linear functional on $V$ again.  It follows
that $\lambda$ is continuous if and only if there are seminorms $N_1,
\ldots, N_l$ in the family $\mathcal{F}$ and positive real numbers
$r_1, \ldots, r_l$ such that $|\lambda(v)| < 1$ for all $v \in V$ such
that $N_i(v) < r_i$, $1 \le i \le l$.  This is equivalent to saying
that there are seminorms $N_1, \ldots, N_l$ in the family $\mathcal{F}$
and nonnegative real numbers $C_1, \ldots, C_l$ such that
\begin{equation}
	|\lambda(v)| \le \max \{C_i \, N_i(v) : 1 \le i \le l\}
\end{equation}
for all $v \in V$.\index{continuous linear functionals on a topological
vector space}

	Assume now that $V$ and $W$ are vector spaces, both real or
both complex, and that $B(v, w)$ is a bilinear form on $V \times W$.
This means that $B(v, w)$ is a function from $V \times W$ to the real
or complex numbers, whichever are being used as scalars, such that $v
\mapsto B(v, w)$ is a linear functional on $V$ for each $w \in W$ and
$w \mapsto B(v,w)$ is a linear functional on $W$ for each $v \in V$.
Let us also assume that $B(v, w)$ satisfies the nondegeneracy property
that for each $v \in V$ with $v \ne 0$ there is a $w \in W$ such that
\begin{equation}
	B(v, w) \ne 0.
\end{equation}

	We can associate to $B$ the family $\mathcal{F}(B)$ of seminorms
on $V$ given by $|B(v, w)|$, $w \in W$, and this defines a topology on $V$
which makes $V$ into a topological vector space.  By construction, each
linear functional on $V$ of the form $B(v, w)$ for some $w \in W$ is
a continuous linear functional with respect to the topology just defined.
Conversely, one can check that every continuous linear functional on $V$
with respect to this topology is of this form.

\chapter{Examples, 3}
\label{chapter on examples, 3}

\section{Power series and spaces of sequences}
\label{section on power series and spaces of sequences}

	In this section we use complex numbers as scalars.  A
\emph{power series}\index{power series} is a series of the form
\begin{equation}
	\sum_{n=0}^\infty a_n \, z^n,
\end{equation}
where the coefficients $a_n$ are complex numbers and $z$ is a complex
variable.  One can add power series, multiply them by complex numbers,
or even multiply two power series to get another one.

	A power series as above converges at some $z_0$ in ${\bf C}$
if the infinite series of scalars
\begin{equation}
	\sum_{n=0}^\infty a_n \, z_0^n
\end{equation}
converges in the usual sense, and the power series converges
absolutely at $z_0$ if this series converges absolutely in the usual
sense, which is to say that
\begin{equation}
	\sum_{n = 0}^\infty |a_n| \, |z_0|^n
\end{equation}
converges.  A well-known fact states that a power series which
converges for some $z_0 \in {\bf C}$ converges absolutely for all $z
\in {\bf C}$ such that $|z| < |z_0|$.  It follows that there is a
radius $R$, $0 \le R \le \infty$, called the \emph{radius of
convergence}\index{radius of convergence of a power series} of the
power series, such that the power series converges absolutely when
$|z| < R$ and does not converge when $|z| > R$, with the behavior for
$|z| = R$ depending on the situation.

	Another basic fact is that on $\{z \in {\bf C} : |z| < R\}$,
the partial sums of the power series converges uniformly on compact
subsets, and hence the power series defines a continuous function
$f(z)$ on this region.  Moreover, $f(z)$ is holomorphic, which means
that
\begin{equation}
	f'(z) = \lim_{h \to 0} \frac{f(z+h) - f(z)}{h}
\end{equation}
exists for all $z \in {\bf C}$ with $|z| < R$.  In fact, the formal
power series for $f'(z)$, given by
\begin{equation}
	\sum_{n = 1}^\infty n \, z^{n-1},
\end{equation}
has the same radius of convergence as the one for $f(z)$, and is equal
to $f'(z)$ for $z \in {\bf C}$, $|z| < R$.  

	If $0 < t \le \infty$, let us define $\mathcal{A}_t$ to be the
complex vector space of power series $\sum_{n=0}^\infty a_n \, z^n$
which converge for $|z| < t$, which is the same as saying that the
radius of convergence is at least $t$.  For each positive real number
$r < t$ and power series $\sum_{n=0}^\infty a_n$ in $\mathcal{A}_r$,
consider the expression
\begin{equation}
\label{sup { |a_n| r^n : n ge 0}}
	\sup \{ |a_n| \, r^n : n \ge 0\}.
\end{equation}
This is finite since $\sum_{n=0}^\infty a_n \, z^n$ is assumed to
converge when $|z| = r$, and conversely $\sum_{n=0}^\infty a_n \, z^n$
converges for $z \in {\bf C}$ with $|z| < t$ when these quantities
are finite for each $r \in (0,t)$.

	 The expression (\ref{sup { |a_n| r^n : n ge 0}}) defines a
seminorm on $\mathcal{A}_t$ for each $r \in (0,t)$, and this family of
seminorms determines a topology on $\mathcal{A}_t$ which makes it into
a locally convex topological vector space.  It suffices to consider
$r$'s in an increasing sequence that tends to $t$, so that in fact one
only needs a countable family of seminorms here.  It is not difficult
to verify that $\mathcal{A}_t$ is complete, so that it is actually a
Fr\'echet space.

	Every power series in $\mathcal{A}_t$ defines a holomorphic
function on $\{z \in {\bf C} : |z| < t\}$, and basic results in complex
analysis imply that every holomorphic function on this region can
be represented by a power series.  If $f(z)$ is a holomorphic function
defined for $z \in {\bf C}$, $|z| < t$, and if $s$ is a positive real
number with $s < t$, then one can consider the quantity
\begin{equation}
	\sup \{|f(z)| : z \in {\bf C}, |z| \le s\}.
\end{equation}
This is another family of seminorms on holomorphic functions on the
region $\{z \in {\bf C} : |z| < t\}$, and one can show that it is
equivalent to the earlier family of seminorms, in the sense that they
define the same topology on $\mathcal{A}_t$.

	The product of two holomorphic functions on $\{z \in {\bf C} :
|z| < t\}$ is also a holomorphic function on this region.  At the
level of power series, we have
\begin{equation}
	\biggl(\sum_{j=0}^\infty a_j \, z^j\biggr)
		\, \biggl(\sum_{k=0}^\infty b_k \, z^k\biggr)
		= \sum_{n=0}^\infty c_m \, z^m
\end{equation}
with
\begin{equation}
	c_m = \sum_{l=0}^m a_l \, b_{m-l},
\end{equation}
and $\sum_{m=0}^\infty c_m \, z^m$ converges when $|z| < t$ if
$\sum_{j=0}^\infty a_j \, z^j$, $\sum_{k=0}^\infty b_k \, z^k$
converge when $|z| < t$.  Multiplication is continuous on
$\mathcal{A}_t$, so that $\mathcal{A}_t$ defines a topological
algebra.

\section{Absolutely convergent Fourier series}
\label{section on absolutely convergent Fourier series}

	In this section we again restrict ourselves to complex numbers
as scalars.  Let us write $\mathcal{AC}$ for the complex vector space
of \emph{absolutely convergent Fourier series},\index{absolutely
convergent Fourier series} by which we mean series of the form
\begin{equation}
\label{sum_{n = -infty}^infty a_n z^n}
	\sum_{n = -\infty}^\infty a_n \, z^n,
\end{equation}
where we think of $z$ as a variable taking values in the unit
circle\index{unit circle}
\begin{equation}
	{\bf T} = \{z \in {\bf C} : |z| = 1\},
\end{equation}
and where the series
\begin{equation}
	\sum_{n = -\infty}^\infty |a_n|
\end{equation}
converges.  Thus (\ref{sum_{n = -infty}^infty a_n z^n}) converges
absolutely for each $z \in {\bf T}$, and the partial sums converge
uniformly to a continuous function on ${\bf T}$.

	On the unit circle the series (\ref{sum_{n = -infty}^infty a_n
z^n}) is equal to
\begin{equation}
\label{sum_{n=1}^infty a_{-n} overline{z}^n + sum_{n=0}^infty a_n z^n}
	\sum_{n=1}^\infty a_{-n} \, \overline{z}^n
		+ \sum_{n=0}^\infty a_n \, z^n,
\end{equation}
where $\overline{z}$ denotes the complex conjugate of $z$.  Under the
assumption that $\sum_{n=-\infty}^\infty |a_n|$ converges, the series
in (\ref{sum_{n=1}^infty a_{-n} overline{z}^n + sum_{n=0}^infty a_n
z^n}) converge uniformly and absolutely on the closed unit disk, which
is the set of $z \in {\bf C}$ such that $|z| \le 1$.  The sum defines
a continuous function $f(z)$ on the closed unit disk which is smooth
and harmonic in the interior, i.e., $\Delta f(z) = 0$ when $|z| < 1$,
where $\Delta$ is the Laplace operator $\partial^2/\partial x^2 +
\partial^2/\partial y^2$, $z = x + y i$.

	On $\mathcal{AC}$ we have a natural norm, namely,
$\sum_{n=-\infty}^\infty |a_n|$.  This is basically the same as
$\ell^1({\bf C})$ as a normed vector space.  In particular,
$\mathcal{AC}$ is complete with respect to this norm, and thus
becomes a Banach space.

	We can also multiply these series, i.e.,
\begin{equation}
	\biggl(\sum_{j=-\infty}^\infty a_j \, z^j \biggr)
		\, \biggl(\sum_{k=-\infty}^\infty b_k \, z^k\biggr)
		= \sum_{m = -\infty}^\infty c_m \, z^m,
\end{equation}
where
\begin{equation}
	c_m = \sum_{l=-\infty}^\infty a_l \, b_{m-l}.
\end{equation}
Assuming that 
\begin{equation}
	\sum_{j=-\infty}^\infty |a_j|, \quad \sum_{k=-\infty}^\infty |b_k|
\end{equation}
converge, one can check that the series defining $c_m$ converges, and that
\begin{equation}
	\sum_{m=-\infty}^\infty |c_m| 
		\le \biggl(\sum_{j=-\infty}^\infty |a_j| \biggr)
			\, \biggl(\sum_{k=-\infty}^\infty |b_k| \biggr).
\end{equation}
In other words, $\mathcal{AC}$ is a normed algebra with respect to
multiplication of series.

	Multiplication of series in this way corresponds exactly to
multiplying the functions on the unit circle defined by the series.
This does not work for the extensions to the closed unit disk
described above, although that is compatible with addition and scalar
multiplication.  Alternatively, under suitable conditions on the
coefficients $a_n$, the series (\ref{sum_{n = -infty}^infty a_n z^n})
can converge on an annular region around the unit circle, defining a
holomorphic function on that region, and products of these ``Laurent''
series correspond to products of the holomorphic functions that they
define.

\section{Spaces of smooth functions}
\label{section on spaces of smooth functions}

	For each positive integer $k$, let $\mathcal{C}^k({\bf R}^n,
{\bf R})$, $\mathcal{C}^k({\bf R}^n, {\bf C})$ denote the real and
complex vector spaces of real and complex-valued functions $f(x)$ on
${\bf R}^n$, respectively, which are continuously differentiable of
order $k$,\index{continuous differentiability of order $k$} i.e.,
which are continuous and for which the partial derivatives of $f(x)$
up to and including order $k$ exist at each point in ${\bf R}^n$ and
are continuous functions.  We also allow $k = \infty$ here, for the
spaces of continuous functions for which partial derivatives of all
orders exist and are continuous.  For $1 \le k \le \infty$ we write
$\mathcal{C}^k_{00}({\bf R}^n, {\bf R})$, $\mathcal{C}^k_{00}({\bf
R}^n, {\bf C})$ for the vector spaces of continuous functions on ${\bf
R}^n$ with compact support which are continuously differentiable of
order $k$.

	By a \emph{multi-index}\index{multi-indices} we mean an
$n$-tuple $\alpha = (\alpha_1, \ldots, \alpha_n)$ of nonnegative
integers, and we write $|\alpha|$ for $|\alpha_1| + \cdots +
|\alpha_n|$.  For each such multi-index $\alpha$, we have the
corresponding partial derivative of order $|\alpha|$ defined by
\begin{equation}
	\partial^\alpha = 
  \frac{\partial^|\alpha|}{\partial x_1^{\alpha_1} 
				\cdots \partial x_n^{\alpha_n}}.
\end{equation}
When $|\alpha| = 0$ this is interpreted as being the identity operator,
so that $\partial^\alpha f(x) = f(x)$.

	Suppose that $1 \le k \le \infty$, $l$ is a positive integer,
and $\alpha$ is a multi-index such that $|\alpha| \le k$.  If $f(x)$
is a real or complex-valued function on ${\bf R}^n$ which is
continuously differentiable of order $k$, consider the quantity
\begin{equation}
	\sup \{|\partial^\alpha f(x)| : x \in {\bf R}^n, |x| \le l\}.
\end{equation}
This defines a seminorm on $\mathcal{C}^k({\bf R}^n, {\bf R})$,
$\mathcal{C}^k({\bf R}^n, {\bf C})$.

	With this family of seminorms, $\mathcal{C}^k({\bf R}^n, {\bf
R})$, $\mathcal{C}^k({\bf R}^n, {\bf C})$ become topological vector
spaces.  As usual, it is not too difficult to show that these spaces
are complete, so that they are in fact Fr\'echet spaces.  Also, the
product of two functions which are continuously differentiable of
order $k$ is also continuously differentiable of order $k$, and these
spaces are topological algebras with respect to ordinary
multiplication.

	The spaces $\mathcal{C}^k_{00}({\bf R}^n, {\bf R})$,
$\mathcal{C}^k_{00}({\bf R}^n, {\bf C})$ have natural topologies in
which they are inductive limit spaces.  Two main points are that
bounded subsets of these spaces consist of functions with support
contained in a fixed compact subset of ${\bf R}^n$, and that when one
restricts to a subspace of functions with support contained in a fixed
compact subset of ${\bf R}^n$, the topology is determined by seminorms
as above, with $l$ chosen large enough so that the supports are
contained in the closed ball in ${\bf R}^n$ with center $0$ and radius
$l$.  Continuous linear functionals on $\mathcal{C}^\infty_{00}({\bf
R}^n, {\bf R})$, $\mathcal{C}^\infty_{00}({\bf R}^n, {\bf R})$ are
often called \emph{distributions}\index{distributions on ${\bf R}^n$}
or \emph{generalized functions}.\index{generalized functions on ${\bf R}^n$}

	Now let us define the \emph{Schwartz classes}\index{Schwartz
classes of rapidly decreasing smooth functions on ${\bf R}^n$}
$\mathcal{S}({\bf R}^n, {\bf R})$, $\mathcal{S}({\bf R}^n, {\bf C})$
of rapidly decreasing smooth functions on ${\bf R}^n$.  Namely,
a real or complex-valued function $f(x)$ on ${\bf R}^n$ is in the 
Schwartz class if $f(x)$ is infinitely-differentiable and if for
each multi-index $\alpha$ the function $\partial^\alpha f(x)$
is rapidly decreasing on ${\bf R}^n$.  This is equivalent to saying
that for each pair of multi-indices $\alpha$, $\beta$, the function
$x^\beta \, \partial^\alpha f(x)$ is bounded on ${\bf R}^n$, where
\begin{equation}
	x^\beta = x_1^{\beta_1} \cdots x_n^{\beta_n},
\end{equation}
and this is interpreted as being equal to $1$ when $|\beta| = 0$.

	Thus we have a natural family of seminorms on the Schwartz
classes, given by
\begin{equation}
	\sup \{|x^\beta \, \partial^\alpha f(x)| : x \in {\bf R}^n\},
\end{equation}
where $\alpha$, $\beta$ run through all multi-indices.  This family
of seminorms define topologies on the Schwartz classes, so that they
become locally convex topological vector spaces.  One can check that
these spaces are complete, and hence are Fr\'echet spaces.

	It is easy to see that the product of two functions in
$\mathcal{S}({\bf R}^n, {\bf R})$, $\mathcal{S}({\bf R}^n, {\bf C})$
lies in the same space, and in fact that the Schwartz classes are
topological algebras.  This is not really the whole picture, since for
instance the product of a function in the Schwartz class and a
polynomial lies in the Schwartz class.  Also, derivatives of functions
in the Schwartz class lie in the Schwartz class.

	Continuous linear functionals on the Schwartz classes are
called \emph{tempered distributions}.\index{tempered distributions on
${\bf R}^n$} In the case of complex-valued functions, a key point
about the Schwartz class is that the Fourier transform takes
$\mathcal{S}({\bf R}^n, {\bf C})$ to itself.  This leads to a nice
theory of Fourier transforms of tempered distributions.

	If $\phi(x)$, $f(x)$ are functions on ${\bf R}^n$, then the
\emph{convolution}\index{convolution of functions on ${\bf R}^n$} is
denoted $\phi * f(x)$ and is basically defined by
\begin{equation}
	\phi * f(x) = \int_{{\bf R}^n} \phi(y) \, f(x - y) \, dy.
\end{equation}
Depending on the circumstances, one should perhaps be careful about
the integral.  For instance, this makes sense if both $\phi$, $f$
are continuous and one of them has compact support, or if one of
them lies in the class of rapidly decreasing continuous functions on
${\bf R}^n$ and the other lies in the class of continuous functions
with polynomial growth.

	More precisely, if $\phi$ is a continuous function on ${\bf R}^n$
with compact support, then
\begin{equation}
	f \mapsto \phi * f
\end{equation}
defines a continuous linear mapping from the vector space of
continuous functions on ${\bf R}^n$ to itself.  One can take $\phi$, $f$
to both be real-valued or complex-valued here.  If both $\phi$ and $f$
have compact support, then so does $\phi * f$.

	Similarly, if $\phi$, $f$ are rapidly decreasing continuous
functions on ${\bf R}^n$, then $\phi * f$ is also a rapidly decreasing
continuous function on ${\bf R}^n$.  If $\phi$ is a rapidly-decreasing
continuous function on ${\bf R}^n$ and $f$ is a bounded continuous
function on ${\bf R}^n$, then $\phi * f$ is a bounded continuous
function on ${\bf R}^n$.  If $\phi$ is a rapidly-decreasing continuous
function on ${\bf R}^n$ and $f$ is a continuous function on ${\bf
R}^n$ with polynomial growth, then $\phi * f$ is a continuous function
on ${\bf R}^n$ with polynomial growth.

	The operation of convolution is a nice commutative and
associative operation, a kind of product.  For instance, if $\phi_1$,
$\phi_2$, $\phi_3$ are rapidly-decreasing continuous functions on
${\bf R}^n$, then
\begin{equation}
	\phi_1 * \phi_2 = \phi_2 * \phi_1
\end{equation}
and
\begin{equation}
	\phi_1 * (\phi_2 * \phi_3) = (\phi_1 * \phi_2) * \phi_3.
\end{equation}
The spaces of rapidly-decreasing real or complex-valued continuous
functions on ${\bf R}^n$ are topological algebras with respect to
convolution, as well as with respect to ordinary multiplication.

	Convolutions also behave well in terms of differentiation,
in the sense that
\begin{equation}
	\partial^{\alpha + \beta} (\phi * f) 
		= (\partial^\alpha \phi) * (\partial^\beta f)
\end{equation}
under suitable conditions on $\phi$, $f$.  For instance, this holds
when $\phi$ is continuously-differentiable of order $l$, $f$ is
continuously-differentiable of order $m$, $|\alpha| \le l$, $|\beta|
\le m$, and at least one of $\phi$, $f$ has compact support.  In
particular, $\phi * f$ is continuously differentiable of order $l + m$
in this case.

\section{Banach algebras}
\label{section on banach algebras}

	Let $\mathcal{A}$ be a real or complex normed algebra, with
norm $\|a\|$, which is complete as a normed vector space.  In this
case we say that $\mathcal{A}$ is a \emph{Banach
algebra}.\index{Banach algebra} We shall make the standing assumption
that $\mathcal{A}$ contains a nonzero multiplicative identity element,
denoted $1$, with norm equal to the real number $1$.

	A basic class of Banach algebras consists of the algebras of
bounded linear operators on a Banach space.  We shall discuss this
case further in the next section.  A basic class of commutative Banach
algebras consists of the algebras of real or complex-valued continuous
functions on compact Hausdorff topological spaces, using the supremum
norm.

	An element $a$ of $\mathcal{A}$ is said to be \emph{invertible}
if there is a $b \in \mathcal{A}$ such that
\begin{equation}
	b \, a = a \, b = 1.
\end{equation}
Such an element $b$ is unique if it exists, and is denoted $a^{-1}$.
If $a_1$, $a_2$ are invertible elements of $\mathcal{A}$, then the
product $a_1 \, a_2$ is invertible, and 
\begin{equation}
	(a_1 \, a_2)^{-1} = a_2^{-1} \, a_1^{-1}.
\end{equation}

	A basic result about Banach algebras is that if $a \in
\mathcal{A}$ and $\|a\| < 1$, then $1 - a$ is an invertible element of
$\mathcal{A}$.  Indeed, we can use the usual formula
\begin{equation}
	(1 - a)^{-1} = \sum_{j = 0}^\infty a^j,
\end{equation}
where $a^0$ is taken to be $1$.  Because $\mathcal{A}$ is a Banach
algebra,
\begin{equation}
	\|a^j\| \le \|a\|^j,
\end{equation}
and this implies that the series above converges when $\|a\| < 1$, 
and that
\begin{equation}
	\|(1 - a)^{-1}\| \le \frac{1}{1 - \|a\|}.
\end{equation}

	More generally, if $x$ is an invertible element of $\mathcal{A}$,
and if $a$ is an element of $\mathcal{A}$ such that
\begin{equation}
	\|a\| < \|x^{-1}\|^{-1},
\end{equation}
then $x - a$ is an invertible element of $\mathcal{A}$, and
\begin{equation}
	\|(x - a)^{-1}\| \le \frac{\|x^{-1}\|}{1 - \|x^{-1}\| \, \|a\|}.
\end{equation}
Indeed, $x - a = x (1 - x^{-1} \, a)$, which is then the product of
two invertible elements of $\mathcal{A}$.  Thus the group of
invertible elements of $\mathcal{A}$, under multiplication, is an open
subset of $\mathcal{A}$.

	Let us restrict our attention for the rest of the section to
the situation where complex numbers are used as scalars.  If $\lambda$
is a complex number, then we can also view $\lambda$ as an element of
$\mathcal{A}$, namely, as a multiple of the multiplicative identity
element.  As is well known, even if one starts with real scalars,
often complex numbers are in the vicinity and play an important
role anyway.

	For each element $a$ of $\mathcal{A}$, the resolvent set
associated to $a$ is the set of complex numbers $\lambda$ such that
$\lambda - a$ is invertible in $\mathcal{A}$.  The spectrum of $a$ is
the complement of this set in ${\bf C}$, which is to say the set of
complex numbers $\lambda$ such that $\lambda - a$ is not invertible.
By the earlier remarks, the resolvent set is always an open subset of
${\bf C}$, and the spectrum is always a compact subset of ${\bf C}$.

	A famous result states that the spectrum of $a$, $a \in
\mathcal{A}$, always contains at least one element.  The basic idea of
the proof is that if the spectrum of $a$ were empty, then $(\lambda -
a)^{-1}$ would define a holomorphic function on all of $\mathcal{C}$,
with values in $\mathcal{A}$.  Because we can analyze $(\lambda -
a)^{-1}$ easily when $|\lambda|$ is large, we see that $(\lambda -
a)^{-1}$ should in fact be bounded, and behave like $1/\lambda$ when
$\lambda$ is large, while a bounded holomorphic function should be
constant.

	In any case, $(\lambda - a)^{-1}$ defines a holomorphic
function on the resolvent of $a$ with values in $\mathcal{A}$.  These
expressions can be used to define $f(a)$ whenever $f(z)$ is a
holomorphic complex-valued function defined on a neighborhood of the
spectrum of $a$.  Namely, one can use the Cauchy integral formula
applied to a family of curves in the resolvent of $a$ which also
lie in the domain of $f$ and surround the spectrum of $a$.

	If $n$ is a positive integer and $\lambda$ is a complex number
such that
\begin{equation}
	|\lambda|^n > \|a^n\|,
\end{equation}
then $\lambda$ lies in the resolvent set of $a$.  When $n = 1$ this
follows from the remarks near the beginning of the section.  In
general, one first applies those remarks to obtain that $\lambda^n -
a^n$ is invertible, and then it follows that $\lambda - a$ is
invertible, since
\begin{equation}
	\lambda^n - a^n = (\lambda - a) (\lambda^{n-1} + \lambda^{n-2} \, a
				+ \cdots + a^{n-1}).
\end{equation}

	The spectral radius of $a$, $a \in \mathcal{A}$, is defined to
be the maximum of $|\lambda|$, where $\lambda$ runs through the
spectrum of $a$.  Thus the spectral radius of $a$ is less than or
equal to the norm of $a$, and in fact the spectral radius of $a$ is
less than or equal to
\begin{equation}
	\|a^n\|^{1/n}
\end{equation}
for every positive integer $n$, by the remarks of the previous
paragraph.  A famous result states that the spectral radius of
$a$ is equal to
\begin{equation}
	\lim_{n \to \infty} \|a^n\|^{1/n}.
\end{equation}

	To prove this one can again use complex analysis.  Let $\rho$
denote the spectral radius of $a$, and let $r$ be any real number such
that $r > \rho$.  Consider the holomorphic function $(\lambda -
a)^{-1}$ again, on the region $\{z \in {\bf C} : |z| > \rho\}$,
and which we know also behaves well at the point at infinity.

	For reasons of compactness and continuity, $(\lambda -
a)^{-1}$ is bounded on the set of $\lambda \in {\bf C}$ such that
$|\lambda| = r$.  Using complex analysis, we can express $a^n$ in
terms of integrals of $(\lambda - a)^{-1}$ on the circle where
$|\lambda| = r$.  This leads to bounds for $a^n$ in terms of $r^n$, as
desired.

	A fundamental case occurs when $M$ is a compact Hausdorff
topological space, and $\mathcal{A}$ is the Banach algebra of
continuous complex-valued functions on $M$, using the supremum norm
\begin{equation}
	\|f\|_{sup} = \sup \{|f(x)| : x \in M\}.
\end{equation}
It is easy to see that a complex number $\lambda$ lies in the spectrum
of a complex-valued continuous function $f$ on $M$ if and only
if $\lambda$ lies in the image of $f$.  As a result, the spectral
radius of $f$ is equal to the norm of $f$ in this situation.

	Another very interesting case occurs when we consider the
Banach algebra of complex-valued continuous functions $f(z)$ on the
closed unit disk $\{z \in {\bf C} : |z| \le 1\}$ which are holomorphic
inside the unit disk.  In other words, we assume that there is a power
series
\begin{equation}
	\sum_{n=0}^\infty a_n \, z^n
\end{equation}
which converges for every $z \in {\bf C}$ with $|z| < 1$, and that the
function $f(z)$ on the open unit disk that it defines extends to a
continuous function on the closed unit disk.  It follows from standard
results that sums and products of such functions are again of the same
type.

	For the norm we use the supremum norm, and in fact by the
maximum principle for holomorphic functions the supremum of $|f(z)|$
for such a function $f$ over the $z \in {\bf C}$ with $|z| \le 1$ is
the same as the supremum over the $z \in {\bf C}$ such that $|z| = 1$.
Also, the coefficients $a_n$ of the power series expansion of $f$ on
the open unit disk can be given through well-known formulae as
integrals of $f$ on the unit circle.  In short, the functions on the
closed unit disk in this algebra are determined by their restrictions
to the unit circle, and we can think of this as a closed subalgebra of
the algebra of complex-valued continuous functions on the unit circle.

	In this case again the spectrum of a function $f(z)$ in the
algebra is given by the image of $f(z)$ as a function on the closed
unit disk.  This is rather striking when we think of the algebra
as a closed subalgebra of the complex-valued continuous functions on
the unit circle.  The spectral radius of a function in the algebra
is equal to the norm of the function.

	As a third example, let us consider the algebra $\mathcal{AC}$
of series of the form
\begin{equation}
	\sum_{n = -\infty}^\infty a_n \, z^n,
\end{equation}
where we think of the variable $z$ as taking values in the unit circle
\begin{equation}
	\{z \in {\bf C} : |z| = 1\},
\end{equation}
and where the coefficients $a_n$ are assumed to be absolutely summable,
so that
\begin{equation}
	\sum_{n = -\infty}^\infty |a_n|
\end{equation}
converges.  Basically our series is then the Fourier series of a
continuous function on the unit circle, and the coefficients $a_n$ can
be given in terms of integrals of the function on the circle.  We
define the norm in this case to be the sum of the moduli of the
coefficients $a_n$, and this is automatically greater than or equal to
the supremum norm of the corresponding function on the unit circle.

	In general this norm is strictly larger than the supremum norm
of the corresponding function $f(z)$ on the unit circle.  Also,
although one can define the Fourier coefficients $a_n$ of any
continuous function $f(z)$ on the unit circle, in general they may not
converge absolutely, as for elements of $\mathcal{AC}$.  However,
under modest additional regularity assumptions on $f(z)$, one can show
that the Fourier coefficients do converge absolutely, and that the
Fourier series converges to $f(z)$ everywhere on the unit circle.

	In this example it is a famous result that the spectrum of an
element of $\mathcal{AC}$ is equal to the image of the corresponding
continuous function on the unit circle, and hence that the spectral
radius is equal to the supremum norm of the function.  In other words,
when a continuous function on the unit circle which has an absolutely
convergent Fourier series is invertible as a simply a continuous function,
then the inverse also has absolutely convergent Fourier series.  In this
connection, let us note more broadly that for many regularity conditions
for a function, if one has a continuous function which satisfies that
regularity condition and which is invertible as a continuous function,
then the inverse also satisfies the same regularity condition.

\section{Bounded linear operators on Banach spaces}
\label{section about bounded linear operators on banach spaces}

	Let $V$ be a real or complex Banach space, with norm $\|v\|$.
If $T$ is a bounded linear operator on $V$, then the operator norm
of $T$ is denoted $\|T\|_{op}$ and defined by
\begin{equation}
	\|T\|_{op} = \sup \{\|T(v)\| : v \in V, \|v\| \le 1\}.
\end{equation}
Of course the identity operator $I$ on $V$ is a bounded linear
operator with norm equal to $1$.

	Using composition of linear operators as multiplication, the
bounded linear operators on $V$ becomes a Banach algebra.  It is easy
to see that the operator norm automatically satisfies the required
conditions, and one has completeness for the space of bounded linear
operators on $V$ because of completeness for $V$ itself.  A bounded
linear operator on $V$ is invertible as an element of this algebra
if and only if it is invertible in the usual sense, as a linear
operator on $V$, where the inverse is also bounded.

	If $T$ is a bounded linear operator on $V$, then $T$ is
invertible if and only if (i) the kernel of $T$, which is the linear
subspace of $V$ consisting of vectors $v$ such that $T(v) = 0$, is the
trivial subspace $\{0\}$, (ii) there is a positive real number $c$
such that
\begin{equation}
	\|T(v)\| \ge c \, \|v\|
\end{equation}
for all $v \in V$, and (iii) the image of $T$ is dense in $V$.
Condition (i) is of course implied by (ii), but it is simpler and
convenient to state separately.  In some cases (iii) can be derived
from (i), and we shall say more about this soon.

	More precisely, under condition (ii) above, the image of $T$
is a closed linear subspace of $V$.  Indeed, let
$\{w_j\}_{j=1}^\infty$ be any sequence in the image of $T$ which
converges to some $w \in V$.  We can write $w_j$ as $T(v_j)$, where
$\{v_j\}_{j=1}^\infty$ is a Cauchy sequence in $V$, because of
condition (ii), and therefore $\{v_j\}_{j=1}^\infty$ converges to some
$v \in V$, and $w = T(v)$, as desired.

	Because the image of $T$ is a closed subspace of $V$ and is
also dense, by condition (iii), we have that the image of $T$ maps $V$
onto $V$.  Of course condition (i) states that $T$ is one-to-one, so
that $T$ is invertible as a linear mapping of $V$ onto $V$.  Condition
(ii) implies that the inverse of $T$ is bounded, with norm less than
or equal to $c^{-1}$.

	There is a general result in the theory of Banach spaces,
called the open mapping theorem, which implies that a bounded linear
mapping $T$ of $V$ onto $V$ with trivial kernel has bounded inverse.
In the proof, one first uses the Baire category theorem to show that
the image of the closed unit ball in $V$ under $T$ contains a
neighborhood of the origin, and then one shows that in fact the image
of the open unit ball in $V$ under $T$ contains a neighborhood of the
origin.  The first step does not give quantitative information, and
for that matter the hypotheses are not very quantitative anyway.

	One should not necessarily take this result too seriously,
although it is interesting that the continuity of the inverse is
connected to surjectivity in principle in this manner.  In practice,
one often derives an estimate anyway.  It is not so easy to show that
a mapping is surjective, and a common way to do this is to establish
the appropriate bound and show that the image is dense.

	Another general result, called the closed graph theorem, says
that a linear mapping $T$ from the Banach space $V$ to itself is
continuous if its graph, as a linear subspace of $V \times V$, is
closed.  To be more explicit, the graph of $T$ is closed if for every
sequence $\{v_j\}_{j=1}^\infty$ in $V$ which converges to $0$ and has
the property that $\{T(v_j)\}_{j=1}^\infty$ converges in $V$, we have
that $\lim_{j \to \infty} T(v_j) = 0$.  One can derive the closed
graph theorem from the open mapping theorem, because the mapping $(v,
T(v)) \mapsto v$ defines a continuous linear mapping from the graph of
$T$ onto $V$, and the assumption that the graph of $T$ is closed
implies that it is a Banach space, so that the inverse mapping $v
\mapsto (v, T(v))$ from $V$ to the graph of $T$ is then continuous
too.

	Again, one should not necessarily take this result too
seriously, although it is interesting again that in principle the
continuity of $T$ is implied by apparently quite mild conditions.  Of
course having $T$ be defined everywhere on the Banach space, and in
some reasonable manner, is already a pretty strong condition.  As
before, in practice this is often because one derives an estimate
anyway, on a dense subspace of $V$, and then extends the operator to
all of $V$ by continuity.

	Let us specialize for the moment to the case where $V$ is a
\emph{Hilbert space}.\index{Hilbert spaces} This means that there is
an \emph{inner product}\index{inner products} $\langle v, w \rangle$
on $V$ which determines the norm on $V$, and that $V$ is complete with
respect to this norm.  More precisely, the inner product $\langle v, w
\rangle$ is a function from $V \times V$ to the real or complex
numbers, as appropriate, such that
\begin{equation}
	v \mapsto \langle v, w \rangle
\end{equation}
is a linear functional on $V$ for each $v \in V$, the symmetry condition
\begin{equation}
	\langle w, v \rangle = \langle v, w \rangle
\end{equation}
in the real case and
\begin{equation}
	\langle w, v \rangle = \overline{\langle v, w \rangle}
\end{equation}
in the complex case holds for all $v, w \in V$, and
\begin{equation}
	\langle v, v \rangle
\end{equation}
is a nonnegative real number for all $v \in V$ which is equal to $0$
if and only if $v = 0$.

	When one has such an inner product $\langle v, w \rangle$,
one can define the associated norm on $V$ by
\begin{equation}
	\|v\| = \langle v, v \rangle^{1/2}.
\end{equation}
As is well-known, one has the Cauchy--Schwarz
inequality\index{Cauchy--Schwarz inequality}
\begin{equation}
	|\langle v, w \rangle| \le \|v\| \, \|w\|
\end{equation}
for all $v, w \in V$.  This inequality can be used to verify the
triangle inequality for $\|v\|$, so that $\|v\|$ is indeed a norm.

	Two elements $v_1$, $v_2$ of $V$ are said to be
\emph{orthogonal}\index{orthogonal elements in a Hilbert space} if
\begin{equation}
	\langle v_1, v_2 \rangle = 0,
\end{equation}
in which case we write $v_1 \perp v_2$.  If $v_1 \perp v_2$, then
\begin{equation}
	\|v_1 + v_2\|^2 = \|v_1\|^2 + \|v_2\|^2.
\end{equation}
For any two elements $v_1$, $v_2$ of $V$, whether or not they are
orthogonal, we have the \emph{parallelogram
identity}\index{parallelogram identity}
\begin{equation}
	\|v_1 + v_2\|^2 + \|v_1 - v_2\|^2 = 2 \, \|v_1\|^2 + 2 \, \|v_2\|^2.
\end{equation}

	Let $W$ be a linear subspace of $V$.  The \emph{orthogonal
complement}\index{orthogonal complement of a linear subspace in a
Hilbert space} of $W$ in $V$ is denoted $W^\perp$ and defined to be
the set of $z \in V$ such that $z \perp w$ for all $w \in W$.  It is
easy to see that $W^\perp$ is automatically a closed subspace of $W$,
and that
\begin{equation}
	W \cap W^\perp = \{0\}.
\end{equation}

	Suppose that $W$ is a linear subspace of $V$, that $v$ is an
element of $V$, and that $w$ is an element of $W$ such that $v - w \in
W^\perp$.  Such an element $w$ of $W$ is uniquely determined by $v$,
since if $w' \in W$ also satisfies $v - w' \in W^\perp$, then $w - w'
\in W$ and $w - w' = (w - v) + (v - w') \in W^\perp$.  Also, if $u$ is
any element of $W$, then
\begin{eqnarray}
\lefteqn{\|v - u\|^2 = \langle v - u, v - u \rangle}		\\
	& & = \langle v - w, v - w \rangle + \langle w - u, w - u \rangle
		= \|v - w\|^2 + \|w - u\|^2,		\nonumber
\end{eqnarray}
and hence $w$ is the unique element of $W$ whose distance to $v$ is as
small as possible.

	Conversely, suppose that $W$ is a closed linear subspace of $V$
and that $v$ is an element of $V$, and let us show that there is a
$w \in W$ such that $\langle v - w, u \rangle = 0$ for all $u \in W$.
We do this by looking for an element of $W$ whose distance to $v$
is as small as possible.  If
\begin{equation}
	\dist(v, W) = \inf \{\|v - z\| : z \in W\}
\end{equation}
is equal to $0$, then $v \in W$, because $W$ is assumed to be closed,
and we can simply take $w = v$.

	For each positive integer $j$, choose $w_j \in W$ so that
\begin{equation}
	\|v - w_j\| \le \dist(v, W) + \frac{1}{j}.
\end{equation}
If $j$ and $l$ are positive integers, then the parallelogram identity
applied to $v - w_j$, $v - w_l$ yields
\begin{equation}
	\|2 v - w_j - w_l\|^2 + \|w_j - w_l\|^2 
		= 2 \, \|v - w_j\|^2 + 2 \, \|v - w_l\|^2.
\end{equation}
Because
\begin{equation}
	\|2 v - w_j - w_l\| = 2 \biggl\|v - \frac{w_j + w_l}{2} \biggr\| 
							\ge 2 \dist(v, W),
\end{equation}
we get that
\begin{equation}
	\|w_j - w_l\|^2 
		\le 4 \dist(v, W) \biggl(\frac{1}{j} + \frac{1}{l}\biggr) 
			+ \frac{2}{j^2} + \frac{2}{l^2}.
\end{equation}

	This shows that $\{w_j\}_{j=1}^\infty$ is a Cauchy sequence in
$V$.  The assumption that $V$ is complete implies that this sequence
converges to a point $w \in W$, since $W$ is supposed to be a closed
subspace of $V$.  Thus
\begin{equation}
	\|v - w\| = \dist(v, W),
\end{equation}
and it is not difficult to derive from this that $v - w \in W^\perp$.

	Let us write $P_W(v)$ for the unique element $w$ of $W$ such
that $v - w \in W^\perp$.  It is easy to see that $P_W$ is a linear
mapping from $V$ to $W$ such that $P_W(w) = w$ when $w \in W$, $P_W(u)
= 0$ when $u \in W^\perp$, and
\begin{equation}
	\|v\|^2 = \|P_W(v)\|^2 + \|v - P_W(v)\|^2.
\end{equation}
Thus $P_W$ is a bounded linear operator on $V$ with norm equal to $1$,
except in the trivial case where $W = \{0\}$, when $P_W$ is the zero
operator.

	This operator $P_W$ is called the \emph{orthogonal
projection}\index{orthogonal projection onto a subspace of Hilbert
space} of $V$ onto $W$.  Notice that for all $v_1, v_2 \in V$ we have
that
\begin{equation}
	\langle P_W(v_1), v_2 \rangle = \langle P_W(v_1), P_W(v_2) \rangle
			= \langle v_1, P_W(v_2) \rangle.
\end{equation}
Also, if $W$ is a closed linear subspace of $V$ which is not equal to
$V$, it follows that there are nonzero elements of $V$ in $W^\perp$.

	In general a bounded linear operator $A$ on $V$ is said to be
\emph{self-adjoint}\index{self-adjoint bounded linear operators on a
Hilbert space} if
\begin{equation}
	\langle A(v), w \rangle = \langle v, A(w) \rangle
\end{equation}
for all $v, w \in V$.  Observe that the sum of two self-adjoint
bounded linear operators on a Hilbert space is also self-adjoint, and
a real number times a self-adjoint operator is again self-adjoint.
The restriction to real multiples is important when $V$ is a complex
Hilbert space.

	Suppose that $A$ is a self-adjoint bounded linear operator
on $V$.  If $v$ is any vector in $V$, then $A(v) = 0$ if and only if
\begin{equation}
	\langle A(v), w \rangle = 0
\end{equation}
for every $w \in V$, and this holds if and only if
\begin{equation}
	\langle v, A(w) \rangle = 0,
\end{equation}
which is to say that $v$ is orthogonal to the image of $A$.  Thus
$A$ has trivial kernel if and only if the image of $A$ is dense in $V$.

	Thus, a self-adjoint bounded linear operator $A$ on $V$ is
invertible if and only if there is a positive real number $c$
such that
\begin{equation}
	\|A(v)\| \ge c \, \|v\|
\end{equation}
for all $v \in V$.  In particular a self-adjoint bounded linear
operator $A$ is invertible if it satisfies the positivity condition
that there is a positive real number $\alpha$ such that
\begin{equation}
	\langle A(v), v \rangle \ge \alpha \, \|v\|^2
\end{equation}
for all $v \in V$.  Let us note that if $B$ is a self-adjoint bounded
linear operator on $V$ which is invertible, then $A = B^2$ satisfies
the kind of positivity condition just mentioned.

	Now let us return to the setting of a general Banach space
$V$.  We would like to consider operators of the form $T + A$, where
$T$ is a bounded linear operator on $V$ which is invertible and $A$ is
a bounded linear operator which is in the closure of the space of
finite-rank bounded linear operators on $V$, with respect to the
operator norm.  Recall that a linear mapping between vector spaces has
\emph{finite rank}\index{finite rank linear mappings between vector
spaces} if its image is finite-dimensional.

	Here is a basic situation.  Let $V$ be the space of real or
complex-valued continuous functions on the unit interval, and suppose
that $a(x, y)$ is a real or complex-valued continuous function on the
unit square $[0, 1] \times [0, 1]$, as appropriate.  One can show that
the integral operator on $V$ defined by
\begin{equation}
	A(f)(x) = \int_0^1 a(x, y) \, f(y) \, dy
\end{equation}
is a bounded linear operator which can be approximated in the operator
norm by bounded finite-rank linear operators on $V$, using the fact
that $a(x, y)$ is uniformly continuous on $[0, 1] \times [0, 1]$.

	If $T$ is a bounded linear operator on $V$ which is invertible
and $A$ can be approximated in the operator norm by finite rank operators,
then we can write $T + A$ as $T_1 + A_1$, where $T_1$ is also a bounded
linear operator on $V$ which is invertible, and $A_1$ is a bounded linear
operator on $V$ with finite rank.  That is, one can subtract a small part 
from $A$ and add it to $T$ to rewrite $T + A$ in this way.  A basic fact
about operators of this form is that they are invertible exactly when
their kernels are trivial.

	Assuming that $V$ is a complex Banach space, the spectrum of a
bounded linear operator $R$ is defined to be the set of complex
numbers $\lambda$ such that $\lambda \, I - R$ is not invertible.  If
$V$ is a complex Hilbert space and $R$ is self-adjoint, then one can
check that the spectrum of $R$ consists of real numbers.  On a complex
Banach space in general, if $R$ can be approximated in the operator
norm by finite rank operators, then a nonzero complex number $\lambda$
lies in the spectrum of $R$ if and only if $\lambda$ is an eigenvalue
of $R$, which is to say that there is a nonzero vector $v \in V$ such
$R(v) = \lambda \, v$.

\backmatter

\newpage

\addcontentsline{toc}{chapter}{Index}
\printindex

\end{document}